\documentclass[final,leqno,subequations,onefignum,onetabnum]{siamltex}

\usepackage{amsmath,amssymb}
\usepackage[sc]{mathpazo}
\usepackage{graphicx}
\usepackage{algorithm,algorithmic}
\usepackage{multirow}
\usepackage{makecell}
\usepackage{textcomp}
\usepackage{rotating}
\usepackage{siunitx}
\usepackage{paralist}
\usepackage[top=1in,bottom=1.5in,left=1in,right=1in]{geometry}
\usepackage{hhline}

\usepackage[table]{xcolor}

\usepackage[pagebackref,colorlinks=true,pdftex,citecolor=blue,urlcolor=blue]{hyperref}

\graphicspath{{./figures/}}

\newcommand{\figref}[1]{Fig.~\ref{#1}}
\newcommand{\tabref}[1]{Tab.~\ref{#1}}
\newcommand{\secref}[1]{\S\ref{#1}}
\newcommand{\algref}[1]{Alg.~\ref{#1}}

\newcommand{\real}{\ensuremath{\operatorname{Re}}}
\newcommand{\imag}{\ensuremath{\operatorname{Im}}}
\newcommand{\tabadjust}{\centering\scriptsize\renewcommand\arraystretch{1.2}}

\definecolor{utorange}   {RGB} {203,96,21}
\definecolor{utblack}    {RGB} {99,102,106}
\definecolor{utbrown}    {RGB} {110,98,89}
\definecolor{utsecbrown} {RGB} {217,200,158}
\definecolor{utsecgreen} {RGB} {208,222,187}
\definecolor{utsecblue}  {RGB} {127,169,174}
\definecolor{dgreen}     {RGB} {0,100,0}

\def\rowcolA{\rowcolor{black!40}}
\def\rowcolB{\rowcolor{black!30}}
\def\rowcolC{\rowcolor{black!20}}
\def\rowcolD{\rowcolor{black!10}}

\def\cellcolA{\cellcolor{black!40}}

\newcommand{\F}[1]{\ensuremath{\mathcal{#1}}}    
\newcommand{\D}[1]{\ensuremath{\mathcal{#1}}}    

\newcommand{\ns}[1]{\ensuremath{\mathbf{#1}}}
\newcommand{\fs}[1]{\ensuremath{\mathcal{#1}}}

\newcommand{\Div}{\ensuremath{\nabla\cdot}}
\newcommand{\Grad}{\ensuremath{\nabla}}
\newcommand{\Lap}{\rotatebox[origin=c]{180}{$\nabla$}}

\newcommand{\idiv}{\Div}
\newcommand{\igrad}{{\Grad}}
\newcommand{\ilap}{\Lap}

\newcommand{\half}[1]{\frac{#1}{2}}
\renewcommand{\d}[1]{\,\mathrm{d}#1}
\newcommand{\dt}{\d{t}}
\newcommand{\dx}{\d{\vect{x}}}

\newcommand{\bigO}{\mathcal{O}}
\newcommand{\p} {\partial}

\newcommand{\secheadskip}{${}$}

\newcommand{\vect}[1]{\boldsymbol{#1}} 
\newcommand{\mat}[1]{\boldsymbol{#1}}  

\newcommand{\bipa}{\begin{inparaenum}[(\itshape i\upshape)]}
\newcommand{\eipa}{\end{inparaenum}}

\newcommand{\bipasub}{\begin{inparaenum}[(\itshape a\upshape)]}
\newcommand{\eipasub}{\end{inparaenum}}

\newcommand{\acr}[1]{#1}
\newcommand{\ipoint}[1]{\textit{#1}:}
\newcommand{\iname}[1]{\emph{#1}}

\newcommand{\iquote}[1]{{\it #1}}

\newcommand{\ccomment}[1]{{\hfill\color{gray}$\rhd$\,#1}}

\newcommand{\bea} {\begin{eqnarray}}
\newcommand{\eea} {\end{eqnarray}}
\newcommand{\barr} {\begin{array}}
\newcommand{\earr} {\end{array}}
\newcommand{\bean} {\begin{eqnarray*}}
\newcommand{\eean} {\end{eqnarray*}}

\newcommand{\iom}[1]{\int_{\Omega}#1\dx}

\newcommand{\iut}[1]{\int_0^1#1\dt}

\newcommand{\ipl}[2]{\langle #1,#2 \rangle}

\newcommand{\defeq}{\ensuremath{\mathrel{\mathop:\!\!=}}}
\newcommand{\eqdef}{\ensuremath{\mathrel{=\!\!\mathop:}}}

\newcommand{\T}{\ensuremath{\mathsf{T}}}
\newcommand{\adj}{\ensuremath{\mathsf{H}}}

\definecolor{sred}{cmyk}{0.01,0.98,0,0.2} 
\reversemarginpar
\newcommand{\mmargin}[1]{{\marginpar{\em\tiny #1}}}\renewcommand{\mmargin}[1]{}

\title{An inexact Newton-Krylov algorithm for constrained diffeomorphic image registration}

\author
{
Andreas Mang\thanks
{
The Institute for Computational Engineering and Sciences,
The University of Texas at Austin,
Austin, Texas, 78712-0027, US
({\tt andreas@ices.utexas.edu}, {\tt gbiros@acm.org})
}
\and George Biros\footnotemark[1]
}

\begin{document}

\maketitle

\begin{abstract}
We propose numerical algorithms for solving large deformation diffeomorphic image registration problems. We formulate the nonrigid image registration problem as a problem of optimal control. This leads to an infinite-dimensional partial differential equation (PDE) constrained optimization problem.

The PDE constraint consists, in its simplest form, of a hyperbolic transport equation for the evolution of the image intensity. The control variable is the velocity field. Tikhonov regularization on the control ensures well-posedness. We consider standard smoothness regularization based on $H^1$- or $H^2$-seminorms. We augment this regularization scheme with a constraint on the divergence of the velocity field (control variable) rendering the deformation incompressible (Stokes regularization scheme) and thus ensuring that the determinant of the deformation gradient is equal to one, up to the numerical error.

We use a Fourier pseudospectral discretization in space and a Chebyshev pseudospectral discretization in time. The latter allows us to reduce the number of unknowns and enables the time-adaptive inversion for nonstationary velocity fields. We use a preconditioned, globalized, matrix-free, inexact Newton-Krylov method for numerical optimization. A parameter continuation is designed to estimate an optimal regularization parameter. Regularity is ensured by controlling the geometric properties of the deformation field. Overall, we arrive at a black-box solver that exploits computational tools that are precisely tailored for solving the optimality system. We study spectral properties of the Hessian, grid convergence, numerical accuracy, computational efficiency, and deformation regularity of our scheme. We compare the designed Newton-Krylov methods with a globalized Picard method (preconditioned gradient descent). We study the influence of a varying number of unknowns in time.

The reported results demonstrate excellent numerical accuracy, guaranteed local deformation regularity, and computational efficiency with an optional control on local mass conservation. The Newton-Krylov methods clearly outperform the Picard method if high accuracy of the inversion is required. Our method provides equally good results for stationary and nonstationary velocity fields for two-image registration problems.
\end{abstract}

\newcommand{\slugmaster}{\slugger{siims}{xxxx}{xx}{x}{x--x}}



\pagestyle{myheadings}
\thispagestyle{plain}
\markboth
{
ANDREAS MANG AND GEORGE BIROS
}
{
A NEWTON-KRYLOV ALGORITHM FOR CONSTRAINED IMAGE REGISTRATION
}

\section{Introduction and Motivation}
\label{s:intro}

Image registration has become a key area of research in computer vision and (medical) image analysis~\cite{Modersitzki:2004a, Sotiras:2013a}. The task is to establish spatial correspondence between two images $m_R : \bar{\Omega} \rightarrow \ns{R}$, $\vect{x}\mapsto m_R(\vect{x})$, and $m_T : \bar{\Omega} \rightarrow \ns{R}$, $\vect{x}\mapsto m_R(\vect{x})$, with compact support on some domain $\Omega\defeq(-\pi,\pi)^d\subset\ns{R}^d$ via a mapping $\vect{y} : \ns{R}^d \rightarrow \ns{R}^d$, $\vect{x}\mapsto\vect{y}(\vect{x})$, such that $m_T \circ \vect{y} \approx m_R$~\cite{Droske:2003a, Fischer:2008a}, where $\circ$ is the function composition. Here, $m_T$ is referred to as the \iname{template image} (the image to be registered), $m_R$ is referred to as the \iname{reference image} (the fixed image), and $d\in\{2,3\}$ is the data dimensionality. We limit ourselves to nonrigid image registration. The search for $\vect{y}$ is typically formulated as a variational optimization problem~\cite{Fischer:2008a, Modersitzki:2004a},
\begin{equation}\label{e:opt-prob-y}
\min_{\vect{y}\in\fs{Y}}
\left\{\F{J}[\vect{y}]
\defeq
\half{1} \|m_R -m_T\circ\vect{y}\|_{L^2(\Omega)}^2
+ \half{\beta}\F{S}[\vect{y}]\right\}.
\end{equation}

\noindent The proximity between $m_R$ and $m_T\circ \vect{y}$ is measured on the basis of an $L^2$-distance (other measures can be considered~\cite{Modersitzki:2004a, Sotiras:2013a}). The functional $\F{S}$ in~\eqref{e:opt-prob-y} is a regularization model that is introduced to overcome ill-posedness. The regularization parameter $\beta>0$ weights the contribution of $\F{S}$. Various regularization models $\F{S}$ have been considered (see~\cite{Broit:1981a, Burger:2013a, Christensen:1996a, Droske:2003a, Fischer:2002a, Fischer:2003a, FrohnSchauf:2008a, Henn:2005a, Rumpf:2009a} for examples).

A key requirement in many image registration problems is that the mapping $\vect{y}$ is a \iname{diffeomorphism}~\cite{Beg:2005a, Burger:2013a, Dupuis:1998a, Trouve:1998a, Vercauteren:2008a, Vercauteren:2009a}. This translates into the necessary condition $\det(\igrad\vect{y})>0$, where $\igrad\vect{y}\in\ns{R}^{d\times d}$ is the \iname{deformation gradient} (also referred to as the \iname{Jacobian matrix}). An intuitive approach is to explicitly safeguard against nondiffeomorphic mappings $\vect{y}$ by adding a constraint to~\eqref{e:opt-prob-y}~\cite{Haber:2007a, Lee:2010a, Rohlfing:2003a}. Another strategy is to perform the optimization on the manifold of diffeomorphic mappings~\cite{Ashburner:2007a, Ashburner:2011a, Beg:2005a, Hernandez:2009a, Mansi:2011a, Miller:2004a, Vercauteren:2008a, Vercauteren:2009a}. The latter models, in general, do not control geometric properties of the deformation field and may result in fields that are close to being nondiffeomorphic. Further, for certain image registration problems, restricting the search space to the manifold of diffeomorphisms does not necessarily guarantee that $\vect{y}$ is \emph{physically meaningful}. Some applications may benefit from extending these types of models by introducing additional constraints. One example for such a constraint is incompressibility (i.e.\ enforcing $\det(\igrad\vect{y})=1$; see also~\cite{Borzi:2002a, Chen:2011a, Mansi:2011a, Ruhnau:2007a}). Incompressibility is a requirement that might be of interest in medical image computing applications. If required, we can modify the incompressibility constraint to control the deviation of $\det(\igrad\vect{y})$ from identity. This will be the topic of a follow-up paper, in which we will extend our formulation. Here, we focus on the algorithmic issues of incorporating the incompressibility constraint. Furthermore, we remark that our optimal control formulation can naturally be extended to account for more complex constraints on the velocity field, for example, constraints related to biophysical models (examples of such models can be found in~\cite{Gooya:2012a, Han:2014a, Hogea:2008a, Mang:2012b, Sundar:2009a}).

In what follows, we outline our method (see \secref{s:outline}), summarize the key contributions (see~\secref{s:contributions}), list the limitations of our method (see~\secref{s:limitations}), and relate our work to existing approaches (see~\secref{s:related-work}).

\subsection{Outline of the Method}
\label{s:outline}

We assume the images are compactly supported functions, defined on the open set $\Omega\defeq(-\pi,\pi)^d \subset \ns{R}^d$ with boundary $\p\Omega$ and closure $\bar{\Omega} \defeq \Omega \cup \p\Omega$. The deformation is modeled in an Eulerian frame of reference. We introduce a pseudotime variable $t > 0$ and solve on
a unit time horizon for a velocity field $\vect{v} : \bar{\Omega} \times [0,1] \rightarrow \ns{R}^d$, $(\vect{x},t) \mapsto \vect{v}(\vect{x},t)$, as follows:
\begin{subequations}\label{e:opt-prob-v}
\begin{equation}\label{e:objective}
\min_{\vect{v}\in\fs{V}} \left\{\F{J}[\vect{v}]
\defeq \half{1}\|m_R - m_1\|_{L^2(\Omega)}^2
+ \int_0^1 \F{S}[\vect{v}]\,\d{t}
\quad \text{subject\,to}\quad \D{C}[m,\vect{v}] = 0\right\},
\end{equation}
\noindent where
\begin{equation}\label{e:constraint}
\D{C}[m,\vect{v}] \defeq \left\{
\begin{array}{c}
\p_t m + \igrad m \cdot \vect{v} \\
m - m_T                         \\
\gamma(\idiv \vect{v})
\end{array}\right.
\begin{array}{l}
{\rm in} \; \Omega \times (0,1], \\
{\rm in} \; \Omega \times \{0\}, \\
{\rm in} \; \Omega \times [0,1],
\end{array}
\end{equation}
\end{subequations}

\noindent with periodic boundary conditions on $\p\Omega$. The parameter $\gamma \in \{0,1\}$ in~\eqref{e:constraint} is introduced to enable or disable this constraint on the control $\vect{v}$. In PDE constrained optimization theory, $m$ is referred to as the \iname{state variable} and $\vect{v}$ as the \iname{control variable}. The first equation in~\eqref{e:constraint}, in combination with its initial condition (second equation), models the flow of $m_T$ subject to $\vect{v}$, where $m:\bar{\Omega}\times[0,1]\rightarrow\ns{R}$, $(\vect{x},t)\mapsto m(\vect{x},t)$, represents the transported intensities of $m_T$. Accordingly, the \iname{final} \iname{state} $m_1 \defeq m(\cdot,1)$, $m_1: \bar{\Omega} \rightarrow \ns{R}$, $\vect{x} \mapsto m_1(\vect{x})$, corresponds to $m_T \circ \vect{y}$ in~\eqref{e:opt-prob-y}. We measure the proximity between the deformed template image $m_1$ and the reference image $m_R$ in terms of an $L^2$-distance. Once we have found $\vect{v}$, we can compute $\vect{y}$ from $\vect{v}$ as a post-processing step (this is also true for the deformation gradient $\igrad\vect{y}$; see~\secref{s:map-and-def-grad} for details). The third equation in~\eqref{e:constraint} is a control on the divergence of $\vect{v}$ and guarantees that the flow is incompressible (Stokes flow), i.e.\ the volume is conserved. This is equivalent to enforcing $\det(\igrad\vect{y}) = 1$ (see \cite[p.~77ff.]{Gurtin:1981a}).

We use either an $H^1$- or an $H^2$-seminorm for the smoothness regularization $\F{S}$ (resulting in Laplacian or biharmonic vector operators, respectively; see~\secref{s:regularization-model}). We report results for standard $H^1$- and $H^2$-regularization (neglecting $\idiv \vect{v} = 0$, i.e.\ $\gamma = 0$ in~\eqref{e:constraint}) and for a Stokes regularization scheme (incompressible flow; enforcing $\idiv \vect{v} = 0$, i.e.\ $\gamma = 1$ in~\eqref{e:constraint}).

In~\secref{s:numerics} we will see that the optimality condition for~\eqref{e:opt-prob-v} is a  system of space-time nonlinear multicomponent PDEs for the transported image $m$, the velocity $\vect{v}$ and the adjoint variables for the transport and the divergence condition. Efficiently solving this system is a significant challenge.  For example, when we include the incompressibility constraint, the equation for the velocity ends up being a linear Stokes equation.

We solve for the first-order optimality conditions using a globalized, matrix-free, preconditioned Newton-Krylov method for the Schur complement of the velocity $\vect{v}$ (a linearized Stokes problem driven by the image mismatch). We first derive the optimality conditions and then discretize using a pseudospectral discretization in space with a Fourier basis. We use a second-order Runge-Kutta method in time. The preconditioner for the Newton-Krylov schemes is based on the exact spectral inverse of the second variation of $\F{S}$.

\subsection{Contributions}
\label{s:contributions}

The fundamental contributions are:
\begin{itemize}
\item We design a numerical scheme for \eqref{e:opt-prob-v} with the following key features:
\begin{itemize}
    \item We use an adjoint-based Newton-type solver.
    \item We provide a fast Stokes solver.
    \item We introduce a spectral Galerkin method in time.
    \item We design a parameter continuation method for automatically selecting the regularization parameter $\beta$.
    \item Our framework guarantees deformation regularity\footnote{Note that controlling the magnitude of $\det(\igrad \vect{y})$ is not sufficient to guarantee that $\vect{y}$ is locally well behaved. Therefore, our framework features geometric constraints that guarantee a nice, locally diffeomorphic mapping $\vect{y}$. In particular, we control the shear angle of the cells within the deformed grid during the parameter continuation in $\beta$.}.
\end{itemize}
\item We provide a numerical study for the designed framework. We compare a globalized Picard method (preconditioned gradient descent) to an inexact Newton-Krylov and a Gauss-Newton-Krylov method. We report results for synthetic and real-world problems. We study spectral properties, grid convergence and numerical accuracy of the proposed scheme. We study the effects of compressible (plain $H^1$- and $H^2$ regularization) and incompressible (Stokes regularization) deformation models. We report results for a varying number of unknowns in time (i.e.\ inverting for stationary and nonstationary velocity fields).
\end{itemize}

The numerical discretization (pseudospectral) allows for an efficient solution of the Stokes-like equations by eliminating the pressure (i.e.\ the adjoint variable for the incompressibility constraint $\idiv \vect{v} = 0$ in~\eqref{e:constraint}).

The $\inf$-$\sup$ condition for pressure spaces~\cite[p.~200ff.]{Brezzi:1991a} is not an issue with our scheme\footnote{The $\inf$-$\sup$ condition is an important requirement when solving Stokes-like problems via the finite element method (see~\cite[p.~200ff.]{Brezzi:1991a}; examples for a finite element discretization of similar problem formulations can be found in~\cite{Chen:2011a, Ruhnau:2007a}). Satisfying the $\inf$-$\sup$ condition ensures that the finite element solution exists, is stable and optimal. Essentially, we require two different finite element spaces for the discretization of the pressure and the velocity. The $\inf$-$\sup$ condition is key for the decision on an adequate pair of spaces. For our scheme, we can use the same basis (Fourier) for the discretization of the pressure and the velocity. Also, it is very efficient since we can eliminate the pressure from the optimality system (see \secref{s:elimination}).}. In fact, for smooth images, our scheme is spectrally accurate in space and second-order accurate in time. We will see that we can numerically enforce incompressibility up to almost machine precision. Also, our scheme allows for efficient preconditioning of the Hessian: at the cost of a diagonal scaling we obtain a problem with a bounded condition number.

Overall, we demonstrate that the designed framework\bipa\item is efficient and accurate, \item features a precise control of the deformation regularity, and \item does not require manual tuning of parameters.\eipa

\subsection{Limitations}
\label{s:limitations}

The main limitations of our method are:
\begin{itemize}
    \item The considered model assumes a constant illumination of $m_R$ and $m_T$ (a consequence of the transport equation and the $L^2$-distance in \eqref{e:opt-prob-v}). Therefore, it is (in its current state) not directly applicable to multi-modal registration problems. Nevertheless, let us remark that the $L^2$-distance is commonly used in practice~\cite{Beg:2005a, Chen:2011a, FrohnSchauf:2008a, Hart:2009a, Hernandez:2009a, Lee:2010a, Museyko:2009a, Vialard:2012a}.
    \item The efficient use of a Fourier discretization for the PDEs requires periodic boundary conditions. If the images are not periodic, we artificially introduce periodic boundary conditions by mollification and zero padding.
\end{itemize}

\subsection{Related Work}
\label{s:related-work}

Due to the vast body of literature, it is not possible to provide a comprehensive review of numerical methods for nonrigid image registration. Background on image registration formulations and numerics can be found in~\cite{Fischer:2008a, Modersitzki:2004a, Sotiras:2013a}.  We limit the discussion to approaches that\bipa\item model the deformation via a velocity field $\vect{v}$, \item view image registration as a problem of optimal control and/or \item constrain $\vect{v}$ to be divergence-free (i.e.\ introduce a mass conservation equation as an additional constraint).\eipa

Fluid mechanical models have been introduced~\cite{Christensen:1994a, Christensen:1996a} to overcome limitations of small deformation models~\cite{Broit:1981a, Fischer:2002a, Fischer:2003a}. The work in~\cite{Christensen:1994a, Christensen:1996a} has been extended in~\cite{Beg:2005a, Dupuis:1998a, Miller:2004a, Trouve:1998a} using concepts from differential geometry. This class of approaches is referred to as \iname{large deformation diffeomorphic metric mapping} (\acr{LDDMM}). Under the assumption that $\vect{v}$ is adequately smooth, it is guaranteed that $\vect{y}$ is a diffeomorphism~\cite{Dupuis:1998a, Trouve:1998a}. The associated smoothness requirements are enforced by the regularization model $\F{S}$ (typically an $H^2$-norm)~\cite{Ashburner:2011a, Beg:2005a}. The optimization is performed on the space of nonstationary velocity fields~\cite{Beg:2005a}. To reduce the number of unknowns, it has been suggested to perform the optimization either on the space of stationary velocity fields~\cite{Arsingny:2006a, Ashburner:2007a, Hernandez:2009a} or with respect to an initial momentum that entirely defines the flow of the map $\vect{y}$~\cite{Ashburner:2011a, Vialard:2012a, Younes:2007a}.

The idea of parameterizing a diffeomorphism $\vect{y}$ via a stationary velocity field~\cite{Arsingny:2006a} has also been introduced to the \iquote{demons} registration framework~\cite{Mansi:2011a, Vercauteren:2008a, Vercauteren:2009a}. Here, optimization is performed in a sequential fashion, alternating between updates resulting from the distance measure (forcing term) and the application of a smoothing operator to regularize the problem (typically through Gaussian smoothing~\cite{Thirion:1998a, Vercauteren:2008a, Vercauteren:2009a}). This scheme is somewhat equivalent to the Picard scheme we discuss in our paper but it is unclear how one couples it with line search or trust region techniques.

Approaches that more closely reflect an optimal control PDE constrained optimization formulation~\eqref{e:opt-prob-v} are described in~\cite{Borzi:2002a, Chen:2011a, Hart:2009a, Lee:2010a, Lee:2011a, Museyko:2009a, Ruhnau:2007a, Vialard:2012a}. The model in~\cite{Chen:2011a} is equivalent to~\eqref{e:opt-prob-v}. The model in~\cite{Ruhnau:2007a} follows the traditional optical flow formulation~\cite{Horn:1981a}. The conceptual difference between our formulation and the latter is that in optical flow, the transport equation constraint appears in the objective (see e.g.~\cite{Horn:1981a, Kalmoun:2011a, Ruhnau:2007a}) and is therefore only fulfilled approximately in an $L^2$ least squares sense. We treat it as a hard constraint instead. In~\cite{Museyko:2009a} an optimal mass transport formulation is described that is based on the Monge-Kantorovich problem. The formulations in~\cite{Hart:2009a, Vialard:2012a} do not account for incompressibility. The optical flow approach in~\cite{Borzi:2002a} treats incompressibility as an $L^2$-penalty. An optimal control formulation for a constant-in-time velocity field was proposed in~\cite{Lee:2010a, Lee:2011a}, in which the divergence of the velocity field is penalized along with smoothness constraints.

What sets our work apart are the numerical algorithms and the discretization scheme. Almost all existing efforts on large deformation diffeomorphic image registration that are closely related to our optimal control formulation exclusively use first-order information for numerical optimization~\cite{Ashburner:2007a, Beg:2005a, Borzi:2002a, Chen:2011a, Christensen:1994a,  Christensen:1996a, Hart:2009a, Hernandez:2009a, Lee:2010a, Lee:2011a, Museyko:2009a, Ruhnau:2007a, Vialard:2012a, Younes:2007a}. We use second-order information. The only work\footnote{After the submission of our work, another contribution on second-order numerical optimization for LDDMM appeared~\cite{Hernandez:2014a}.} in the context of large deformation diffeomorphic image registration that to our knowledge uses second-order information is~\cite{Ashburner:2011a}. The model in~\cite{Ashburner:2011a} is based on the LDDMM framework~\cite{Beg:2005a, Dupuis:1998a, Lee:2010a, Lee:2011a, Miller:2004a, Hernandez:2009a, Trouve:1998a}. The inversion is, likewise to~\cite{Vialard:2012a}, performed with respect to an initial momentum. No additional constraints on $\vect{v}$ are considered. Another difference is that we use a Galerkin method in time to reduce the number of unknowns. This allows us to invert for stationary~\cite{Arsingny:2006a, Ashburner:2007a, Hernandez:2009a, Mansi:2011a, Vercauteren:2008a, Vercauteren:2009a} as well as time-varying~\cite{Beg:2005a, Borzi:2002a, Chen:2011a, Hart:2009a} velocity fields. Nothing changes in our formulation other than the number of unknowns. Furthermore, we globalize our methods with a line search strategy (i.e.\ we guarantee a sufficient decrease of the objective $\F{J}$). This is a standard---yet important---ingredient for guaranteeing convergence, which is often not accounted for~\cite{Ashburner:2007a, Chen:2011a, Chen:2012a, Hart:2009a, Mansi:2011a, Vercauteren:2008a, Vercauteren:2009a}.

We, likewise to~\cite{Borzi:2002a, Chen:2011a, Lee:2010a, Lee:2011a, Mansi:2011a, Ruhnau:2007a}, consider incompressibility as an optional constraint (see \eqref{e:constraint}). Operating with divergence-free velocity fields is equivalent to enforcing $\det(\igrad\vect{y})=1$ up to numerical accuracy (see \cite[p.~77ff.]{Gurtin:1981a}; other formulations for controlling $\det(\igrad\vect{y})$ can be found in~\cite{Burger:2013a, Haber:2010a, Haber:2004a, Haber:2007a, Loeckx:2004a, Modersitzki:2008a, Nielsen:2002a, Pennec:2005a, Rohlfing:2003a, Sdika:2008a, Yanovsky:2007a}). Unlike \cite{Borzi:2002a, Lee:2010a, Lee:2011a}, which penalize the divergence of the velocity, we treat it exactly. We are not arguing that this approach is better per se. The use of penalties is adequate unless one has reasons to insist on an incompressible velocity field. In that case, a penalty method results in ill-conditioning.

Finally our pseudospectral formulation in space allows us to resolve several numerical difficulties related to the incompressibility constraint. For example, the $\inf$-$\sup$ condition for pressure spaces is not an issue with our scheme. Regarding accuracy, for smooth images, our scheme is spectrally accurate in space and second-order accurate in time. We do not have to use different discretization models~\cite{Chen:2011a, Ruhnau:2007a} for solving the individual subsystems of the mixed-type (hyperbolic-elliptic) optimality conditions. Since we use second-order explicit time stepping in combination with Fourier spectral methods, we have at hand a scheme that displays minimal numerical diffusion and does not require flux-limiters~\cite{Borzi:2002a, Chen:2011a, Hart:2009a, Vialard:2012a}. As we will see, the conditioning of the Hessian that appears in our Newton-Krylov scheme can be quite bad. Although the literature for preconditioners for PDE constrained optimization problem is quite rich (e.g.~\cite{Benzi:2011a, Borzi:2002a, Biros:2005b, Biros:2005b, Haber:2001a}), none of these methods directly applies to our formulation. Developing effective preconditioning schemes for our formulation is ongoing work in our group.

\section{Outline}
\label{s:outline-paper}

In~\secref{s:continuous-formulation} the mathematical model is developed. The numerical strategies are described in~\secref{s:numerics}. In particular, we specify\bipa\item the optimality conditions (see~\secref{s:optimality-conditions}), \item strategies for numerical optimization (see~\secref{s:numerical-optimization}) and \item implementation details (see~\secref{s:implementation})\eipa. Numerical experiments on synthetic and real-world data are reported in~\secref{s:numerical-experiments}. Final remarks can be found in~\secref{s:conclusions}.

\section{Continuous Problem Formulation}
\label{s:continuous-formulation}

We provide a summary of the basic notation in~\tabref{t:notation}. The original problem formulation is stated in~\secref{s:outline}. The only missing building block is the considered choices for $\F{S}$ in~\eqref{e:objective}. This is what we discuss next. Note that we neglect any technicalities with respect to the associated function spaces; we assume that the considered functions are adequately smooth (i.e.\ sufficiently many derivatives exist and are bounded).

\begin{table}
\label{t:notation}
\caption
{
Notation (frequently used acronyms and symbols).
}
\tabadjust
\centering
\begin{tabular}{|l|l|}
\hline
\cellcolor[gray]{0.8}
Notation
& \cellcolor[gray]{0.8}
Description
\\
\hline
GN
& Gauss-Newton
\\
KKT system
& Karush-Kuhn-Tucker system
\\
PCG
& preconditioned conjugate gradient method
\\
PDE
& partial differential equation
\\
PDE solve
& solution of hyperbolic PDEs of optimality systems
\eqref{e:first-order-opt-cond} and \eqref{e:second-order-opt-cond}
\\
\hline
$m_R$
& reference (fixed) image
($m_R : \ns{R}^d \rightarrow \ns{R}$)
\\
$m_T$
& template image
(image to be registered; $m_T : \ns{R}^d \rightarrow \ns{R}$)
\\
$\vect{y}$
& mapping (deformation; $\vect{y} : \ns{R}^d \rightarrow \ns{R}^d$)
\\
$\vect{v}$
& velocity field
(control variable; $\vect{v} : \ns{R}^d \times [0,1] \rightarrow \ns{R}^d$)
\\
$m$
& state variable
(transported image; $m : \ns{R}^d \times [0,1] \rightarrow \ns{R}$)
\\
$m_1$
& state variable at $t=1$
(deformed template image; $m_1 : \ns{R}^d \rightarrow \ns{R}$)
\\
$\lambda$
& adjoint variable
(transported mismatch; $\lambda : \ns{R}^d\times[0,1]\rightarrow\ns{R}$)
\\
$\vect{f}$
& body force (drives the registration;
$\vect{f} : \ns{R}^d \times [0,1]\rightarrow \ns{R}^d$)
\\
$\mat{F}_1$
& deformation gradient (tensor field) at $t=1$
($\mat{F}_1 : \ns{R}^d \rightarrow \ns{R}^{d\times d}$)
\\
$\F{J}$
& objective functional
\\
$\F{S}$
& regularization functional
\\
$\D{A}$
& differential operator (first and second variation of $\F{S}$)
\\
$\beta$
& regularization parameter
\\
$\gamma$
& parameter that enables ($\gamma=1$) or disables ($\gamma=0$) the incompressibility constraint
\\
$n_t$
& number of time points (discretization)
\\
$n_c$
& number of coefficient fields (spectral Galerkin method in time)
\\
$\vect{n}_x$
& number of grid points (discretization; $\vect{n}_x = (n_x^1, \ldots, n_x^d)^\T$)
\\
$\vect{g}$
& reduced gradient
(first variation of Lagrangian with respect to $\vect{v}$)
\\
$\D{H}$
& reduced Hessian
(second variation of Lagrangian with respect to $\vect{v}$)
\\
\hline
\end{tabular}
\end{table}

\subsection{Regularization Models}
\label{s:regularization-model}

In contrast to~\cite{Borzi:2002a} we do not explicitly enforce continuity in time. We relax the model to an $L^2$-integrability instead (see~\eqref{e:objective}). This relaxation still yields a velocity field that varies smoothly in time~\cite{Chen:2011a}.

Quadratic smoothing regularization models are commonly used in nonrigid image registration~\cite{Fischer:2008a, Modersitzki:2004a, Modersitzki:2009a, Sotiras:2013a}. They can be defined as
\begin{equation}\label{e:quadratic-regularization}
\F{S}[\vect{v}] \defeq \half{\beta}\|\vect{v}\|^2_{\fs{W}}
= \half{\beta}\langle\D{B}[\vect{v}],\D{B}[\vect{v}]\rangle_{L^2(\Omega)},
\end{equation}

\noindent where $\D{B}$ is a differential operator that (together with its dual) defines the function space $\fs{W}$ and $\beta>0$ is a regularization parameter that balances the contribution of $\F{S}$.

As images are functions of bounded variation, regularity requirements on $\vect{v}\in \fs{V}$, $\fs{V}\defeq L^2([0,1];\fs{W})$ (i.e.\ the choice of $\fs{W}$ in~\eqref{e:quadratic-regularization}) have to be considered with care (for an analytical result see~\cite{Chen:2011a}). Experimental analysis suggests that an $H^1$-seminorm is appropriate if incompressibility is considered (i.e.\ $\gamma=1$ in~\eqref{e:constraint}; see also~\cite{Chen:2011a}). Thus,
\[
\D{B}[\vect{v}] = \igrad \vect{v}.
\]

\noindent This choice is motivated from continuum mechanics and yields a viscous model of \iname{linear Stokes flow} (see~\secref{s:optimality-conditions}; Stokes regularization). If we neglect the incompressibility constraint (i.e.\ $\gamma=0$ in~\eqref{e:constraint}), we use a vectorial Laplacian operator,
\[
\D{B}[\vect{v}] = \ilap \vect{v},
\]
\noindent instead. This choice is motivated by the fact that $H^2$-norm based quadratic regularization is commonly used in large deformation diffeomorphic image registration~\cite{Beg:2005a, Hart:2009a, Hernandez:2009a}.

\section{Numerics}
\label{s:numerics}

We describe the numerical methods used to solve~\eqref{e:opt-prob-v} next. Whenever discretized quantities are considered, a superscript $h$ is added to the continuous variables and operators (i.e.\ the discretized representation of $\vect{v}$ is denoted by $\vect{v}^h$). Likewise, if we refer to a discrete variable at a particular iteration, we will add the iteration index as a subscript (i.e.\ $\vect{v}^h$ at iteration $k$ is denoted by $\vect{v}^h_k$).

We discretize the data on a \iname{nodal grid} in space and time. The number of spatial grid points is denoted by $\vect{n}_x \defeq (n_x^1,\ldots,n_x^d)^\T\in\ns{N}^d$ with spatial step size $\vect{h}_x \defeq (h^1_x,\ldots,h_x^d)\in\ns{R}^d_{>0}$. The number of time points is denoted by $n_t\in\ns{N}$ with step size $h_t = 1/n_t$, $h_t>0$.

We use the \iname{method of Lagrange multipliers} to numerically solve~\eqref{e:opt-prob-v} with Lagrange multipliers $\lambda : \bar{\Omega}\times[0,1]\rightarrow\ns{R}$, $(\vect{x},t)\mapsto\lambda(\vect{x},t)$ (for the hyperbolic transport equation in~\eqref{e:constraint}), and $p:\bar{\Omega}\times[0,1]\rightarrow\ns{R}$, $(\vect{x},t)\mapsto p(\vect{x},t)$ (pressure; for the incompressibility constraint in~\eqref{e:constraint}). We use an \iname{optimize-then-discretize} approach (for a discussion on advantages and disadvantages see~\cite[p.~57ff.]{Gunzburger:2003a}). The resulting optimality conditions are described next.

\subsection{Optimality Conditions}
\label{s:optimality-conditions}

Computing variations of the Lagrangian with respect to perturbations of the state ($m$), adjoint ($\lambda$ and $p$) and control ($\vect{v}$) variables, respectively, yields the (necessary) first-order optimality (\acr{KKT}) conditions (in strong form)
\begin{subequations}\label{e:first-order-opt-cond}
\begin{align}
\p_t m + \igrad m \cdot \vect{v} & = 0
&&{\rm in}\quad\Omega\times(0,1],
\label{e:first-order-forward}
\\
-\p_t\lambda - \idiv(\vect{v}\lambda) & = 0
&&\text{in}\quad \Omega \times [0,1),
\label{e:first-order-adjoint}
\\
\gamma(\idiv\vect{v}) & = 0
&&\text{in}\quad \Omega \times [0,1],
\label{e:first-order-incomp}
\\
\vect{g}\defeq
\beta\D{A}[\vect{v}] + \gamma\igrad p + \vect{f} & = 0
\quad &&\text{in}\quad \Omega \times [0,1],
\label{e:first-order-control}
\end{align}
\end{subequations}

\noindent subject to the initial and terminal conditions
\[
m = m_T\;\;\text{in}\;\;\Omega \times\{0\}
\qquad\text{and}\qquad
\lambda = -(m_1 - m_R)\;\;\text{in}\;\;\Omega \times\{1\}
\]

\noindent and periodic boundary conditions on $\p\Omega$; \eqref{e:first-order-control} is referred to as the \iname{reduced gradient}, where $\vect{f}\defeq\lambda\igrad m$, $\vect{f} : \bar{\Omega} \times [0,1]\rightarrow \ns{R}^d$, $(\vect{x},t)\mapsto\vect{f}(\vect{x},t)$, is the applied \iname{body force} and $\D{A} = \D{B}\D{B}^\adj$ is the G\^ateaux derivative of $\F{S}$. In particular, we have
\begin{subequations}
\begin{align}
\D{A}[\vect{v}]
&= - \ilap\vect{v}
& (\text{$H^1$-regularization};
\gamma = 1 \text{ in \eqref{e:first-order-opt-cond}}),
\label{e:first-variation-H1}
\\
\D{A}[\vect{v}]
&= \ilap^2\vect{v}
& (\text{$H^2$-regularization};
\gamma = 0 \text{ in \eqref{e:first-order-opt-cond}}),
\label{e:first-variation-H2}
\end{align}
\end{subequations}

\noindent respectively. We refer to~\eqref{e:first-order-forward} (hyperbolic initial value problem) as the \iname{state equation}, to~\eqref{e:first-order-adjoint} as the \iname{adjoint equation} (hyperbolic final value problem) and to~\eqref{e:first-order-control} as the \iname{control equation} (elliptic problem). Note that the adjoint equation~\eqref{e:first-order-adjoint} is, likewise to~\eqref{e:first-order-forward}, a scalar conservation law that flows the mismatch between $m_R$ and $m_1$ backward in time. If we neglect the incompressibility constraint in~\eqref{e:constraint}, $\gamma$ in~\eqref{e:first-order-opt-cond} is set to zero (i.e.~\eqref{e:first-order-opt-cond} consists only of~\eqref{e:first-order-forward}, \eqref{e:first-order-adjoint} and~\eqref{e:first-order-control}).

Taking second variations of the Lagrangian yields the system
\begin{subequations}\label{e:second-order-opt-cond}
\begin{align}
\p_t\tilde{m} + \igrad\tilde{m} \cdot \vect{v}
              + \igrad m \cdot \vect{\tilde{v}}
              & = 0 &&\text{in}\quad \Omega \times (0,1],
\label{e:sec-order-istate}\\
- \p_t \tilde{\lambda} - \idiv(\vect{v}\tilde{\lambda}) -
                         \idiv(\vect{\tilde{v}}\lambda) & = 0
                         &&\text{in}\quad \Omega \times [0,1),
\label{e:sec-order-iadj}\\
\gamma(\idiv \vect{\tilde{v}}) & = 0
                               &&\text{in}\quad \Omega \times [0,1],
\label{e:sec-order-incomp}\\
\D{H}[\vect{\tilde{v}}]\defeq
\beta\D{A}[\vect{\tilde{v}}]
+ \gamma\igrad \tilde{p} + \vect{\tilde{f}} & = -\vect{g}
  \quad &&\text{in}\quad \Omega \times [0,1],
\label{e:sec-order-icontrol}
\end{align}
\end{subequations}

\noindent subject to initial and terminal conditions
$\tilde{m}_0 \defeq \tilde{m}(\cdot, 0) = 0$,
$\tilde{m}_0 : \bar{\Omega}\rightarrow\ns{R}$,
$\vect{x} \mapsto\tilde{m}_0(\vect{x})$,
and
$\tilde{\lambda}_1 \defeq \tilde{\lambda}(\cdot,1) = -\tilde{m}_1$,
$\tilde{\lambda}_1 : \bar{\Omega}\rightarrow\ns{R}$,
$\vect{x} \mapsto\tilde{\lambda}_1(\vect{x})$,
$\tilde{m}_1 \defeq \tilde{m}(\cdot, 1)$,
$\tilde{m}_1 : \bar{\Omega}\rightarrow\ns{R}$,
$\vect{x} \mapsto\tilde{m}_1(\vect{x})$,
respectively, and periodic boundary conditions on $\p\Omega$. Here, \eqref{e:sec-order-istate}, \eqref{e:sec-order-iadj} and \eqref{e:sec-order-icontrol} are referred to as the \iname{incremental} state, adjoint and control equations, respectively; the incremental variables are denoted with a tilde. Further, $\D{H}$ in~\eqref{e:sec-order-icontrol} is referred to as the \iname{reduced Hessian} and $\vect{\tilde{f}} \defeq \tilde{\lambda} \igrad m + \lambda \igrad\tilde{m}$, $\vect{\tilde{f}}: \bar{\Omega}\times [0,1]\rightarrow \ns{R}$, $(\vect{x},t)\mapsto\vect{\tilde{f}}(\vect{x},t)$, is the \iname{incremental body force}. The operator $\D{A}$ in \eqref{e:sec-order-icontrol} represents the second variation of $\F{S}$ with respect to the control $\vect{v}$. We use the same symbol as in~\eqref{e:first-order-opt-cond}, since the second variation of $\F{S}$ with respect to $\vect{v}$ is identical to its first variation (the corresponding vectorial differential operators are given in~\eqref{e:first-variation-H1} and~\eqref{e:first-variation-H2}, respectively).

\subsection{Numerical Optimization}
\label{s:numerical-optimization}

We discuss strategies for numerical optimization next. We consider second-order Newton-Krylov methods (see \secref{s:rsqp}) and a first-order Picard method (see \secref{s:picard}).

We use a backtracking line search subject to the \iname{Armijo condition} with search direction $\vect{s}_k \in \ns{R}^n$ and step size $\alpha_k > 0$ at (outer) iteration $k \in\ns{N}_0$ to ensure a sequence of monotonically decreasing objective values $\F{J}^h$ (we use default parameters; see~\cite[Algorithm~3.1,\,p.~37]{Nocedal:2006a}). Note that each evaluation of $\F{J}^h$ requires a forward solve (i.e.\ the solution of~\eqref{e:first-order-forward} to obtain $m_1^h$ given some trial solution $\vect{v}^h_k\in\ns{R}^n$). Therefore, it is desirable to keep the number of line search steps at minimum.

\subsubsection{Inexact Newton-Krylov Method}
\label{s:rsqp}

Applying Newton's method to~\eqref{e:first-order-opt-cond} yields a large KKT system that has to be solved numerically at each \emph{outer iteration} $k$. We will refer to the iterative solution of this system as \emph{inner iterations}\footnote{As opposed to the steps for updating $\vect{v}^h_k$, which we refer to as \emph{outer iterations}; see \algref{a:outer-iteration}.}. In \iname{reduced space methods}, incremental adjoint and state variables are eliminated from the system via block elimination (under the assumption that state and adjoint equations are fulfilled exactly)~\cite{Biros:2005a, Biros:2005b}. We obtain the reduced KKT system
\begin{equation}\label{e:rs-kkt-sys}
\F{H}^h_k \vect{\tilde{v}}^h_k = -\vect{g}^h_k, \quad k \in \ns{N},
\end{equation}

\noindent where $\F{H}^h_k \in \ns{R}^{n \times n}$ corresponds to the reduced Hessian in~\eqref{e:sec-order-icontrol} (i.e.\ the Schur complement of the full Hessian for the control variable $\vect{v}^h$) and $\vect{\tilde{v}}^h_k\in\ns{R}^n$ to the incremental control variable in~\eqref{e:second-order-opt-cond} (which is nothing but the search direction $\vect{s}_k$ mentioned earlier). Further, the right-hand side $\vect{g}_k^h\in\ns{R}^n$ corresponds to the reduced gradient in~\eqref{e:first-order-control}.

The numerical scheme amounts to a sequential solution of the optimality conditions~\eqref{e:first-order-opt-cond} and~\eqref{e:second-order-opt-cond}. \algref{a:outer-iteration} illustrates a realization of an outer iteration\footnote{Note that the scheme in \algref{a:outer-iteration} also applies to the Picard method (see \secref{s:picard}). The only difference is the way we compute the search direction $\vect{s}_k$ in line~\ref{a:outer-iteration:line:search-direction}.}. Note that we eliminate~\eqref{e:first-order-incomp} and~\eqref{e:sec-order-incomp} from the optimality conditions (see \secref{s:elimination}). The inner iteration (i.e.\ the solution of~\eqref{e:rs-kkt-sys}) is what we discuss next.

\begin{algorithm}
\caption
{
Outer iteration of the designed inexact Newton-Krylov method.
}
\label{a:outer-iteration}
\begin{algorithmic}[1]
\STATE
{
$\vect{v}^h_0 \leftarrow 0$;
compute $m^h_0$, $\lambda^h_0$, $\F{J}^h(\vect{v}^h_0)$ and $\vect{g}^h_0$;
$k\leftarrow0$
}
\WHILE{true}
    \STATE
    {
    stop $\leftarrow$ \eqref{e:converged}
    }
    \STATE
    {
    \textbf{if} stop \textbf{break}
    }
    \STATE
    {
    $\vect{s}_k \leftarrow$ solve~\eqref{e:rs-kkt-sys}
    given $m^h_k$, $\lambda^h_k$, $\vect{v}^h_k$ and $\vect{g}^h_k$
    \ccomment{Newton step} \label{a:outer-iteration:line:search-direction}
    }
    \STATE
    {
    $\alpha_k \leftarrow$ perform line search on $\vect{s}_k$
    }
    \STATE
    {
    $\vect{v}^h_{k+1} \leftarrow \vect{v}^h_k + \alpha_k \vect{s}_k$,
    }
    \STATE
    {
    $m^h_{k+1}(t=0) \leftarrow m_T^h$
    }
    \STATE
    {
    $m^h_{k+1} \leftarrow$ solve~\eqref{e:first-order-forward} forward in time
    given $\vect{v}^h_{k+1}$
    \ccomment{forward solve}
    }
    \STATE
    {
    $\lambda^h_{k+1}(t=1) \leftarrow (m^h_R-m^h_{k+1}(t=1))$
    }
    \STATE
    {
    $\lambda^h_{k+1} \leftarrow$ solve~\eqref{e:first-order-adjoint}
    backward in time given $\vect{v}^h_{k+1}$ and $m^h_{k+1}$
    \ccomment{adjoint solve}
    }
    \STATE
    {
    compute $\F{J}^h(\vect{v}_{k+1}^h)$ and $\vect{g}^h_{k+1}$
    given $m^h_{k+1}$, $\lambda^h_{k+1}$ and $\vect{v}^h_{k+1}$
    }
    \STATE
    {
    $k \leftarrow k + 1$
    }
\ENDWHILE
\end{algorithmic}
\end{algorithm}

Forming or storing $\D{H}^h$ in~\eqref{e:rs-kkt-sys} is computationally prohibitive. Therefore, it is desirable to use an iterative solver for which $\D{H}^h$ does not have to be assembled in practice. Krylov-subspace methods are a popular choice~\cite{Benzi:2005a, Biros:2005a, Biros:2005b, Byrd:2008a, Haber:2001a}, as they only require matrix-vector products. We use a \acr{PCG} method, exploiting the fact that $\D{H}^h$ is positive definite (i.e.\ $\D{H}^h \succ 0$; see \secref{s:gauss-newton} for a discussion) and symmetric.

Solving~\eqref{e:rs-kkt-sys} exactly can be prohibitively expensive and might not be justified if an iterate is far from the (true) solution~\cite{Dembo:1982a}. A common strategy is to perform \emph{inexact solves}. That is, starting with a large tolerance for the Krylov-subspace method we successively reduce the tolerance and by that solve more accurately for the search direction, as we approach a (local) minimizer~\cite{Dembo:1983a, Eisenstat:1996a}. This can be achieved with the termination criterion
\begin{equation}\label{e:termination-criterion-ksm}
\|\D{H}^h_\iota \vect{\tilde{v}}^h_{\iota} + \vect{g}_k^h\|_2
\eqdef \|\vect{r}_\iota\|_2
\leq \eta_k \|\vect{g}_k^h\|_2
\end{equation}

\noindent for the Krylov-subspace method. Here, $\eta_k \defeq \min(0.5,\sqrt{\|\vect{g}_k^h\|_2/\|\vect{g}_0^h\|_2})$, $\eta_k\in [0,1)$, is referred to as \iname{forcing sequence} (assuming superlinear convergence; details can be found in~\cite[p.~165ff.]{Nocedal:2006a}); $\iota\in\ns{N}$ in~\eqref{e:termination-criterion-ksm} is the iteration index of the inner iteration (i.e.\ for the iterative solution of~\eqref{e:rs-kkt-sys}) at a given outer iteration $k$.

The course of an inner iteration follows the standard PCG steps (see e.g.~\cite[p.~119, Alg.~5.3]{Nocedal:2006a}). During each inner iteration $\iota$ we have to apply $\D{H}^h$ in~\eqref{e:sec-order-icontrol} to a vector. We summarize this matrix-vector product in~\algref{a:hessian-matvec}. As can be seen, each application of $\D{H}^h$ requires an additional forward and adjoint solve (i.e.\ the solution of the incremental state and adjoint equations~\eqref{e:sec-order-istate} and~\eqref{e:sec-order-iadj}, respectively). This is a direct consequence of the block elimination in reduced space methods.

The number of inner iterations essentially depends on the spectrum of the operator $\D{H}^h$. Typically, $\D{H}^h$ displays poor conditioning. An optimal preconditioner $\mat{P} \in\mathbf{R}^{n \times n}$ renders the number of iterations independent of $n$ and $\beta$. The design of such a preconditioner is an open area of research~\cite{Benzi:2011a, Biros:2008a, Biros:2005a, Biros:2005b, Haber:2001a}. Standard techniques like incomplete factorizations or algebraic multigrid are not applicable, as they require the assembling of $\D{H}^h$ in~\eqref{e:rs-kkt-sys}. Geometric, matrix-free preconditioners are a valid option. This is something we will investigate in the future. Here, we consider a left preconditioner based on the exact spectral inverse of the regularization part of $\D{H}^h$. That is, $\mat{P} \defeq \D{A}^h$ (implementation details can be found in \secref{s:inversion-reg-op}). Note that the PCG method only requires the action of $\mat{P}^{-1}$ on a vector (i.e.\ a matrix-free implementation is in place). Since we use a Fourier spectral method, the cost of our preconditioning amounts to a spectral diagonal scaling. We will refer to this algorithm as Newton-PCG (\acr{N-PCG}) method.

\begin{algorithm}
\caption
{
Hessian matrix-vector product of the designed inexact Newton-Krylov algorithm at outer iteration $k\in\ns{N}$. We illustrate the computational steps required for applying $\D{H}^h$ in~\eqref{e:sec-order-icontrol} to the PCG search direction at inner iteration index $\iota\in\ns{N}$.
}
\label{a:hessian-matvec}
\begin{algorithmic}[1]
    \STATE
    {
    $\tilde{m}^h_\iota(t=0) \leftarrow 0$
    }
    \STATE
    {
    $\tilde{m}^h_\iota \leftarrow$ solve~\eqref{e:sec-order-istate}
    forward in time
    given $m_k^h$, $\vect{v}^h_k$ and $\vect{\tilde{v}}^h_\iota$
    \ccomment{incremental forward solve}
    }
    \STATE{
    $\tilde{\lambda}^h_\iota(t=1) \leftarrow -\tilde{m}^h_\iota(t=1)$
    }
    \STATE
    {
    $\tilde{\lambda}^h_\iota \leftarrow$ solve~\eqref{e:sec-order-iadj}
    backward in time
    given $\lambda_k^h$, $\vect{v}^h_k$ and $\vect{\tilde{v}}^h_\iota$
    \ccomment{incremental adjoint solve}
    }
    \STATE
    {
    apply $\D{H}^h_\iota$ in~\eqref{e:sec-order-icontrol} to the PCG
    search direction given $\lambda_k^h$, $\tilde{\lambda}_\iota^h$,
    $m^h_k$ and $\tilde{m}^h_\iota$
    }
\end{algorithmic}
\end{algorithm}

\subsubsection{GN Approximation}
\label{s:gauss-newton}

Even though $\D{H}^h$ is in the proximity of a (local) minimum by construction positive semidefinite (i.e.\ $\D{H}^h\succeq 0$) it can be indefinite or singular far away from the solution. Accordingly, the search direction is not guaranteed to be a descent direction. One remedy is to terminate the inner iteration whenever negative curvature occurs~\cite{Dembo:1983a}. Another approach is to use a Quasi-Newton approximation. We consider a \acr{GN} approximation $\D{H}_{\rm GN}^h$ instead; we drop certain expressions of $\D{H}^h$, which in turn guarantees that $\D{H}_{\rm GN}^h \succ 0$. In particular, we drop all expressions in~\eqref{e:second-order-opt-cond} in which $\lambda$ appears. Accordingly, we obtain the (continuous) system
\begin{subequations}
\label{e:gauss-newton}
\begin{align}
\p_t \tilde{m} + \igrad\tilde{m} \cdot \vect{v}
+ \igrad m \cdot \vect{\tilde{v}} & = 0
&&\text{in}\quad \Omega \times (0,1],\\
-\p_t \tilde{\lambda}  - \idiv(\vect{v}\tilde{\lambda}) & = 0
&&\text{in}\quad\Omega\times[0,1),\\
\gamma(\idiv \vect{\tilde{v}}) & = 0
&&\text{in}\quad\Omega\times[0,1],\\
\D{H}_{\rm GN}[\vect{\tilde{v}}]\defeq
\beta\D{A}[\vect{\tilde{v}}] + \gamma\igrad p + \tilde{\lambda}\igrad m
& = -\vect{g} \quad &&\text{in}\quad \Omega \times [0,1].
\end{align}
\end{subequations}

\noindent We expect the rate of convergence to drop from quadratic to {(super-)}linear when turning to~\eqref{e:gauss-newton}. However, if the $L^2$-distance can be driven to zero, we recover fast local convergence close to the true solution $\vect{v}^\star$, even if the adjoint variable is neglected. This is due to the fact that~\eqref{e:first-order-adjoint} models the flow of the mismatch backward in time, such that $\lambda\rightarrow 0$ for $\vect{v} \rightarrow \vect{v}^\star$. We refer to this method as (\emph{inexact}) \acr{GN-PCG} method~\cite{Biros:2005a, Biros:2005b}. We remark that all algorithmic details described in this note apply to both Newton-Krylov methods.

\subsubsection{Picard Method}
\label{s:picard}

We consider a globalized Picard iteration (fixed point iteration) in addition to the described Newton-Krylov methods. Based on~\eqref{e:first-order-control} we have
\begin{equation}\label{e:picard}
\vect{v}^h_{k+1} = -(\beta\D{A}^{h})^{-1}[\gamma\nabla^h p^h_k + \vect{f}^h_k].
\end{equation}

\noindent Since we use Fourier spectral methods, the inversion of $\D{A}^h$ in~\eqref{e:picard} comes at the cost of a diagonal scaling (implementation details can be found in \secref{s:inversion-reg-op}). Accordingly, this scheme does not require the (iterative) solution of a linear system. However, it potentially results in a larger number of outer iterations until convergence as we expect the optimization problem to be poorly conditioned.

We do not directly use the solution of~\eqref{e:picard} as a new iterate but compute a search direction $\vect{s}_k$ instead. This in turn allows us to perform a line search on $\vect{s}_k$. That is, we subtract the new from the former iterate. This scheme can be viewed as a gradient descent in the function space induced by $\fs{W}$ (i.e.\ a \emph{preconditioned gradient descent scheme}; see \secref{s:relation-to-lddmm}).

Note that $\vect{s}_k$ is, in contrast to Newton methods, arbitrarily scaled. Therefore, we provide an augmented implementation that tries to estimate an optimal scaling during the course of optimization. Details can be found in \secref{s:parameters}.

\subsubsection{Termination Criteria}
\label{s:stopping-criteria}

The termination criteria are in accordance with \cite{Modersitzki:2009a} (see \cite[305\,ff.]{Gill:1981a} for a discussion) given by
\begin{equation}
\label{e:stopping-criteria}
\begin{array}{r@{\quad} r@{\;} l@{\;} l}
(C1)
& \F{J}^h(\vect{v}^h_{k-1}) - \F{J}^h(\vect{v}^h_k)
& < \tau_{\F{J}} (1 + \F{J}^h(\vect{v}^h_0)),
\\
(C2)
& \|\vect{v}^h_{k-1} - \vect{v}^h_k\|_{\infty}
& < \sqrt{\tau_{\F{J}}} (1 + \|\vect{v}^h_k\|_{\infty}),
\\
(C3)
& \|\vect{g}_k^h\|_{\infty}
& < \sqrt[3]{\tau_{\F{J}}} (1 + \F{J}^h(\vect{v}^h_0)),
\\
(C4)
& \|\vect{g}_k^h\|_{\infty}
& < \num{1E3}\,\epsilon_{\rm mach},
\\
(C5)
& k
& > n_{\rm opt}.
\end{array}
\end{equation}

\noindent Here, $\tau_{\F{J}} > 0$ is a user-defined tolerance, $\epsilon_{\rm mach} > 0$ is the machine precision and $n_{\rm opt} \in\ns{N}$ is the maximal number of outer iterations. The algorithm is terminated if
\begin{equation}\label{e:converged}
    \left\{(C1) \land (C2) \land (C3)\right\} \; \lor \; (C4) \; \lor \; (C5),
\end{equation}

\noindent where $\lor$ denotes the \iname{logical or} and $\land$ the \iname{logical and} operator, respectively.

\subsection{Algorithmic Details}
\label{s:implementation}

This section provides additional specifics on the implementation. In particular, we describe\bipa\item the numerical discretization (see \secref{s:discretization}), \item the parameterization in time (see \secref{s:expansion}), \item the inversion of the operator $\D{A}^h$ (see \secref{s:inversion-reg-op}), and \item strategies for the parameter selection (see \secref{s:parameters}).\eipa

\subsubsection{Numerical Discretization}
\label{s:discretization}

We use a (regular) nodal grid for the discretization in space and time. The problem is defined on the space-time interval $\Omega\times[0,1]$, where $\Omega\defeq (-\pi,\pi)^d$. Accordingly, we obtain the time step size via $h_t= 1/n_t$. The cell size (pixel or voxel size) $\vect{h}_x \defeq (h_x^1,\ldots, h_x^d)\in\ns{R}^d_{\geq0}$ for a spatial grid cell can be computed via $h_x^i = 2\pi/n_x^i$, $i=1,\ldots,d$, where $n_x^i$ is the number of grid points along the $i$th spatial direction $x^i$.

The derivative operators are discretized via Fourier spectral methods~\cite{Boyd:2000a}. The time integrator for the forward and adjoint solves is an explicit second-order Runge-Kutta (RK2) method, which, in connection with Fourier spectral methods, displays minimal numerical diffusion.

Following standard numerical theory for hyperbolic equations, the step size $h_t>0$ is bounded from above by $h_{t,\max} \defeq \epsilon_{\text{CFL}} / \max(\|\vect{v}^h\|_\infty \oslash \vect{h}_x)$, $h_{t,\max} > 0$ (Courant-Friedrich-Lewy (\acr{CFL}) condition). Here, $\oslash$ denotes a Hadamard division and $\epsilon_{\text{CFL}}>0$ is the CFL number. The theoretical bound for $h_{t,\max}$ is attained for $\epsilon_{\text{CFL}} = 1$. We use $\epsilon_{\text{CFL}}= 0.2$ for all experiments. Since we use a spectral Galerkin method in time (see \secref{s:expansion}), we can adaptively adjust $n_t$ (and therefore $h_t$) for the forward and adjoint solves as required by the CFL condition.

\subsubsection{Spectral Galerkin Method}
\label{s:expansion}

To reduce the number of unknowns, $\vect{v}$ is expanded \emph{in time} in terms of basis functions $b_l : [0,1] \rightarrow \ns{R}$, $t \mapsto b_l(t)$, $l = 1,\ldots,n_c$,
\begin{equation}\label{e:expansion}
    \vect{v}(\vect{x},t) = \sum_{l=1}^{n_c} b_l(t) \vect{v}_l(\vect{x}),
\end{equation}

\noindent where $\vect{v}_l : \ns{R}^d \rightarrow\ns{R}^d$, $\vect{x} \mapsto \vect{v}_l(\vect{x})$, is a coefficient field. The coefficients $\vect{v}_l$ are the new unknowns of our problem. This reduces the number of unknowns in time from $n_t$ to $n_c$, where $n_c \ll n_t$. Thus, we can invert for a stationary ($n_c=1$) or a nonstationary velocity field as required. Nothing changes in our formulation---just the number of unknowns.

We use \iname{Chebyshev polynomials} as basis functions $b_l$ on account of their excellent approximation properties as well as their orthogonality (see \secref{s:appendix-expansion} for details). The expansion~\eqref{e:expansion} solely affects $\F{S}$ and the (incremental) control equation (i.e.~\eqref{e:first-order-control} and~\eqref{e:sec-order-icontrol}); $\vect{v}$ is computed from the coefficient fields $\vect{v}_l$, $l=1,\ldots,n_c$, during the forward and adjoint solves according to~\eqref{e:expansion}.

\subsubsection{Inversion: Regularization Operators}
\label{s:inversion-reg-op}

The Picard iteration in~\eqref{e:picard} as well as the preconditioning of~\eqref{e:rs-kkt-sys} require the inversion of the differential operator $\D{A}^h$. Since we use Fourier spectral methods this inversion can be accomplished at the cost of a spectral diagonal scaling. However, $\D{A}^h$ has a nontrivial kernel (which only includes constant functions due to the periodic boundary conditions). We make $\D{A}^h$ invertible by setting the base frequency of the inverse of $\D{A}^h$ (including the scaling by $\beta$) to one. This ensures not only invertibility but also that the constant part of the (incremental) body force $\vect{f}$ (or $\vect{\tilde{f}}$, respectively) remains in the kernel of our regularization scheme. This in turn allows us to invert for constant velocity fields.

\subsubsection{Elimination of $p$ and $\tilde{p}$}
\label{s:elimination}

In our numerical scheme, we eliminate $p$ and by that~\eqref{e:first-order-incomp} from~\eqref{e:first-order-opt-cond}. Details on the derivation can be found in \secref{s:elimination-details}. We obtain
\begin{equation}\label{e:eliminated-control}
\vect{\hat{g}} \defeq \beta\D{A}[\vect{v}]
+ \D{K}[\vect{f}] = 0,
\quad
\D{K}[\vect{f}]
\defeq
-\igrad(\ilap^{-1}(\idiv \vect{f}))
+\vect{f},
\end{equation}

\noindent to replace~\eqref{e:first-order-control}, where $\vect{f}$ is the body force as defined in \secref{s:optimality-conditions} and $\D{A}$ the first variation of $\F{S}$ with respect to $\vect{v}$ (see~\eqref{e:first-variation-H1} and~\eqref{e:first-variation-H2}, respectively). It immediately follows that we obtain
\begin{equation}\label{e:eliminated-icontrol}
\hat{\D{H}}[\vect{\tilde{v}}]
\defeq\beta\D{A}[\vect{\tilde{v}}]
+ \D{K}[\vect{\tilde{f}}] = -\vect{\hat{g}}
\end{equation}

\noindent to replace~\eqref{e:sec-order-icontrol}; $\vect{\tilde{f}}$ denotes the incremental body force defined in \secref{s:optimality-conditions}.

\subsubsection{Parameter Selection}
\label{s:parameters}

To the extent possible, it is desirable to design a numerical scheme that does not require a selection of parameters (\iquote{black-box solver}). This is challenging for previously unseen data. In general, the user should only be required to decide on the following:
\begin{itemize}
    \item The desired accuracy of the inversion (controlled by the tolerance $\tau_{\F{J}}$; see \secref{s:stopping-criteria}).
    \item The desired properties of the mapping $\vect{y}$ (controlled by $\epsilon_{\theta}$ or $\epsilon_F$, respectively; see below).
    \item The budget we are willing to assign to the computation (controlled by $n_{\textrm{opt}}$; see \secref{s:stopping-criteria}).
\end{itemize}

\noindent For the purpose of this numerical study we proceed as follows:

\textbf{Optimization}: We set the maximum number of iterations $n_{\text{opt}}$ (see~\eqref{e:stopping-criteria}) to $\num{1E6}$, as we do not want our algorithm to terminate early (i.e.\ we make sure that we terminate only if either we reach the defined tolerances or we no longer observe a decrease in $\F{J}^h$). For the convergence study in \secref{s:numerical-experiments}, we use the relative change of the $\ell^\infty$-norm of the gradient $\vect{g}^h$ as a stopping criterion, as we are interested in studying convergence properties. This enables an unbiased comparison in terms of the required work to solve an optimization problem up to a desired accuracy. In particular, we terminate the optimization if the relative change of the reduced gradient $\vect{g}^h$ is larger than or equal to three orders of magnitude.

Following standard textbook literature~\cite{Gill:1981a, Modersitzki:2009a} we use the stopping criteria in~\eqref{e:stopping-criteria} for the remainder of the experiments. We set the tolerance to $\tau_{\F{J}} = \num{1E-3}$. We \emph{qualitatively} did not observe significant differences in the final results for the experiments performed in this study, when turning to smaller tolerances. We will further elaborate on the required accuracy for the inversion  (i.e.\ the registration quality) in a follow-up paper.

The tolerance of the PCG method is set as discussed in \secref{s:rsqp} (see~\eqref{e:termination-criterion-ksm}). The maximal number of iterations for the PCG method is set to $n$ (order of the reduced KKT system in~\eqref{e:rs-kkt-sys}). In theory, this guarantees that the PCG method converges to a solution. This choice not only ensures that we provide an unbiased study (i.e.\ we do not terminate early) but also makes sure that we do not miss any issues in the implementation or parameter selection. We converged for all experiments conducted in this study after only a fraction of $n$ inner iterations. This statement is confirmed by the reported number of PDE solves\footnote{`'PDE solve`' refers to the solution of one of the hyperbolic PDEs that appear in the optimality system~\eqref{e:first-order-opt-cond} and~\eqref{e:second-order-opt-cond}.}.

For all our experiments we initialized the line search with a factor of $\alpha_k = 1$ (see~\secref{s:numerical-optimization}). This is a sensible choice, as search directions obtained from second-order methods are nicely scaled (i.e.\ we expect $\alpha_k$ to be 1). However, this is not the case for the Picard scheme (i.e. the preconditioned gradient descent). Our implementation features an option to memorize the scaling of $\vect{s}_k$ for the next outer iteration. That is, we introduce an additional scaling factor $\tilde{\alpha}_k > 0$ that is applied to $\vect{s}_k$ before entering the line search (initialized with $\tilde{\alpha}_k=1$). If the line search kicks in, we downscale $\tilde{\alpha}_k$ by $\alpha_k$. On the contrary, we upscale $\tilde{\alpha}_k$ by a factor of two if $\alpha_k = 1$.

\textbf{PDE Solver}: The number of time steps $n_t$ is bounded from below due to stability requirements (see \secref{s:discretization}). Since we use an expansion in time (see \secref{s:expansion}), it is possible to adaptively adjust $n_t$, so that numerical stability is attained.

However, we fix $n_t$ for the numerical experiments in \secref{s:numerical-study} as we are interested in studying the convergence behavior with respect to the employed grid size. We set $n_t$ to $4\max(\vect{n}_x)$. This is a pessimistic choice. If we still encounter instabilities (as judged by monitoring the CFL condition (see \secref{s:discretization})), $\vect{s}_k$ is scaled by a factor of 0.5 until numerical stability is attained, \emph{before} entering the line search. For all numerical experiments conducted in this study, we did not observe any instabilities for the Newton-Krylov methods. However, for the Picard method we observed instabilities in the case when we did not consider the rescaling procedure detailed above. This is due to the fact that $\vect{s}_k$ is arbitrarily scaled for first-order methods (as opposed to second-order methods). By introducing the additional scaling parameter $\tilde{\alpha}_k$ we could stabilize the Picard method---we did not observe a violation of the CFL condition for any of the conducted experiments (for $n_t$ fixed).

\textbf{Regularization}: Estimating an optimal value for $\beta$ is an area of research by itself. A variety of methods has been designed (see e.g.\ \cite{Vogel:2002a}). A key difficulty is computational complexity. Methods based on the assumption that differences between model output and observed data are associated with random noise (such as generalized cross validation) might not be reliable in the context of nonrigid image registration. This is due to the fact that the noise in the images is likely to be highly structured~\cite{Haber:2006a}. Another possibility is to estimate the regularization parameter on the basis of the spectral properties of the Hessian (see \secref{s:res-spectral-analysis}). That is, we can estimate the condition number of the problem during the PCG solves for the unregularized problem using the Lanczos algorithm (see \cite[p.~528]{Golub:1996a}). We can do this very efficiently by initializing the problem with a zero velocity field. Given $\vect{v}$ is zero, the application of the Hessian within the PCG is computationally inexpensive, as a lot of the terms in the optimality systems drop (see \secref{s:optimality-conditions}). However, the level of regularization depends not only on properties of the data but also on regularity requirements on $\vect{y}$.

Another common strategy is to perform a \iname{parameter continuation} in $\beta$ (see e.g.~\cite{Haber:2000a, Haber:2006a}). In~\cite{Haber:2006a} it has been suggested to inform the algorithm about the required regularity of a solution on the basis of a lower bound on the $L^2$-distance between the reference and the deformed template image. The decision on such bound, however, might not be intuitive for practitioners. Further, one is ultimately interested not only in a small residual but also in a bounded determinant of the deformation gradient. Therefore, we propose to inform the algorithm on regularity requirements in terms of a lower bound $\epsilon_F\in (0,1)$ on $\det(\mat{F}_1^h)$ (i.e.\ a bound on the tolerable compression of a volume element). If the Stokes regularization scheme ($\gamma=1$ in \eqref{e:constraint} and $\F{S}$ is an $H^1$-seminorm) is considered, bounds on geometric constraints of the deformation of a volume element can be used. In particular, we use a lower bound $\epsilon_{\theta} > 0$ on the acute angle of a grid element. The upper bound on the obtuse angle is given by $2\pi - \epsilon_{\theta}$. Note that it is actually necessary to monitor geometric properties to guarantee a local diffeomorphism; a lower bound on $\det(\mat{F}_1^h)$ is not sufficient.

Our algorithm proceeds as follows: In a first step, the registration problem is solved for a large value of $\beta$ ($\beta=1$ in our experiments) so that we underfit the data\footnote{Note that for large $\beta$ the optimization problem is almost quadratic, so that Newton-Krylov methods converge quickly.}. Subsequently, $\beta$ is reduced by one order of magnitude until we reach $\epsilon_F$ (or $\epsilon_{\theta}$). From there on a binary search is performed. The algorithm is terminated if the change in $\beta$ is below 5\% of the value for $\beta$, for which $\epsilon_F$ (or $\epsilon_{\theta}$) was breached. We add a lower bound of \num{1E-6} on $\beta$ as well as a lower bound for the relative change of the $L_2$-distance of \num{1E-2} to ensure that we do not perform unnecessary work. We never reached these bounds for the experiments conducted in this study.

\textbf{Presmoothing}: A numerical challenge in image computing is that images are functions of bounded variation. Therefore, an accurate computation of the derivatives becomes more involved. A common approach to ensure numerical stability and avoid the Gibbs phenomena is to reduce high-frequency information in the data. We use a Gaussian smoothing, which is parametrized by a user-defined standard deviation $\sigma > 0$. We experimentally found a value of $\sigma=2\pi/\min(\vect{n}_x)$ to be adequate for the problems at hand. However, we note that we increased $\sigma$ by a factor of 2 for one set of experiments in \secref{s:res-convergence}. We also note that we implemented a method for grid and scale continuation for the images. This avoids the problem of deciding on $\sigma$. We will additionally investigate an automatic selection strategy for $\sigma$ in a follow-up paper.

It is important to note that the sensitivity of second-order derivatives to noise in the data is problematic. Therefore, we refrain from applying the N-PCG method to nonsmooth images.

\section{Numerical Experiments}
\label{s:numerical-experiments}

We report results only in two dimensions. We test the algorithm on real-world and synthetic registration problems (see \secref{s:data}). The measures to analyze the registration results are summarized in \secref{s:measures}. We conduct a numerical study (see~\secref{s:numerical-study}), which includes an analysis of\bipa\item the spectral properties of the Hessian (see~\secref{s:res-spectral-analysis}), \item grid convergence (see~\secref{s:res-convergence}) and \item the effects of varying the number of the unknowns in time (see~\secref{s:res-expansion})\eipa. We additionally report results for a fully automatic registration on high-resolution images based on the designed parameter continuation in $\beta$ (see \secref{s:res-parameter-continuation}).

\subsection{Data}
\label{s:data}

\begin{figure}
\centering
\includegraphics[width=0.9\textwidth]
{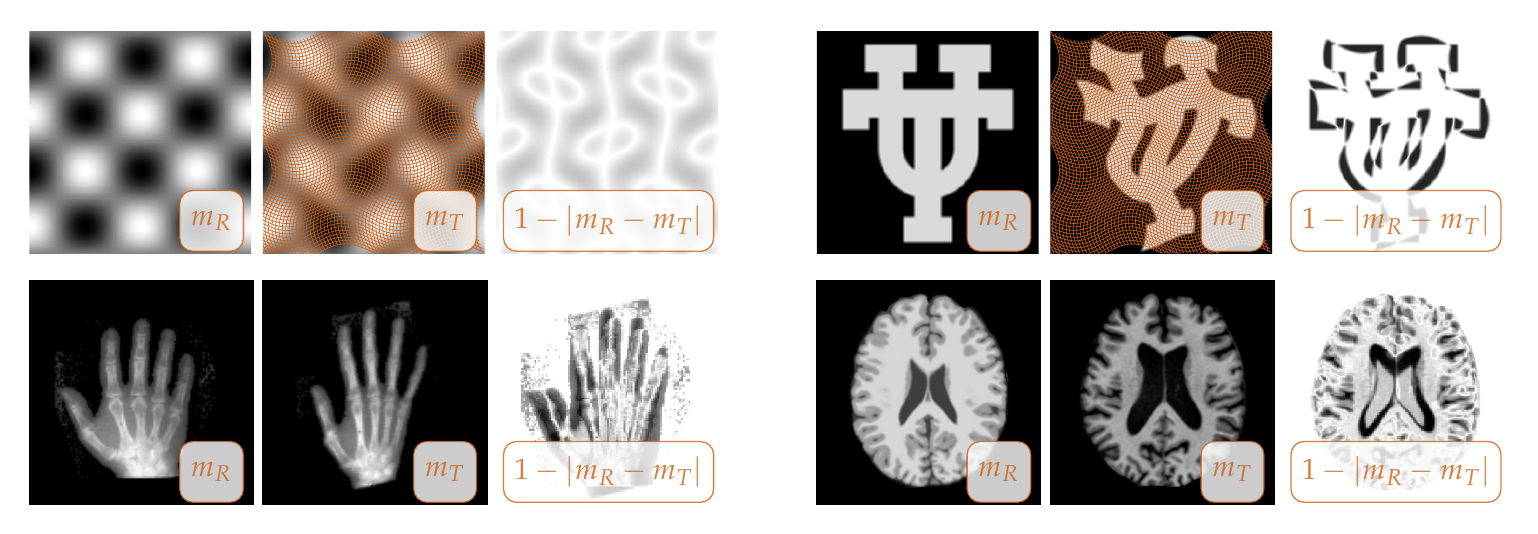}
\caption
{
Synthetic and real-world registration problems. All images have been normalized to an intensity range of $[0,1]$. The registration problems are referred to as \iquote{sinusoidal images} (top row, left), \iquote{UT images} (top row, right), \iquote{hand images} (bottom row, left) and \iquote{brain images} (bottom row, right). Each row displays the reference image $m_R$, the template (deformable) image $m_T$ and a map of their pointwise difference (from left to right as identified by the inset in the images). We provide an illustration of the deformation pattern $\vect{y}$ (overlaid onto $m_T$) for the synthetic problems. This mapping is computed from $\vect{v}^\star$ in~\eqref{e:synvelocity} via~\eqref{e:displacement}.
}
\label{f:regprob}
\end{figure}

We consider synthetic and real-world registration problems\footnote{The \iquote{hand images} are taken from \cite{Modersitzki:2009a}.}. These are illustrated in \figref{f:regprob}. All images have been normalized to an intensity range of $[0,1]$. The synthetic problems are constructed by solving the forward problem to create an artificial template image $m_T$ given some image $m_R$ and some velocity field $\vect{v}^\star$ (\iquote{sinusoidal images} and \iquote{UT images} in~\figref{f:regprob}). Here, $\vect{v}^\star$ is chosen to live on the manifold of divergence-free velocity fields to provide a testing environment for the Stokes regularization scheme. Further, $\vect{v}^\star$ is by construction assumed to be constant in time (i.e.\ $n_c=1$). In particular, we have
\begin{equation}\label{e:synvelocity}
v^{i,\star}(x^i,t)
= 0.5\sin(x^i) \cos(x^i)
\quad\forall t \in [0,1],
\quad i=1,\ldots,d.
\end{equation}

\subsection{Measures of Performance}
\label{s:measures}

We report the number of the (outer) iterations and the number of the hyperbolic PDE solves to assess the computational work load. The latter is a good proxy for the wall clock time. It provides a transparent comparison between the designed first and second-order methods, given that the number of the hyperbolic PDE solves varies between these methods: Two hyperbolic PDE solves are required during each iteration of the Picard method (we have to solve~\eqref{e:first-order-forward} and~\eqref{e:first-order-adjoint} to evaluate the reduced gradient in~\eqref{e:first-order-control}; see also \algref{a:outer-iteration}). For the Newton-Krylov methods, we require an additional two hyperbolic PDE solves per inner iteration (we have to solve \eqref{e:sec-order-istate} and \eqref{e:sec-order-iadj} to compute the Hessian matrix-vector product given in~\eqref{e:sec-order-icontrol}; see also \algref{a:hessian-matvec}). Each evaluation of $\F{J}^h$ in \eqref{e:opt-prob-v} (i.e.\ each line search step) requires an additional hyperbolic PDE solve (we have to compute the deformed template image at $t=1$ by solving \eqref{e:first-order-forward}). Note that the solution of the forward and adjoint problems is the key bottleneck of our algorithm. We will report the wall clock times in a follow-up paper, in which we extend the current framework to a 3D implementation. This study focuses on algorithmic features.

We report the relative change of\bipa\item the $L^2$-distance, \item the objective functional $\F{J}^h$ and \item the $\ell^{\infty}$-norm of the (reduced) gradient $\vect{g}^h$ to assess the quality of the inversion\eipa. We additionally report values for the determinant of the deformation gradient to study local deformation properties. These measures are defined more explicitly in~\tabref{t:quantities-of-interest}.

We visually support this quantitative analysis on basis of snapshots of the registration results. Information on the reconstruction accuracy can be obtained from pointwise maps of the residual difference between $m_R$ and $m_1$. The deformation regularity and the mass conservation can be assessed via images of the pointwise determinant of the deformation gradient and/or of a deformed grid overlaid onto $m_1$. Details on how these are obtained and on how to interpret them can be found in \secref{s:map-and-def-grad}.

\begin{table}
\caption
{
Overview of the quantitative measures that are used to assess the registration performance. We report the number of outer iterations (steps for updating the control variable $\vect{v}^h$) and the number of the hyperbolic PDE solves (i.e.\ how often we have to solve one of the hyperbolic PDEs \eqref{e:first-order-forward}, \eqref{e:first-order-adjoint}, \eqref{e:sec-order-istate} and \eqref{e:sec-order-iadj} that appear in the optimality systems) to assess the work load. We report the relative change of the $L^2$-distance, the objective and the reduced gradient to asses the quality of the inversion. We report values for the determinant of the deformation gradient to assess the regularity of the computed deformation map. We report the relative power spectrum of the coefficients $\vect{v}^h_l$ to assess, which of the coefficients of the expansion in \eqref{e:expansion} are significant.
}
\label{t:quantities-of-interest}
\tabadjust
\begin{tabular}{|lc|}\hline
\rowcolA
  Description
& Definition
\\\hline\rowcolC
\# of required outer iterations
& $k^\star$
\\\rowcolD
\# of required hyperbolic PDE solves
& $n_{\rm PDE}$
\\\rowcolC
relative change of $L^2$-distance
& $\|m_1^h-m_R^h\|_{2,\text{rel}}
\defeq \|m_1^h-m_R^h\|_2^2/\|m_T^h-m_R^h\|_2^2$
\\\rowcolD
relative change of objective value
& $\delta\F{J}^h_{\text{rel}}
\defeq\F{J}^h(\vect{v}^h_{k^\star})/\F{J}^h(\vect{v}^h_0)$
\\\rowcolC
relative change of reduced gradient
& $\|\vect{g}^h\|_{\infty,\text{rel}}
\defeq\|\vect{g}^h_{k^\star}\|_{\infty}/\|\vect{g}^h_0\|_{\infty}$
\\\rowcolD
determinant of deformation gradient
& $\det(\mat{F}_1^h)$
\\\rowcolC
relative power spectrum of $\vect{v}^h_l$
& $\|\vect{v}^h_l\|_{2,{\rm rel}}
\defeq \|\vect{v}_l^h\|_2 / \|\{\vect{v}^h_{l'}\}_{l'=1}^{n_c}\|_{2}$
\\\hline
\end{tabular}
\end{table}

\subsection{Numerical Study}
\label{s:numerical-study}

We study the spectral properties of the Hessian (see~\secref{s:res-spectral-analysis}), grid convergence (see~\secref{s:res-convergence}) as well as the influence of an increase in the number of the unknowns in time (see~\secref{s:res-expansion}).

\subsubsection{Spectral Analysis}
\label{s:res-spectral-analysis}
\secheadskip

\ipoint{Purpose} We study the ill-posedness and the ill-conditioning of the problem at hand. We report spectral properties of $\mathcal{H}^h$. We study the eigenvalues and the eigenvectors with respect to different choices for $\beta$. We study the differences between plain $H^2$-regularization and the Stokes regularization scheme ($H^1$-regularization).

\ipoint{Setup} This study is based on the \iquote{UT images} (the true solution $\vect{v}^{h,\star}$ is divergence-free; see \secref{s:data} for more details on the construction of this synthetic registration problem; $\vect{n}_x = (64,64)^\T$ and $n_c = 1$ so that $n = 8192$). The eigendecomposition $\mat{V} \mat{\Lambda} \mat{V}^{-1}$, $\mat{V} = (\vect{\nu}_i)_{i=1}^n$, $\vect{\nu}_i\in\ns{R}^n$, $\|\vect{\nu}_i\|_2 = 1$, $\mat{\Lambda} = \operatorname{diag}(\Lambda_{11}, \ldots, \Lambda_{nn})$, $\Lambda_{ii}>0$, is computed at the true solution $\vect{v}^{h,\star}$ to guarantee that $\mathcal{H}^h \succ 0$. The spectrum is computed for three different choices of $\beta$: for the unregularized problem ($\beta=\num{0}$), an empirically determined (moderate) value ($\beta=\num{1E-3}$) and solely for the regularization model ($\beta=\num{1E6}$).

\begin{figure}
\centering
\includegraphics[width=\textwidth]
{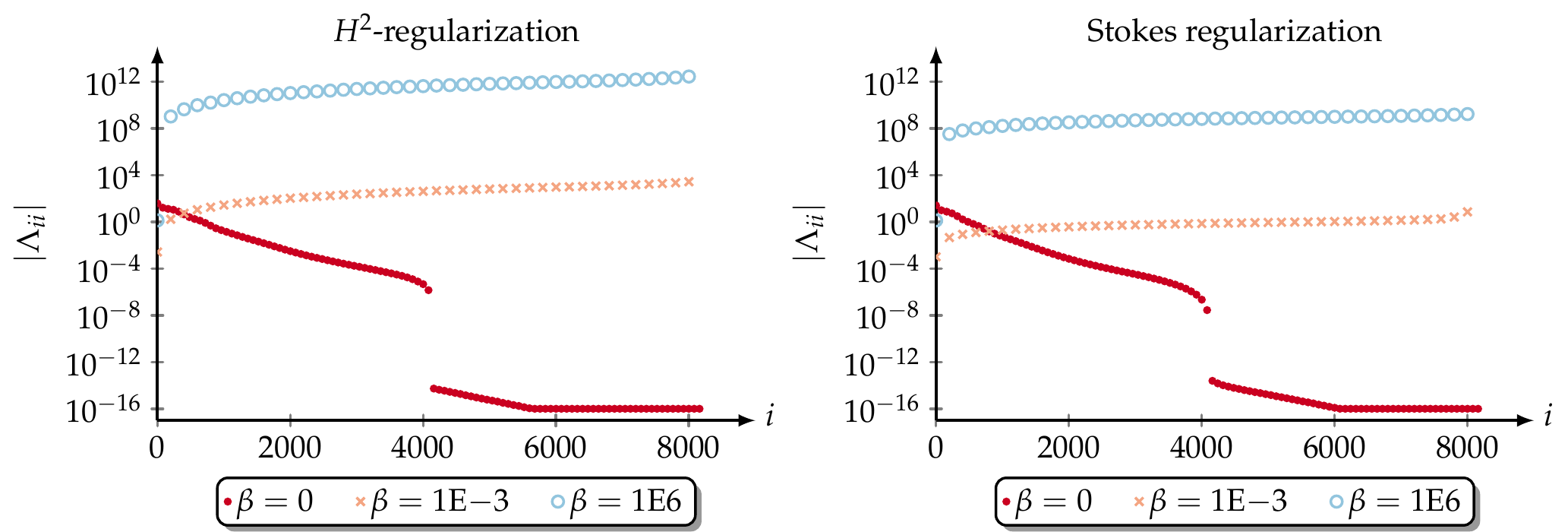}
\caption
{
Trend of the absolute value of the eigenvalues $\Lambda_{ii}$, $i=1,\ldots,8192$, of the reduced Hessian $\D{H}^h$ for plain $H^2$-regularization ($\gamma=0$; left) and the Stokes regularization scheme ($\gamma=1$; $H^1$-regularization; right) for $\beta\in\{0,\num{1E-2},\num{1E6}\}$ (as indicated in the legend of the plots). We report the trend of the entire set of 8192 eigenvalues. The test problem is the \iquote{UT images} (see \secref{s:data} for details on the construction of this synthetic registration problem; $\vect{n}_x = (64,64)^\T$ and $n_c=1$). The Hessian is computed at the true solution $\vect{v}^{h,\star}$ to ensure that $\D{H}^h \succ 0$ (this statement is confirmed by the values reported in \tabref{t:eigenvalues}). The eigenvalues (absolute value) are sorted in descending order for the unregularized problem (i.e.\ for $\beta=0$) and in ascending order otherwise (i.e.\ for $\beta=\num{1E-3}$ and $\beta=\num{1E6}$).
}
\label{f:eigenvalues}
\end{figure}

\begin{figure}
\centering
\includegraphics[width=0.85\textwidth]
{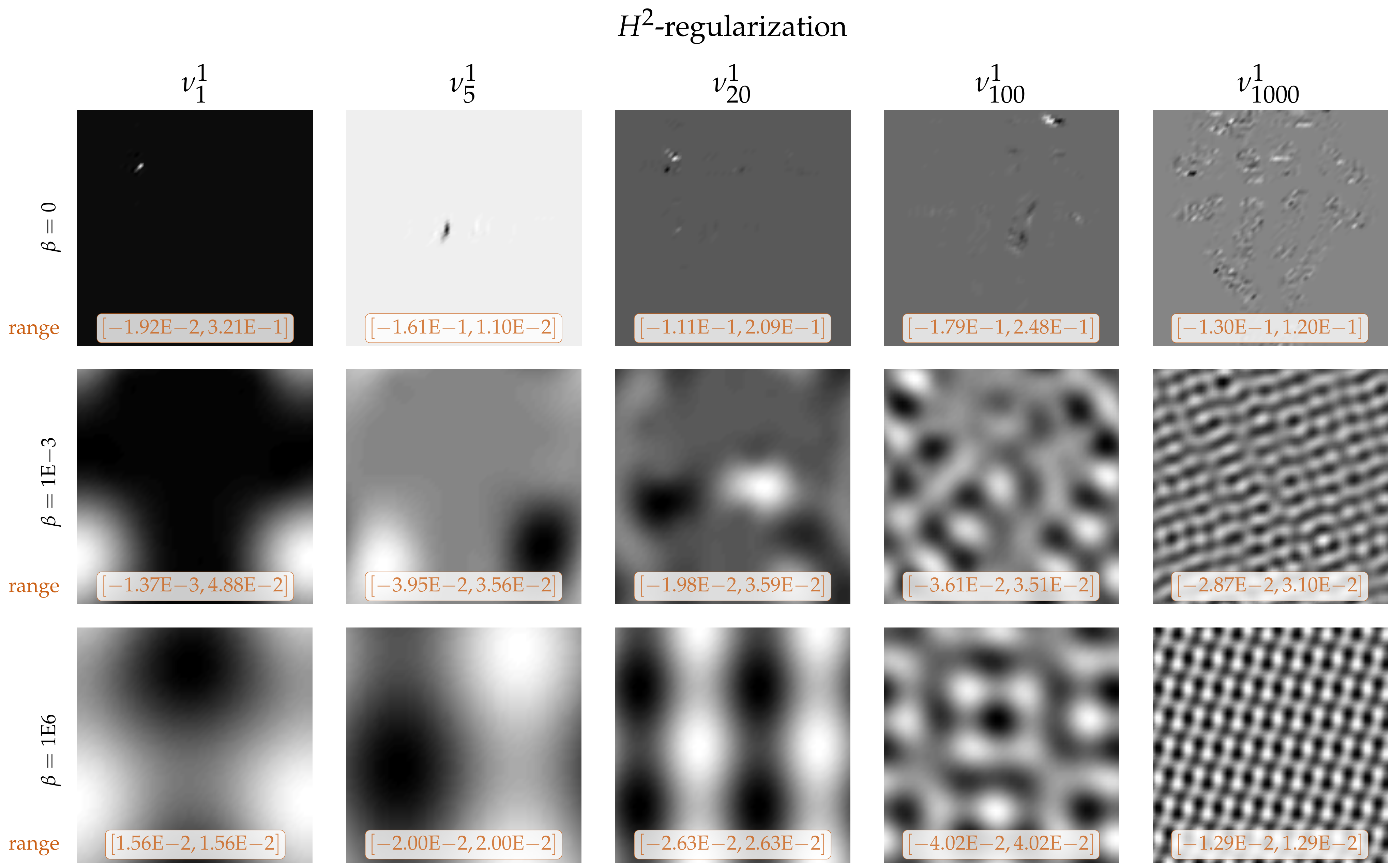}
\\\bigskip
\includegraphics[width=0.85\textwidth]
{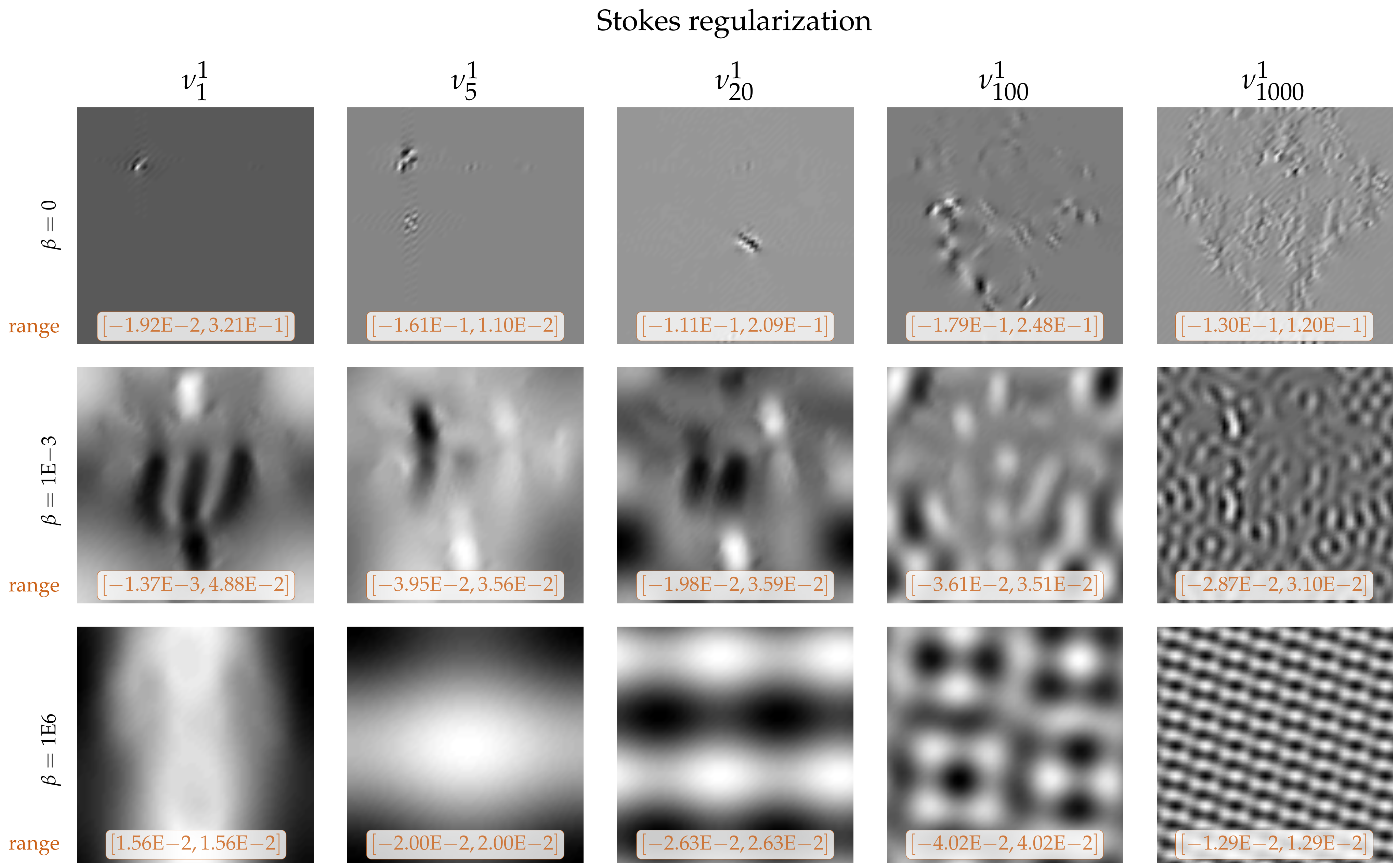}
\caption
{
Eigenvector plots of the reduced Hessian $\D{H}^h \in\ns{R}^{n \times n}$, $n=8192$, for $\beta\in\{0,\num{1E-3},\num{1E6}\}$ for plain $H^2$-regularization ($\gamma=0$; top) and for the Stokes regularization scheme ($\gamma=1$; $H^1$-regularization; bottom). The results correspond to the eigenvalue plots reported in \figref{f:eigenvalues}. We refer to \figref{f:eigenvalues} and the text for details on the experimental setup. Each plot provides the spatial variation of the portion of an eigenvector $\vect{\nu}_i\in\ns{R}^n$ associated with the first component of the coefficient field $\vect{v}_l^h$, $l=n_c$, $n_c=1$. The individual plots correspond to the eigenvalues $\Lambda_{ii}>0$, $i=1,5,20,100,1000$ in \figref{f:eigenvalues}. The range of the values for $\vect{\nu}_i^1$ is provided below each plot.
}
\label{f:eigenvectors}
\end{figure}

\ipoint{Results} \figref{f:eigenvalues} displays the trend of the absolute value of the eigenvalues $\Lambda_{ii}$, $i = 1,\ldots,n$. They are sorted in \emph{descending} order for $\beta=0$ and in \emph{ascending} order otherwise. If an eigenvalue drops below machine precision (i.e.\ $\num{1E-16}$), it is set to $\num{1E-16}$ (only for visualization purposes). The extremal real and imaginary part of the eigenvalues is summarized in \tabref{t:eigenvalues}. \figref{f:eigenvectors} provides the spatial variation of the eigenvectors $\vect{\nu}_i\in\mathbf{R}^n$ that correspond to the eigenvalues $\Lambda_{ii}$, $i\in\{1,5,20,100,1000\}$, in \figref{f:eigenvalues} with respect to different choices for $\beta$ and different regularization schemes. We only display the first component $\vect{\nu}_i^1$ of the coefficient field $\vect{\nu}_i \defeq (\vect{\nu}_i^1,\vect{\nu}_i^2)$. The pattern for the second component is (qualitatively) alike.

\begin{table}
\caption
{
Extrema of the eigenvalues $\Lambda_{ii}$, $i=1,\ldots,8192$, of the reduced Hessian reported in \figref{f:eigenvalues}. We report values for plain $H^2$-regularization ($\gamma=0$; top block) and the Stokes regularization scheme ($\gamma=1$; $H^1$-regularization; bottom block). We refer to \figref{f:eigenvalues} and the text for details on the experimental setup. We report the smallest and the largest real part as well as the largest absolute value of the imaginary part of the eigenvalues $\Lambda_{ii}$ with respect to different choices of the regularization parameter $\beta\in\{0,\num{1E-2},\num{1E6}\}$.
}
\label{t:eigenvalues}
\tabadjust
\begin{tabular}{|llll|}\hline
\multicolumn{4}{|c|}{\cellcolA $H^2$-regularization ($\gamma=0$)}\tabularnewline\hline
  $\beta$
& $\min\{(\real(\Lambda_{ii}))_{i=1}^n\}$
& $\max\{(\real(\Lambda_{ii}))_{i=1}^n\}$
& $\max\{(|\imag(\Lambda_{ii})|)_{i=1}^n\}$
\\ \hline \rowcolD
  $\num{0}$
& $\num{-8.350057e-15}$
& $\num{3.720556e+01}$
& $\num{3.256531e-07}$
\\ \rowcolC
  $\num{1E-3}$
& $\num{2.593468e-03}$
& $\num{4.194738e+03}$
& $\num{0.000000e+00}$
\\ \rowcolB
  $\num{1E6}$
& $\num{1.301983e+00}$
& $\num{4.194304e+12}$
& $\num{0.000000e+00}$ \\
\hline
\hline
\multicolumn{4}{|c|}
{
\cellcolA Stokes regularization ($\gamma=1$)
}
\\
\hline
  $\beta$
& $\min\{(\real(\Lambda_{ii}))_{i=1}^n\}$
& $\max\{(\real(\Lambda_{ii}))_{i=1}^n\}$
& $\max\{(|\imag(\Lambda_{ii})|)_{i=1}^n\}$
\\\hline\rowcolD
  $\num{0}$
& $\num{-3.743833e-13}$
& $\num{2.481540e+01}$
& $\num{5.920341e-06}$
\\ \rowcolC
  $\num{1E-3}$
& $\num{1.000000e-03}$
& $\num{2.560970e+01}$
& $\num{2.594663e-06}$
\\ \rowcolB
  $\num{1E6}$
& $\num{1.294150e+00}$
& $\num{2.048000e+09}$
& $\num{1.425380e-07}$
\\
\hline
\end{tabular}
\end{table}

\ipoint{Observations} The most important observations are that\bipa\item the regularization schemes display a very similar behavior (as judged by the clustering of the eigenvalues as well as the spatial variation of the eigenvectors for $\beta=\num{1E6}$), \item the smoothness of the eigenvectors \emph{decreases} with a \emph{decreasing} regularization parameter $\beta$ and \emph{increasing} eigenvalues (for $\beta=\num{1E-3}$ and $\beta=\num{1E6}$) for both regularization schemes and \item it is less clear how to identify the smooth eigenvectors within the eigenspace of the Stokes regularization scheme\eipa.

The eigenvalues $\Lambda_{ii}$, $i=1,\ldots,8192$, drop rapidly for the unregularized problem, approaching almost machine precision for $i\approx 4000$ (see \figref{f:eigenvalues}). This demonstrates ill-conditioning and ill-posedness. The eigenvalues are bounded away from zero for the regularized problem. Increasing $\beta$ shifts the trend of $|\Lambda_{ii}|$ to larger numbers. The values in~\tabref{t:eigenvalues} confirm that $\D{H}^h\succ 0$ (up to almost machine precision) at the true solution $\vect{v}^{h,\star}$.

Turning to the eigenvector plots, we can see that the first eigenvector displays a delta peak like structure for both regularization schemes, since there is no local coupling of the spatial information. For the regularized problem we can observe a smooth spatial variation for the eigenvectors associated with large eigenvalues for both regularization schemes. The first eigenvector plot is almost constant for $\beta=\num{1E6}$ (bottom row of each block in~\figref{f:eigenvectors}). The structure of the pattern for $\beta=\num{1E6}$ is analogue for both schemes, which indicates similarities in the behavior of both schemes. As the index $i$ increases, the eigenvectors become more oscillatory. We can observe strong differences between the two schemes for a moderate regularization ($\beta=\num{1E-3}$; middle row of each block in \figref{f:eigenvectors}). Also, we can observe a more complex structure for the Stokes regularization scheme for small eigenvalues (i.e.\ it is difficult to identify where the smoothest eigenvectors are located within the eigenspace).

\ipoint{Conclusion} We conclude that the Hessian operator is singular if we do not include a smoothness regularization model for the control variable. For practical values of the regularization parameter, the Hessian behaves as a compact operator; larger eigenvalues are associated with smooth eigenvectors. It is well known that designing a preconditioner for such operators is challenging.

\subsubsection{Convergence Study}
\label{s:res-convergence}

We study the grid convergence of the considered iterative optimization methods on the basis of synthetic registration problems. We use a rigid setting to prevent bias originating from adaptive changes during the computations. That is, the results are computed on a single resolution level. No grid, scale or parameter continuation is applied. The number of the time points is fixed to $n_t = 4\max(\vect{n}_x)$ for all experiments. Further, we use empirically determined values for the regularization parameter $\beta$, namely $\beta \in \{\num{1E-2}, \num{1E-3}\}$. Since we are interested in studying the convergence properties of our method, we consider the relative change of the $\ell^\infty$-norm of the reduced gradient $\vect{g}^h$ as a stopping criterion. This yields a fair comparison between the different optimization methods, as a reduction in the norm of $\vect{g}^h$ directly reflects how well an optimization problem is solved (i.e.\ we exploit that $\vect{g}^h=0$ is a necessary condition for a minimizer). We terminate if the relative change of the $\ell^{\infty}$-norm of $\vect{g}^h$ is at least three orders of magnitude. However, since the Picard method tends to converge slowly for low tolerances with respect to the gradient, we stop if we detect a stagnation in the objective. In particular, we terminate the optimization if the change in the objective in ten consecutive iterations was equal or below $\num{1E-6}$. We solve for a stationary velocity field (i.e.\ $n_c=1$).

\paragraph{$C^\infty$ Registration Problem}
\label{s:res-convergence-cinfty}
\secheadskip

\ipoint{Purpose} We study the numerical behavior for smooth registration problems. We report results for grid convergence and deformation regularity. We compare the Picard, the GN-PCG and the N-PCG method.

\ipoint{Setup} This experiment is based on the \iquote{sinusoidal images} (see \secref{s:data} for more details on the construction of this synthetic registration problem). Therefore, $m_T, m_R \in C^\infty(\Omega)$ and $\vect{v}^\star \in L^2([0,1]; C^\infty(\Omega)^d)$ so that the excellent convergence properties of Fourier spectral methods are expected to pay off. Additionally, it is not problematic to apply the N-PCG method. We report results for different grid sizes $\vect{n}_x = (n_x^1, n_x^2)^\T$, $n_x^i \in \{64,128,256\}$, $i = 1,2$, $n_t = 4\max(\vect{n}_x)$. No pre-smoothing is applied. We use an experimentally determined value of $\beta=\num{1E-3}$ for all experiments. The remainder of the parameters are chosen as stated in the introduction to this section as well as in~\secref{s:parameters}.

\begin{figure}
\centering
\includegraphics[width=\textwidth]
{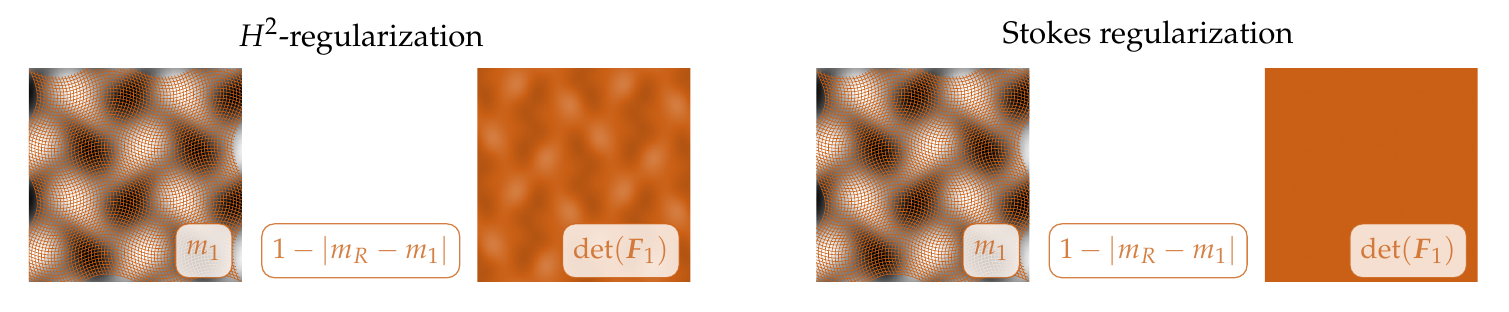}
\caption
{
Qualitative comparison of exemplary registration results of the convergence study reported in \tabref{t:convergence-cinfty}. In particular, we display the results for the N-PCG method for a grid size of $\vect{n}_x = (256,256)^\T$. We refer to \tabref{t:convergence-cinfty} and the text for details on the experimental setup. We report results for plain $H^2$-regularization ($\gamma=0$; images to the left) and the Stokes regularization scheme ($\gamma=1$; $H^1$-regularization; images to the right). We display the deformed template image $m_1$, a pointwise map of the residual differences between $m_R$ and $m_1$ (which appears completely white, as the residual differences are extremely small) as well as a pointwise map of the determinant of the deformation gradient $\det(\mat{F}_1)$ (from left to right as identified by the inset in the images). The values for the $\det(\mat{F}_1)$ are reported in \tabref{t:defgrad-cinfty}. Information on how to interpret these images can be found in \secref{s:map-and-def-grad}.
}
\label{f:convergence-cinfty}
\end{figure}

\ipoint{Results} The grid convergence results are summarized in~\tabref{t:convergence-cinfty}. Values derived from the deformation gradient $\mat{F}_1^h$ are reported in~\tabref{t:defgrad-cinfty}. Exemplary result for the plain $H^2$-regularization ($\gamma=0$) and the Stokes regularization scheme ($\gamma=1$; $H^1$-regularization) are displayed in~\figref{f:convergence-cinfty}. The definitions of the quantitative measures reported in~\tabref{t:convergence-cinfty} and~\tabref{t:defgrad-cinfty} can be found in \tabref{t:quantities-of-interest}.

\begin{table}
\caption
{
Quantitative analysis of the convergence of the Picard, N-PCG and GN-PCG methods. The test problem is the \iquote{sinusoidal images} (see~\secref{s:data} for details on the construction of this synthetic registration problem). We compare convergence results for plain $H^2$-regularization ($\gamma=0$; top block) and the Stokes regularization scheme ($\gamma=1$; $H^1$-regularization; bottom block). We report results for different grid sizes $\vect{n}_x = (n_x^1, n_x^2)^\T$, $n_x^i \in \{64,128,256\}$, $i = 1,2$. We invert for a stationary velocity field (i.e. $n_c = 1$). We terminate the optimization if the relative change of the $\ell^{\infty}$-norm of the reduced gradient $\vect{g}^h$ is at least three orders of magnitude or if the change in $\F{J}^h$ between ten successive iterations is below or equal to $\num{1E-6}$ (i.e.\ the algorithm stagnates). The regularization parameter is empirically set to $\beta=\num{1E-3}$. The number of the (outer) iterations ($k^\star$), the number of the hyperbolic PDE solves ($n_{\rm PDE}$) and the relative change of ($i$) the $L^2$-distance ($\|m_R^h - m_1^h\|_{2,\text{rel}}$), ($ii$) the objective ($\delta\F{J}^h_{\text{rel}}$), and ($iii$) the (reduced) gradient ($\|\vect{g}^h\|_{\infty,\text{rel}}$) and the average number of required line search steps $\bar{\alpha}$ are reported. Note that we introduced a memory for the step size into the Picard method (preconditioned gradient descent) to stabilize the optimization (see \secref{s:parameters} and the description of the results). The definitions for the reported measures are summarized in \tabref{t:quantities-of-interest}.
}
\label{t:convergence-cinfty}
\tabadjust
\begin{tabular}{|l|c|l|l|l|l|l|l|l|}\cline{3-9}
\multicolumn{2}{c|}{}
& $n_x^i$
& $k^\star$
& $n_{\rm PDE}$
& $\|m_R^h - m_1^h\|_{2,\text{rel}}$
& $\delta\F{J}^h_{\text{rel}}$
& $\|\vect{g}^h\|_{\infty,\text{rel}}$
& $\bar{\alpha}$
\\\hline
\parbox[c]{2mm}{\multirow{9}{*}{
\rotatebox[origin=c]{90}{\!\!\!$H^2$-regularization}
}}
& \multirow{3}{*}{Picard}
& 64
& 30
& 420
& \num{4.783097e-03}
& \num{6.051744e-02}
& \num{7.862748e-03}
& \num{1.724138e+00}
\\
&
& 128
& 34
& 414
& \num{4.685544e-03}
& \num{6.043272e-02}
& \num{6.481379e-03}
& \num{1.818182e+00}
\\
&
& 256
& 33
& 414
& \num{4.676883e-03}
& \num{6.033027e-02}
& \num{7.417384e-03}
& \num{1.750000e+00}
\\\cline{2-9}
& \multirow{3}{*}{GN-PCG}
& 64
& 5
& 78
& \num{4.631096e-03}
& \num{6.046260e-02}
& \num{5.566868e-05}
& \num{1.000000e+00}
\\
&
& 128
& 4
& 45
& \num{4.586292e-03}
& \num{6.039090e-02}
& \num{5.603391e-04}
& \num{1.000000e+00}
\\
&
& 256
& 4
& 45
& \num{4.579153e-03}
& \num{6.028773e-02}
& \num{5.103933e-04}
& \num{1.000000e+00}
\\\cline{2-9}
& \multirow{3}{*}{N-PCG}
& 64
& 5
& 63
& \num{4.629180e-03}
& \num{6.046261e-02}
& \num{1.939781e-04}
& \num{1.000000e+00}
\\
&
& 128
& 5
& 61
& \num{4.573073e-03}
& \num{6.039083e-02}
& \num{2.834682e-04}
& \num{1.000000e+00}
\\
&
& 256
& 5
& 61
& \num{4.565849e-03}
& \num{6.028765e-02}
& \num{2.859478e-04}
& \num{1.000000e+00}
\\\hline\hline
\parbox[c]{2mm}{\multirow{9}{*}{
\rotatebox[origin=c]{90}{\!\!\!Stokes regularization}
}}
& \multirow{3}{*}{Picard}
& 64
& 57
& 269
& \num{6.607365e-04}
& \num{1.559613e-02}
& \num{2.554601e-03}
& \num{1.678571e+00}
\\
&
& 128
& 52
& 250
& \num{5.785180e-04}
& \num{1.551730e-02}
& \num{3.547777e-03}
& \num{1.647059e+00}
\\
&
& 256
& 50
& 233
& \num{5.671171e-04}
& \num{1.548508e-02}
& \num{2.310146e-03}
& \num{1.714286e+00}
\\\cline{2-9}
& \multirow{3}{*}{GN-PCG}
& 64
& 5
& 88
& \num{5.861948e-04}
& \num{1.558336e-02}
& \num{7.295983e-05}
& \num{1.000000e+00}
\\
&
& 128
& 5
& 86
& \num{4.869533e-04}
& \num{1.549839e-02}
& \num{8.093116e-05}
& \num{1.000000e+00}
\\
&
& 256
& 5
& 86
& \num{4.864613e-04}
& \num{1.547083e-02}
& \num{8.333276e-05}
& \num{1.000000e+00}
\\\cline{2-9}
& \multirow{3}{*}{N-PCG}
& 64
& 5
& 75
& \num{5.863705e-04}
& \num{1.558337e-02}
& \num{2.596726e-04}
& \num{1.000000e+00}
\\
&
& 128
& 5
& 75
& \num{4.873101e-04}
& \num{1.549840e-02}
& \num{2.694273e-04}
& \num{1.000000e+00}
\\
&
& 256
& 5
& 75
& \num{4.866677e-04}
& \num{1.547084e-02}
& \num{2.514276e-04}
& \num{1.000000e+00}
\\\hline
\end{tabular}
\end{table}

\begin{table}
\caption
{
Obtained values for $\det(\mat{F}_1^h)$ of exemplary registration results of the convergence study reported in \tabref{t:quantities-of-interest} with respect to different iterative optimization methods (Picard, N-PCG and GN-PCG). We refer to \tabref{t:convergence-cinfty} and the text for details on the experimental setup. We report results for plain $H^2$-regularization ($\gamma=0$; top block) and for the Stokes regularization scheme ($\gamma=1$; $H^1$-regularization; bottom block) for a grid size of $\vect{n}_x=(256,256)^\T$.
}
\label{t:defgrad-cinfty}
\tabadjust
\begin{tabular}{|c|c|cccc|}\cline{3-6}
\multicolumn{2}{c|}{}
&$\min(\det(\mat{F}_1^h))$
&$\max(\det(\mat{F}_1^h))$
&$\operatorname{mean}(\det(\mat{F}_1^h))$
&$\operatorname{std}(\det(\mat{F}_1^h))$
\\
\hline
\multirow{3}{*}
{
$H^2$-regularization
}
& Picard
& \num{8.814271e-01}
& \num{1.187278e+00}
& \num{1.002817e+00}
& \num{7.577581e-02}
\\\cline{2-6}
& GN-PCG
& \num{8.830590e-01}
& \num{1.186437e+00}
& \num{1.002804e+00}
& \num{7.559544e-02}
\\\cline{2-6}
& N-PCG
& \num{8.831770e-01}
& \num{1.186366e+00}
& \num{1.002787e+00}
& \num{7.537390e-02}
\\\hline\hline
\multirow{3}{*}
{
Stokes regularization
}
& Picard
& \num{1.000000e+00}
& \num{1.000000e+00}
& \num{1.000000e+00}
& \num{4.547229e-12}
\\\cline{2-6}
& GN-PCG
& \num{1.000000e+00}
& \num{1.000000e+00}
& \num{1.000000e+00}
& \num{4.520745e-12}
\\\cline{2-6}
& N-PCG
& \num{1.000000e+00}
& \num{1.000000e+00}
& \num{1.000000e+00}
& \num{4.524517e-12}
\\\hline
\end{tabular}
\end{table}

\ipoint{Observations} The most important observations are:\bipa\item there are significant differences in computational work between the Picard and the Newton-Krylov methods with the latter being much more efficient, \item the differences between the Newton-Krylov methods are insignificant, \item the rate of convergence is independent of the grid resolution and \item the numerical accuracy is almost at the order of machine precision.\eipa

The registered images are quantitatively (see \tabref{t:convergence-cinfty}) and qualitatively (see \figref{f:convergence-cinfty}) in excellent agreement. For the considered tolerance (reduction of the $\ell^\infty$-norm of the reduced gradient by three orders of magnitude) we can reduce the $L^2$-distance between three (compressible deformation) and four (incompressible deformation) orders of magnitude (see \tabref{t:convergence-cinfty}). The search direction of the Newton-Krylov methods is nicely scaled. No additional line search steps are necessary. We require 1.57 to 1.71 line search steps for the Picard iteration (on average). Note that we prescale the search direction of the Picard method by an additional parameter $\tilde{\alpha}_k$, which is estimated during the computation (see \secref{s:parameters} for details). Otherwise, the number of the line search steps would be seven to eight on average for the Picard method. The Picard method did stagnate during the computations. This is why the gradient has not been reduced by three orders of magnitude for the Picard method. However, it is in general possible to reduce the gradient accordingly. We decided to report only until stagnation for the Picard method as the number of the iterations would significantly increase without making any real progress.

The Newton-Krylov methods display quick convergence. Only five outer iterations are necessary to reduce the gradient by more than four orders of magnitude. The results demonstrate a significant difference in the computational work between first and second-order methods for the considered tolerance.

The reconstruction quality improves by approximately one order of magnitude when switching from plain $H^2$-regularization ($\gamma=0$) to a Stokes regularization scheme ($\gamma=1$; $H^1$-regularization), as judged by the relative change in the $L^2$-distance. This is expected, since the synthetic problem has been created under the assumption of mass conservation (i.e.\ $\idiv \vect{v}^\star = 0$). Secondly, we expect a smaller contribution of the $H^1$-regularization model on the solution for the same values of $\beta$.

From a theoretical point of view, we expect N-PCG to outperform GN-PCG (quadratic vs.\ super-linear convergence). The reported results demonstrate an almost identical performance. This is due to the fact that we can drive the residual almost to zero, such that we can recover fast local convergence for the GN-PCG method (see~\secref{s:gauss-newton} for details).

The Picard method converges faster for the Stokes regularization scheme. However, the differences between the Picard and Newton-Krylov methods are still significant with an approximately four-fold difference in $n_{\rm PDE}$. For the Newton-Krylov methods, we can globally observe a slight increase in the number of inner iterations when switching from plain $H^2$-regularization to the Stokes regularization scheme. These differences have to be attributed to a varying relative change in the reduced gradient\footnote{Note that the tolerance of the Krylov-subspace method and therefore the number of inner iterations depends on the gradient (see~\eqref{e:termination-criterion-ksm}).}.

The results reported in \tabref{t:defgrad-cinfty} demonstrate an excellent numerical accuracy for the mass conservation for all numerical schemes. The error in the determinant of the deformation gradient is $\bigO(\num{1E-12})$, i.e.\ we achieve an accuracy that is almost at the order of machine precision for a grid resolution of $\vect{n}_x=(256,256)^\T$ and $n_t = 1024$.

\ipoint{Conclusion} We conclude that we can interchangeably use the Newton-Krylov methods. Therefore, given that N-PCG is more sensitive to noise and discontinuities in the data, we will exclusively consider GN-PCG for the remainder of the experiments. Also, if we require an inversion with high accuracy, the Newton-Krylov methods clearly outperform the Picard method (i.e.\ the preconditioned gradient descent).

\paragraph{Images with Sharp Features}
\label{s:res-convergence-ut}
\secheadskip

\ipoint{Purpose} We study the grid convergence and deformation regularity for an image with sharp features. We compare the Picard and the GN-PCG method.

\ipoint{Setup} We consider the \iquote{UT images} (see \secref{s:data} for details on the construction of this synthetic registration problem). We report results for experimentally determined values of $\beta\in\{\num{1E-2},\num{1E-3}\}$ with respect to different grid resolution levels $\vect{n}_x = (n_x^1, n_x^2)^\T$, $n_x^i \in \{64,128,256\}$, $i = 1,2$, $n_t = 4\max(\vect{n}_x)$. The remainder of the parameters are chosen as stated in the introduction of this section and in~\secref{s:parameters}. Both, plain $H^2$-regularization ($\gamma=0$) as well as the Stokes regularization scheme ($\gamma=1$; $H^1$-regularization) are considered.

For images of size $\vect{n}_x=(256,256)^\T$ and a Stokes regularization scheme, we observed difficulties in the inversion (only the number of outer and inner iterations increased; the algorithm still converges to the same solution), due to a strong forcing (i.e.\ the sharp features pushed the solver at an early stage to a solution that was far away from the final minimizer). We increased the smoothing by a factor of two as a remedy. This is not an issue for the practical application of our algorithm, as our framework features a method for performing a scale continuation as well as a continuation in the regularization parameter. Therefore, the user does not have to decide on $\sigma$ nor on $\beta$. In addition to that, we currently investigate adaptive approaches to automatically detect insufficient smoothness during the course of the optimization to prevent a deterioration in the convergence behavior.

The remainder of the parameters are chosen as stated in the introduction of this section as well as in~\secref{s:parameters}.

\begin{figure}
\centering
\includegraphics[width=\textwidth]
{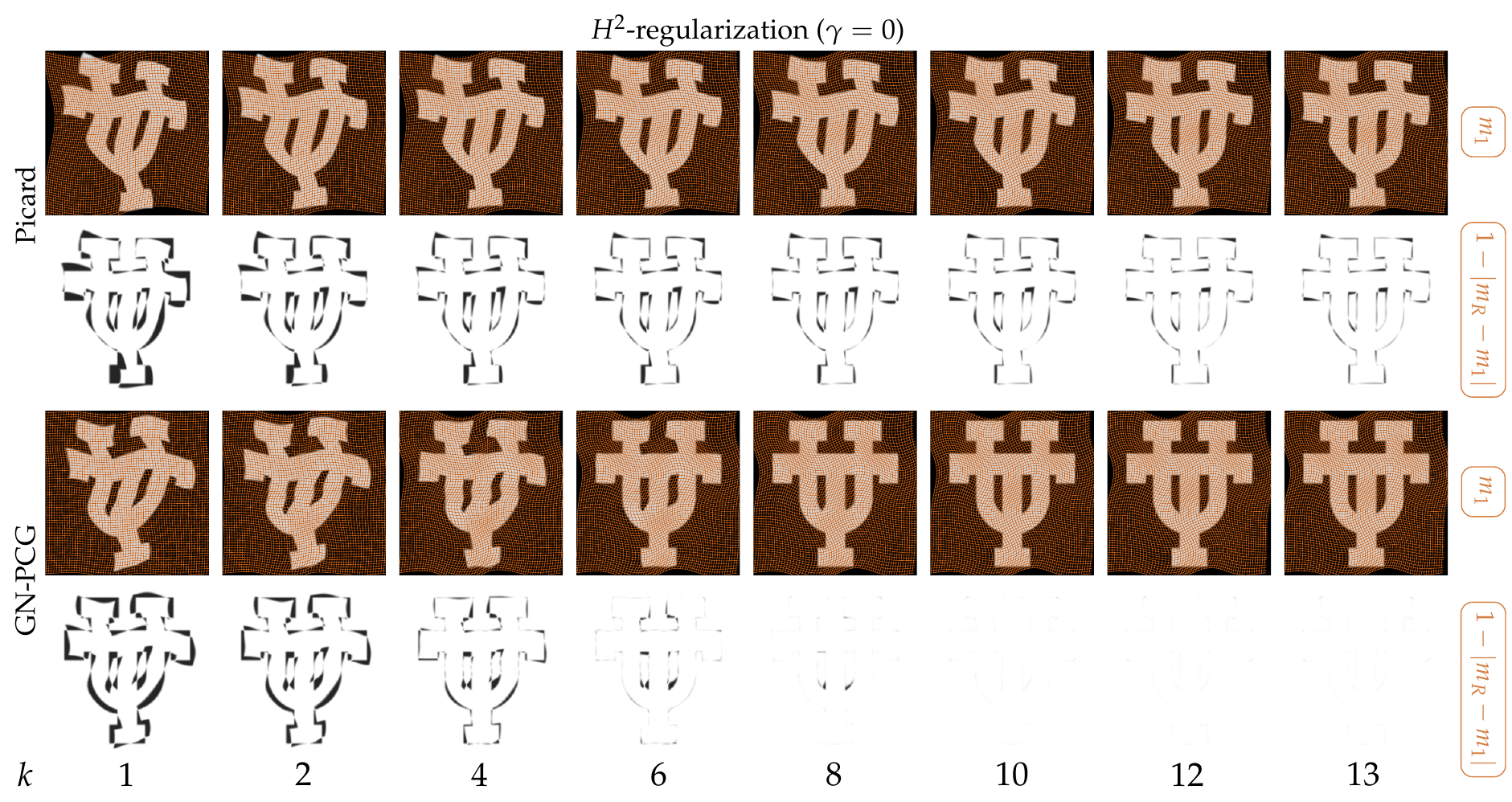}
\caption
{
Illustration of the course of the optimization for the Picard (top block) and the GN-PCG (bottom block) method with respect to the (outer) iteration index $k$ for exemplary results of the convergence study reported in \tabref{t:convergence-ut}. We refer to \tabref{t:convergence-ut} and the text for details on the experimental setup. We report results for plain $H^2$-regularization ($\gamma=0$) with an empirically chosen regularization parameter of $\beta = \num{1E-3}$ for images of grid size $\vect{n}_x=(256,256)^\T$. We report results until convergence of the GN-PCG method ($k^\star=13$). We display the deformed template $m_1$ (top row) and a map of the pointwise difference between $m_R$ and $m_1$ (bottom row) for both iterative optimization methods (as identified by the inset on the right of the images). Information on how to interpret these images can be found in~\secref{s:map-and-def-grad}.
}
\label{f:convergence-ut-iter}
\end{figure}

\begin{figure}
\centering
\includegraphics[width=0.95\textwidth]
{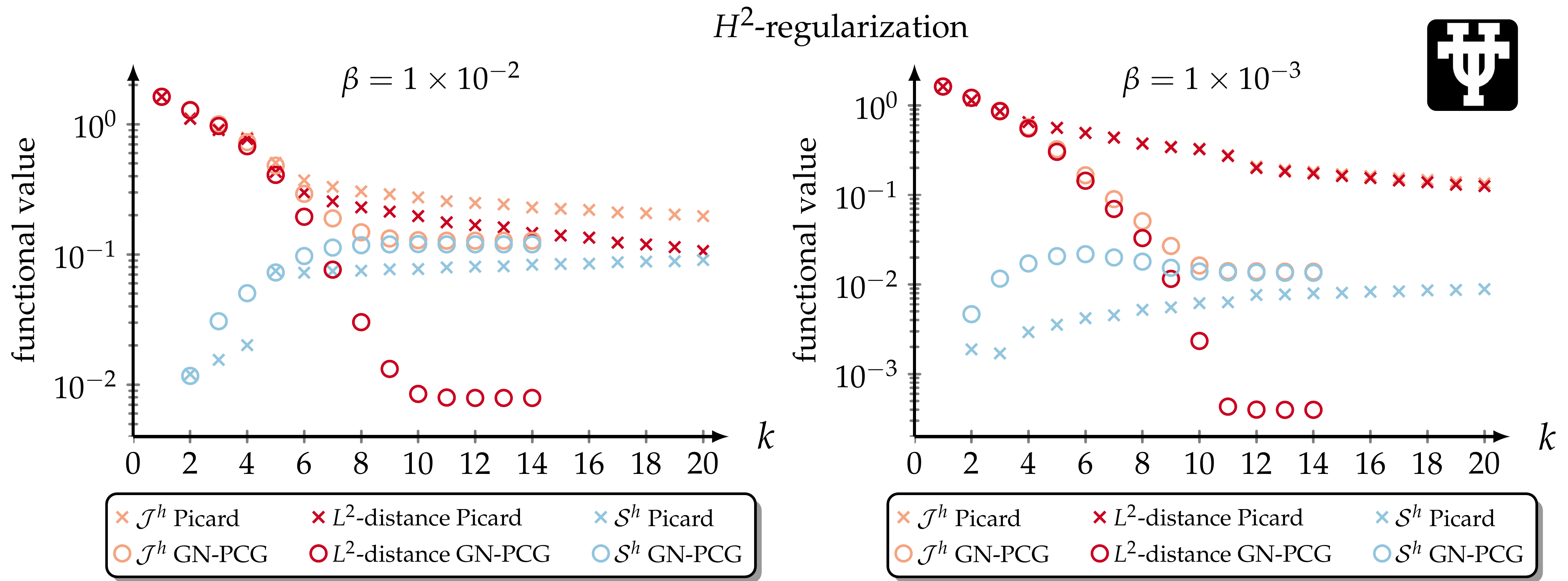}
\\\bigskip
\includegraphics[width=0.95\textwidth]
{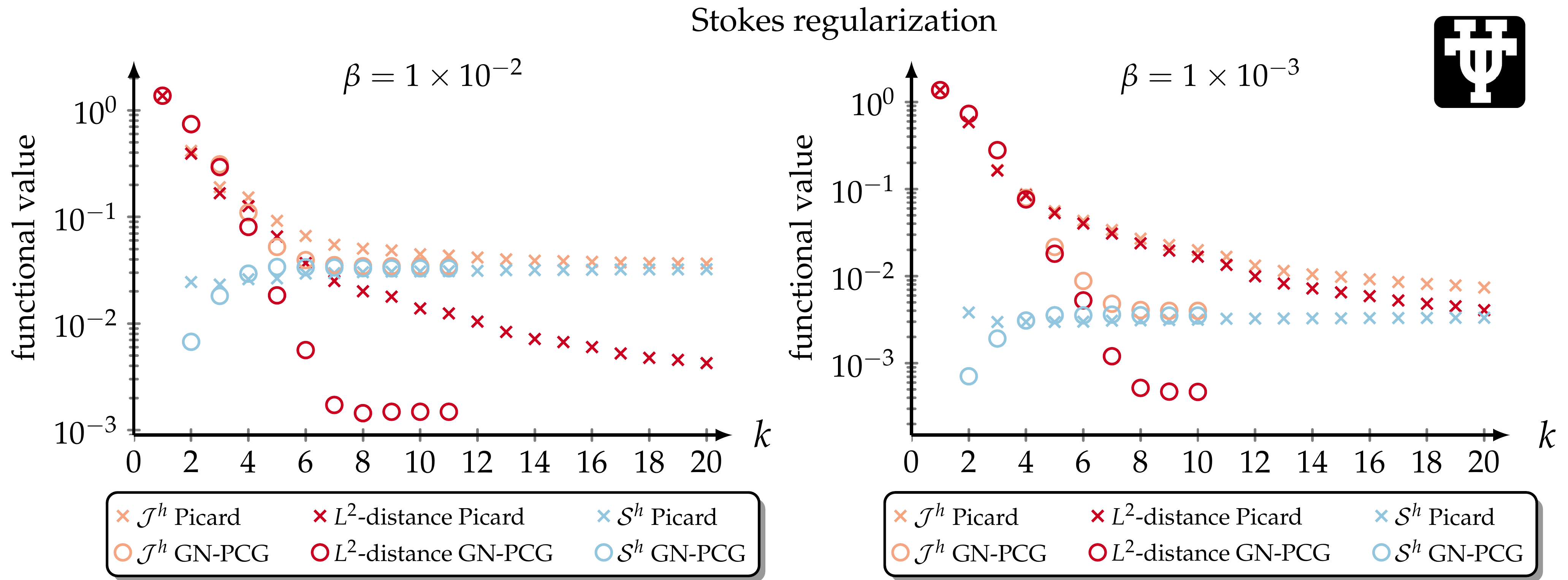}
\caption
{
Trend of the objective $\F{J}^h$, the $L^2$-distance and the regularization model $\F{S}^h$ (logarithmic scale) for the Picard and the GN-PCG method with respect to the (outer) iteration index $k$ for exemplary results of the convergence study reported in \tabref{t:convergence-ut}. We refer to \tabref{t:convergence-ut} and the text for details on the experimental setup. The trend of the functionals is plotted for different (empirically determined) choices of $\beta$ (left column: $\beta=\num{1E-2}$; right column: $\beta=\num{1E-3}$) and a grid size of $\vect{n}_x = (256,256)^\T$. We report results for plain $H^2$-regularization ($\gamma=0$; top row) and the Stokes regularization scheme ($\gamma=1$; $H^1$-regularization; bottom row).
}
\label{f:convergence-ut-f}
\end{figure}

\ipoint{Results} \tabref{t:convergence-ut} summarizes the results of the convergence study. We illustrate intermediate results with respect to the first 13 (outer) iterations $k$ in~\figref{f:convergence-ut-iter} (plain $H^2$-regularization; $\gamma=0$). We report the trend of the individual building blocks of $\F{J}^h$ (contribution of the $L^2$-distance and the regularization model $\F{S}^h$) in~\figref{f:convergence-ut-f}. We report measures of deformation regularity in~\tabref{t:defgrad-ut}.

\begin{table}
\caption
{
Quantitative analysis of the convergence for the Picard and the GN-PCG method. The test problem is the \iquote{UT images} (see~\secref{s:data} for more details on the construction of this synthetic registration problem). We compare convergence results for plain $H^2$-regularization ($\gamma=0$; top block) and the Stokes regularization scheme ($\gamma=1$; bottom block; $H^1$-regularization) for empirically chosen regularization parameters $\beta\in\{\num{1E-2}, \num{1E-3}\}$. We report results for different grid sizes $\vect{n}_x = (n_x^1,n_x^2)^\T$, $n_x^i\in\{64,128,256\}$, $i=1,2$, $n_t = 4\max(\vect{n}_x)$. We invert for a stationary velocity field (i.e. $n_c=1$). We terminate the optimization if the relative change of the $\ell^{\infty}$-norm of the reduced gradient $\vect{g}^h$ is larger than or equal to three orders of magnitude or if the change in $\F{J}^h$ between ten successive iterations is below or equal to $\num{1E-6}$ (i.e.\ the algorithm stagnates). We report the number of the (outer) iterations ($k^\star$), the number of the hyperbolic PDE solves ($n_{\rm PDE}$) and the relative change of ($i$) the $L^2$-distance ($\|m_R^h - m_1^h\|_{2,\text{rel}}$), ($ii$) the objective ($\delta\F{J}^h_{\text{rel}}$), and ($iii$) the (reduced) gradient ($\|\vect{g}^h\|_{\infty,\text{rel}}$) as well as the average number of line search steps $\bar{\alpha}$. Note that we introduced a memory for the step size into the Picard method to stabilize the optimization (see \secref{s:parameters} and the description of the results). The definitions for the reported measures can be found in \tabref{t:quantities-of-interest}. This study directly relates to the results for the smooth registration problem (see \secref{s:res-convergence-cinfty}, in particular \tabref{t:convergence-cinfty}).
}
\label{t:convergence-ut}
\tabadjust
\begin{tabular}{|l|l|c|l|l|l|l|l|l|l|}\cline{4-10}
\multicolumn{3}{c|}{}
& $n_x^i$
& $k^\star$
& $n_{\rm PDE}$
& $\|m_R^h - m_1^h\|_{2,\text{rel}}$
& $\delta\F{J}^h_{\text{rel}}$
& $\|\vect{g}^h\|_{\infty,\text{rel}}$
& $\bar{\alpha}$
\\
\hline
\parbox[c]{2mm}{\multirow{12}{*}{
\rotatebox[origin=c]{90}{\!\!\!$H^2$-regularization}
}}
&  \multirow{6}{*}{\rotatebox[origin=c]{90}{\!\!\!$\beta=\num{1E-2}$}}
&  \multirow{3}{*}{Picard}
& 64
& 130
& 752
& \num{1.987920e-02}
& \num{1.183112e-01}
& \num{2.060494e-02}
& \num{1.682171e+00}
\\
&&
& 128
& 231
& 1589
& \num{8.287684e-03}
& \num{8.473031e-02}
& \num{4.026569e-02}
& \num{1.717391e+00}
\\
&&
& 256
& 388
& 3022
& \num{5.052726e-03}
& \num{7.920616e-02}
& \num{9.856513e-02}
& \num{1.682171e+00}
\\\cline{3-10}
&& \multirow{3}{*}{GN-PCG}
& 64
& 7
& 282
& \num{1.935631e-02}
& \num{1.159537e-01}
& \num{3.296619e-04}
& \num{1.000000e+00}
\\
&&
& 128
& 9
& 450
& \num{8.088135e-03}
& \num{8.368557e-02}
& \num{6.066507e-04}
& \num{1.000000e+00}
\\
&&
& 256
& 13
& 789
& \num{4.871005e-03}
& \num{7.848235e-02}
& \num{7.283920e-04}
& \num{1.000000e+00}
\\
\cline{2-10}
&  \multirow{6}{*}{\rotatebox[origin=c]{90}{\!\!\!$\beta=\num{1E-3}$}}
&  \multirow{3}{*}{Picard}
& 64
& 339
& 4410
& \num{1.842256e-03}
& \num{1.538902e-02}
& \num{2.455415e-02}
& \num{1.562130e+00}
\\
&&
& 128
& 466
& 9671
& \num{7.190131e-04}
& \num{9.843405e-03}
& \num{5.978488e-02}
& \num{1.643011e+00}
\\
&&
& 256
& 632
& 8690
& \num{5.346946e-04}
& \num{8.805996e-03}
& \num{1.080744e-01}
& \num{1.640254e+00}
\\\cline{3-10}
&& \multirow{3}{*}{GN-PCG}
& 64
& 7
& 670
& \num{1.440760e-03}
& \num{1.496298e-02}
& \num{2.586141e-04}
& \num{1.000000e+00}
\\
&&
& 128
& 8
& 929
& \num{4.468053e-04}
& \num{9.595323e-03}
& \num{4.761574e-04}
& \num{1.000000e+00}
\\
&&
& 256
& 13
& 1744
& \num{2.449139e-04}
& \num{8.553653e-03}
& \num{4.043572e-04}
& \num{1.000000e+00}
\\\hline\hline
\parbox[c]{2mm}{\multirow{12}{*}{
\rotatebox[origin=c]{90}{\!\!\!Stokes regularization}
}}
&  \multirow{6}{*}{\rotatebox[origin=c]{90}{\!\!\!$\beta=\num{1E-2}$}}
&  \multirow{3}{*}{Picard}
& 64
& 92
& 448
& \num{2.732051e-03}
& \num{3.623628e-02}
& \num{7.537590e-03}
& \num{1.637363e+00}
\\
&&
& 128
& 136
& 958
& \num{8.924824e-04}
& \num{2.384272e-02}
& \num{1.630243e-02}
& \num{1.718519e+00}
\\
&&
& 256
& 143
& 1004
& \num{1.092540e-03}
& \num{2.529728e-02}
& \num{1.745555e-02}
& \num{1.683099e+00}
\\\cline{3-10}
&& \multirow{3}{*}{GN-PCG}
& 64
& 7
& 313
& \num{2.722711e-03}
& \num{3.592536e-02}
& \num{4.697776e-04}
& \num{1.000000e+00}
\\
&&
& 128
& 9
& 437
& \num{8.880099e-04}
& \num{2.371641e-02}
& \num{5.587307e-04}
& \num{1.000000e+00}
\\
&&
& 256
& 10
& 514
& \num{1.086716e-03}
& \num{2.515382e-02}
& \num{5.491655e-04}
& \num{1.000000e+00}
\\\cline{2-10}
&  \multirow{6}{*}{\rotatebox[origin=c]{90}{\!\!\!$\beta=\num{1E-3}$}}
&  \multirow{3}{*}{Picard}%
& 64
& 216
& 2823
& \num{9.990236e-04}
& \num{4.693308e-03}
& \num{1.336220e-02}
& \num{1.558140e+00}
\\
&&
& 128
& 179
& 5465
& \num{3.565392e-04}
& \num{2.765639e-03}
& \num{1.785651e-02}
& \num{1.685393e+00}
\\
&&
& 256
& 175
& 5569
& \num{5.271465e-04}
& \num{3.069135e-03}
& \num{2.224109e-02}
& \num{1.678161e+00}
\\\cline{3-10}
&& \multirow{3}{*}{GN-PCG}
& 64
& 8
& 1162
& \num{7.573473e-04}
& \num{4.546110e-03}
& \num{8.227993e-04}
& \num{1.000000e+00}
\\
&&
& 128
& 9
& 1069
& \num{2.145082e-04}
& \num{2.648537e-03}
& \num{5.667271e-04}
& \num{1.000000e+00}
\\
&&
& 256
& 9
& 769
& \num{3.425650e-04}
& \num{2.923476e-03}
& \num{6.435597e-04}
& \num{1.000000e+00}
\\\hline
\end{tabular}
\end{table}

\begin{table}
\caption
{
Values for the determinant of the deformation gradient $\det(\mat{F}_1^h)$ for exemplary results of the convergence study reported in \tabref{t:convergence-ut}. We report results for the Picard and the GN-PCG method. We refer to \tabref{t:convergence-ut} and the text for details on the experimental setup. We report results for the Stokes regularization scheme ($\gamma=1$; $H^1$-regularization). The regularization parameter is set to $\beta=\num{1E-3}$. The grid size is $\vect{n}_x = (256,256)^\T$. These results directly relate to those reported for the smooth registration problem (see \secref{s:res-convergence-cinfty}, in particular \tabref{t:defgrad-cinfty}).
}
\label{t:defgrad-ut}
\tabadjust
\begin{tabular}{|c|cccc|}\cline{2-5}
\multicolumn{1}{c|}{}
& $\min(\det(\mat{F}_1^h))$
& $\max(\det(\mat{F}_1^h))$
& $\operatorname{mean}(\det(\mat{F}_1^h))$
& $\operatorname{std}(\det(\mat{F}_1^h))$
\\
\hline
Picard
& \num{9.998504e-01}
& \num{1.000135e+00}
& \num{1.000000e+00}
& \num{2.025202e-05}
\\\hline
GN-PCG
& \num{9.998338e-01}
& \num{1.000162e+00}
& \num{1.000000e+00}
& \num{2.285846e-05}
\\\hline
\end{tabular}
\end{table}

\ipoint{Observations} The most important observations are\bipa\item the GN-PCG method displays a quicker convergence than the Picard method, \item we cannot achieve the same inversion accuracy with the Picard method as compared to the GN-PCG method, and \item the number of the (inner and outer) iterations increases and is no longer independent of the resolution level\eipa.

The rate of convergence decreases compared to the results reported in the former section (see~\tabref{t:convergence-cinfty}). Overall, we require more outer and inner iterations to solve the registration problem.

The residual differences between $m_R$ and $m_1$ clearly depend on the choice of $\beta$ (see~\tabref{t:convergence-cinfty} and~\figref{f:convergence-ut-f}). We achieve a similar reduction in the $L^2$-distance for both the Picard and the GN-PCG method (two to four orders of magnitude). The residual differences are less pronounced when switching from plain $H^2$-regularization to the Stokes regularization scheme as compared to the results reported in \secref{s:res-convergence-cinfty}.

We cannot guarantee that it is possible to reduce the gradient by three orders of magnitude if we use the Picard method. Even if we do not include a condition to terminate if we observe stagnation (i.e.\ the change in $\F{J}^h$ is below or equal to $\num{1E-6}$ for 10 consecutive iterations), it is for some of the experiments not possible to reduce the gradient by three orders of magnitude as the changes of the objective hit our numerical accuracy (which causes the line search to fail). We do not observe this issue when considering the GN-PCG method. Further, there are significant differences in terms of the computational work. If we do not account for the stagnation of the Picard method we have observed a number of hyperbolic PDE solves that is well above $\bigO(\num{1E4})$. Clearly, in a practical application we terminate the Picard method at an earlier stage, as we no longer make significant progress. However, in this part of the study we are interested in the convergence properties. This experiment demonstrates that we cannot guarantee a high inversion accuracy (i.e.\ a significant reduction in the gradient) when turning to first-order methods. Note that we have stabilized the Picard method by introducing an additional scaling parameter for the search direction that prevents additional line search steps (see \secref{s:parameters}). If we neglect this scaling, we observe seven to nine line search steps on average  (results not included in this study) for the considered problem; also, the optimization fails at an early stage. The search direction obtained via the GN-PCG method is nicely scaled; no additional line search steps are necessary.

The trend of the $\F{J}^h$, the $L^2$-distance and $\F{S}^h$ in \figref{f:convergence-ut-f} confirm these observations. The plots in \figref{f:convergence-ut-f} illustrate that the Picard and the GN-PCG method perform similarly during the first few outer iterations. However, after about four outer iterations the differences between the methods manifest, in particular with respect to the reduction of the $L^2$-distance. This observation confirms standard numerical optimization theory on convergence properties of the Picard and the inexact Newton-Krylov methods.

Focusing on the GN-PCG method we can observe that the number of outer iterations is almost constant across different grid sizes. However, the effectiveness of the spectral preconditioner decreases with an increasing grid size as well as with a reduction of the regularization parameter (as judged by an increase in the number of inner iterations). This demonstrates that the preconditioner is not optimal. A similar behavior can be observed for the Picard method\footnote{Note that the Picard method is a gradient descent scheme in the function space induced by $\fs{W}$. We can interpret the inverse of $\F{A}^h$ as a preconditioner acting on the body force $\vect{f}$. This operator is exactly the spectral preconditioner we use for the Newton-Krylov methods, which explains the similar behavior.}.

The numerical accuracy of the incompressibility constraint deteriorates (slightly but not significantly) as compared to the results reported in the former section. In particular, we obtain a numerical accuracy of $\bigO(\num{1E-5})$ for the GN-PCG method (see~\tabref{t:defgrad-ut}).

\ipoint{Conclusion} We conclude that the GN-PCG is less sensitive, provides a better inversion accuracy and overall displays quicker convergence if a high accuracy of the inversion is required and, therefore, is to be preferred.

\subsubsection{Number of Unknowns in Time}
\label{s:res-expansion}
\secheadskip

\ipoint{Purpose} It is not immediately evident how the number of the coefficient fields $\vect{v}_l^h : \ns{R}^d\rightarrow\ns{R}^d$, $l=1,\ldots,n_c$, affects the registration quality. We study the effects of varying $n_c$ on the reconstruction quality and the rate of convergence. We also provide advice on how to decide on $n_c$.

\ipoint{Setup} We report results for registration problems of varying complexity. The analysis is limited to the GN-PCG method. We consider the \iquote{hand images} ($\vect{n}_x=(128,128)^\T$) and the \iquote{brain images} ($\vect{n}_x=(200,200)^\T$). The number of the time steps is fixed to $n_t=2\max(\vect{n}_x)$. The regularization parameter is empirically set to $\beta=\num{1E-3}$ and $\beta=\num{2E-2}$, respectively. We consider the full set of stopping conditions in \eqref{e:stopping-criteria} with $\tau_{\F{J}}=\num{1E-3}$, as we no longer compare different methods. The remainder of the parameters are set as stated in~\secref{s:parameters}.

One possibility to estimate an adequate number of coefficients for the registration of unseen images $m_R$ and $m_T$ is to compute the relative spectral power (see \tabref{t:quantities-of-interest}) of an individual coefficient field $\vect{v}_l^h$ for different choices of $n_c$. If only a small number of coefficients is necessary to recover the deformation, this energy should decrease rapidly with an increasing $l$. The problem is stationary for $n_c=1$.

\begin{figure}
\centering
\includegraphics[width=\textwidth]
{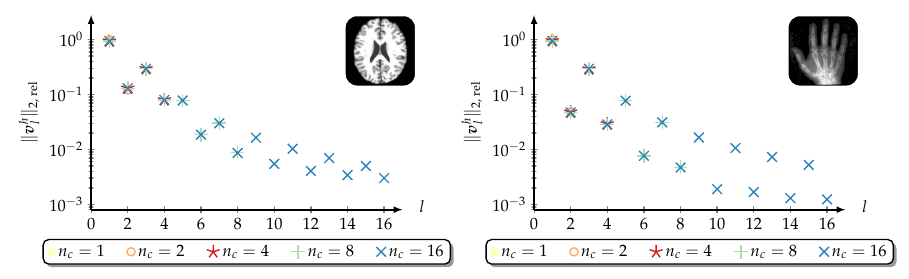}
\caption
{
Relative power spectrum of the individual coefficient fields $\vect{v}_l^h:\ns{R}^d\rightarrow \ns{R}^d$, $l=1,\ldots,n_c$, for different choices of $n_c$ used to solve the considered registration problems. We report results for plain $H^2$-regularization. The reported results correspond to \tabref{t:conv-num-coeff}. We refer to \tabref{t:conv-num-coeff} and the text for details on the experimental setup. We report exemplary results for the the \iquote{brain images} ($\vect{n}_x=(200,200)^\T$; left) and the \iquote{hand images} ($\vect{n}_x=(128,128)^\T$; right). We choose $n_c$ to be in $\{1,2,4,8,16\}$ as indicated in the legend of each plot. The definition of the relative $\ell^2$-norm (relative power spectrum) of $\vect{v}^h_l$ can be found in \tabref{t:quantities-of-interest}.
}
\label{f:energy-num-coeff}
\end{figure}

\begin{figure}
\centering
\includegraphics[width=\textwidth]
{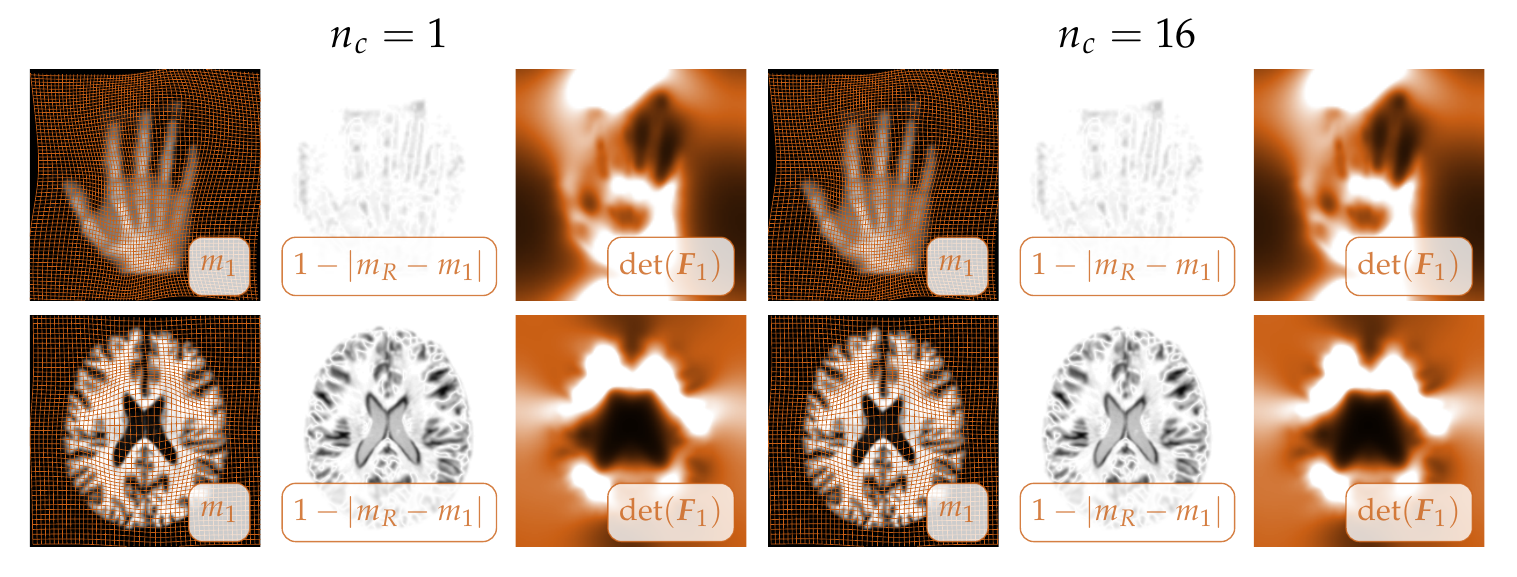}
\caption
{
Qualitative comparison of exemplary registration results for $n_c=1$ (images to the left) and $n_c=16$ (images to the right) of the study reported in \tabref{t:conv-num-coeff}. We refer to \tabref{t:conv-num-coeff} and the text for details on the experimental setup. We report results for the \iquote{hand images} ($\vect{n}_x = (128,128)^\T$) and the \iquote{brain images} ($\vect{n}_x=(200,200)^\T$) for plain $H^2$-regularization. We display (for each experiment) the deformed template image $m_1$, a pointwise map of the absolute difference between $m_R$ and $m_1$ and a map of the determinant of the deformation gradient $\det(\mat{F}_1)$ (from left to right as identified by the inset in the images). Information on how to interpret these images can be found in \secref{s:map-and-def-grad}.
}
\label{f:im-coll-num-coeff}
\end{figure}

\ipoint{Results} The trend of the relative $\ell^2$-norm (i.e.\ the spectral power) of an individual coefficient field $\vect{v}_l^h$ for different choices of $n_c \in \{1,2,4,8,16\}$ is plotted in~\figref{f:energy-num-coeff}. A qualitative comparison of the registration results for different choices of $n_c$ can be found in~\figref{f:im-coll-num-coeff}. Convergence results are reported in~\tabref{t:conv-num-coeff}.

\begin{table}
\caption
{
Comparison of the inversion results for the GN-PCG method for a varying number of the spatial coefficient fields $\vect{v}_l^h:\ns{R}^d\rightarrow \ns{R}^d$, $l = 1,\ldots,n_c$ (i.e.\ we change the number of the unknowns in time). We report results for plain $H^2$-regularization. We consider the \iquote{hand images} ($\vect{n}_x = (128,128)^\T$; top block) and the \iquote{brain images} ($\vect{n}_x = (200,200)^\T$; bottom block). We consider the full set of stopping conditions in \eqref{e:stopping-criteria} with $\tau_{\F{J}}=\num{1E-3}$, as we no longer study grid convergence and/or compare different optimization methods. We report the number of the (outer) iterations ($k^\star$), the number of the hyperbolic PDE solves ($n_{\rm PDE}$) and the relative change of ($i$) the $L^2$-distance ($\|m_R^h - m_1^h\|_{2,\text{rel}}$), ($ii$) the objective ($\delta\F{J}^h_{\text{rel}}$), and ($iii$) the (reduced) gradient ($\|\vect{g}^h\|_{\infty,\text{rel}}$), as well as the minimal and maximal values of the determinant of the deformation gradient. The definitions of these measures can be found in \tabref{t:quantities-of-interest}. The number of the coefficient fields $n_c$ used to solve the individual registration problems is chosen to be in $\{1,2,4,8,16\}$.
}
\label{t:conv-num-coeff}
\tabadjust
\begin{tabular}{|c|l|l|l|l|l|l|l|l|}\cline{2-9}
\multicolumn{1}{c|}{}
& $n_c$
& $k$
& $n_{\rm PDE}$
& $\|m_R^h - m_1^h\|_{2,\text{rel}}$
& $\delta\F{J}^h_{\text{rel}}$
& $\|\vect{g}^h\|_{\infty,\text{rel}}$
& $\min(\det(\mat{F}^h_1))$
& $\max(\det(\mat{F}^h_1))$
\\\hline
\parbox[t]{2mm}
{
\multirow{5}{*}{\rotatebox[origin=c]{90}{\iquote{hand images}}}
}
& 1
& 7
& 279
& \num{6.653706e-02}
& \num{8.523025e-02}
& \num{2.269248e-02}
& \num{2.140051e-01}
& \num{6.595905e+00}
\\
& 2
& 7
& 279
& \num{6.595351e-02}
& \num{8.476257e-02}
& \num{2.294904e-02}
& \num{2.151936e-01}
& \num{6.437201e+00}
\\
& 4
& 7
& 283
& \num{6.414751e-02}
& \num{8.244427e-02}
& \num{2.289360e-02}
& \num{2.067973e-01}
& \num{6.489354e+00}
\\
& 8
& 7
& 277
& \num{6.406930e-02}
& \num{8.241078e-02}
& \num{2.509605e-02}
& \num{2.075160e-01}
& \num{6.455157e+00}
\\
& 16
& 7
& 281
& \num{6.407123e-02}
& \num{8.240913e-02}
& \num{2.561400e-02}
& \num{2.079336e-01}
& \num{6.446671e+00}
\\\hline\hline
\parbox[t]{2mm}
{
\multirow{5}{*}{\rotatebox[origin=c]{90}{\iquote{brain images}}}
}
& 1
& 20
& 669
& \num{5.492716e-01}
& \num{6.675573e-01}
& \num{3.709919e-02}
& \num{4.930483e-02}
& \num{6.465830e+00}
\\
& 2
& 20
& 667
& \num{5.466437e-01}
& \num{6.658711e-01}
& \num{3.720185e-02}
& \num{4.959731e-02}
& \num{6.447210e+00}
\\
& 4
& 21
& 710
& \num{5.307371e-01}
& \num{6.512849e-01}
& \num{3.656659e-02}
& \num{4.551107e-02}
& \num{7.327586e+00}
\\
& 8
& 21
& 710
& \num{5.300295e-01}
& \num{6.506183e-01}
& \num{3.523743e-02}
& \num{4.533632e-02}
& \num{7.369147e+00}
\\
& 16
& 21
& 708
& \num{5.299740e-01}
& \num{6.505564e-01}
& \num{3.539308e-02}
& \num{4.531570e-02}
& \num{7.372251e+00}
\\\hline
\end{tabular}
\end{table}

\ipoint{Observations}: The most important observation is that we obtain the same results for stationary as well as time varying velocity fields for two-image registration problems. Qualitatively, we cannot observe any differences for a varying number of coefficient fields (see~\figref{f:im-coll-num-coeff}). This observation is confirmed by the values for the relative reduction in the $L^2$-distance in~\tabref{t:conv-num-coeff}. Increasing the number of the coefficients slightly reduces the $L^2$-distance. These differences, however, are practically insignificant. In particular, we (on average) observe a relative change in the $L^2$-distance of $\num{6.495600e-02} \pm \num{1.195483e-03}$ (\iquote{hand images}, plain $H^2$-regularization) and $\num{5.373200e-01} \pm \num{9.754845e-03}$ (\iquote{brain images}, plain $H^2$-regularization). Also, we obtain identical deformation patterns as judged by careful visual inspection (see \figref{f:im-coll-num-coeff}) and the variations in the determinant of the deformation gradient. We obtain identical results for the \iquote{UT images} (for plain $H^2$-regularization and the Stokes regularization scheme; results are not included in this study).

Turning to the required work load, we observe that the differences are also insignificant. The number of outer iterations is almost constant; just the number of inner iterations varies. In particular, we require 7 outer iteration with (on average) $\approx$$280$ inner iterations (\iquote{hand images}; plain $H^2$-regularizatoin; $\gamma=0$) and 20-21 outer iterations with (on average) $\approx$$693$ inner iterations (\iquote{brain images}; plain $H^2$-regularization; $\gamma=0$). However, we have to keep in mind that each application of the reduced Hessian is slightly more expensive and we require more memory as $n_c$ increases. That is, we have to store more coefficient fields $\vect{v}_l^h$ (to all of which the regularization operator has to be applied). The cost of the forward and adjoint solves (which is the key bottleneck), however, is (almost) the same, since we expand $\vect{v}^h$ (note that this expansion is not necessary for $n_c=1$; see~\secref{s:expansion}).

The power spectrum of the coefficient fields drops quickly (see~\figref{f:energy-num-coeff}). This also indicates that only a small number of coefficients is required to obtain an excellent agreement between the images. However, we expect the differences to manifest, when registering time series of images (multiple time frames). Here, we might benefit from being able to invert for a time varying velocity field.

\ipoint{Conclusion} We conclude that it is sufficient to use stationary velocity fields for two-image registration problems.

\subsection{Parameter Continuation to Estimate $\beta$}
\label{s:res-parameter-continuation}
\secheadskip

\ipoint{Purpose} We study the stability and accuracy of the designed parameter continuation method (see \secref{s:parameters}) and the associated control over the properties of the mapping. That is, we study how the quantities of interest (determinant of the deformation gradient and $L^2$-distance) behave during the course of the parameter continuation and how close we actually approach the given bounds.

\ipoint{Setup} The registration problems are solved on images with a grid size of $\vect{n}_x = (512,512)^\T$. The number of the time points is adapted as required by monitoring the CFL condition (see \secref{s:discretization}). We use the full set of stopping conditions (see \secref{s:stopping-criteria}) with a tolerance of $\tau_{\F{J}} = \num{1E-3}$. We consider the \iquote{hand images} and the \iquote{brain images} (see \figref{f:regprob}). We invert for a stationary velocity field (i.e.\ $n_c=1$). In case we consider a plain $H^1$- and $H^2$-regularization (smoothness regularization; $\gamma=0$), we set the lower bound on $\det(\mat{F}_1^h)$ to \num{1E-1} (\iquote{hand images}) and \num{5E-2} (\iquote{brain images}), respectively. For the case of the Stokes regularization ($\gamma=1$; $H^1$-regularization) we set the bound on the grid angle to $\epsilon_{\theta} = \pi/16$ (11.25\textdegree). The remainder of the parameters are set as described in~\secref{s:parameters}.

\begin{figure}
\centering
\includegraphics[width=\textwidth]
{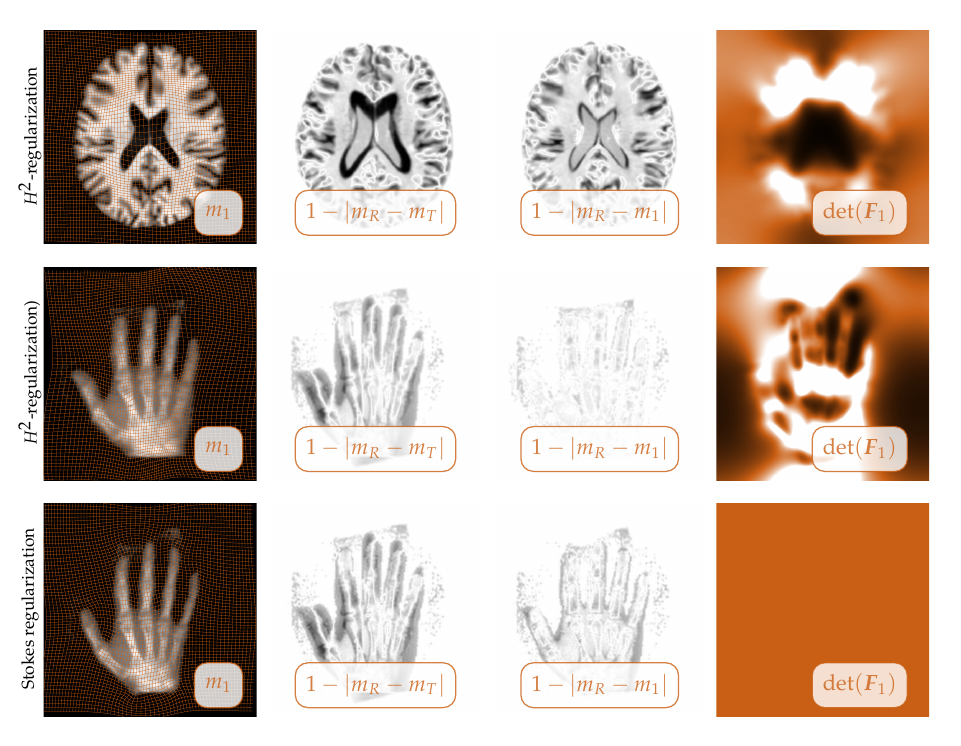}
\caption
{
Qualitative illustration of exemplary registration results of the results for the parameter continuation in $\beta$ reported in \tabref{t:para-cont}. We refer to \tabref{t:para-cont} and the text for details on the experimental setup. We report results for the \iquote{brain images} (top row; plain $H^2$-regularization; $\gamma=0$) and the \iquote{hand images} (middle row: plain $H^2$-regularization ($\gamma=0$); bottom row: Stokes regularization scheme ($\gamma=1$; $H^1$-regularization)). We display the deformed template image $m_1$, a map of the absolute difference between $m_R$ and $m_T$ and between $m_R$ and $m_1$ and a map of the determinant of the deformation gradient $\det(\mat{F}_1)$ (from left to right as indicated in the inset in the images). Information on how to interpret these images can be found in \secref{s:map-and-def-grad}.
}
\label{f:para-cont-im}
\end{figure}

\begin{figure}
\centering
\includegraphics[width=\textwidth]
{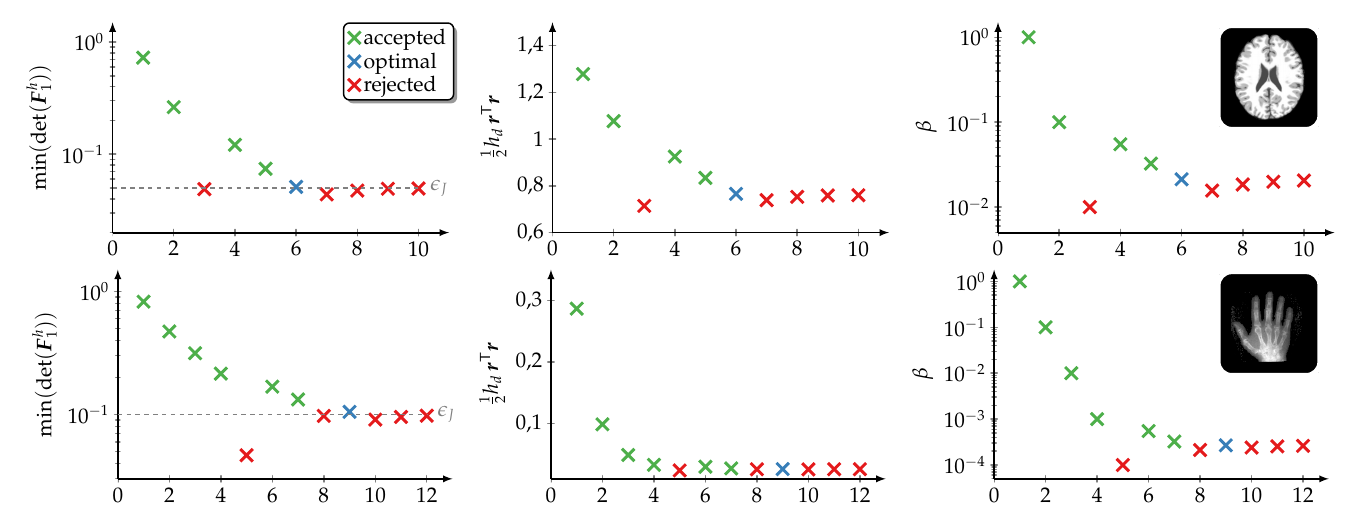}
\caption
{
Exemplary illustration of the course of the parameter continuation in $\beta$ for the quantitative results reported in \tabref{t:para-cont}. We refer to \tabref{t:para-cont} and the text for details on the experimental setup. We report results for the \iquote{brain images} (top row) and the \iquote{hand images} (bottom row) for plain $H^2$-regularization. For each experiment, we display (from left to right) ($i$) the trend of the minimal value of the determinant of the deformation gradient (the dashed line indicates the user-defined lower bound on $\det(\mat{F}_1^h)$), ($ii$) the trend of the $L^2$-distance ($h_d$ is the grid cell volume and $\vect{r} \defeq m_R^h - m_T^h$, $\vect{r} \in \ns{R}^{\tilde{n}}$, $\tilde{n} = \prod_{i=1}^2 n_x^i$), and ($iii$) the trend of $\beta$, all with respect to the parameter continuation step. We indicate our judgment on the results in color. That is, if a result is accepted (i.e.\ $\min(\det(\mat{F}_1^h)) \geq \epsilon_F$) we plot the marker in green and if a result is rejected (i.e.\ $\min(\det(\mat{F}_1^h)) < \epsilon_F$) we plot the marker in red. The optimal value is plotted in blue. The plots correspond to the results reported in \tabref{t:para-cont}.
}
\label{f:para-cont-plots}
\end{figure}

\ipoint{Results} We report the obtained estimates for $\beta$ as well as results for the reconstruction quality and deformation regularity in \tabref{t:para-cont}. We provide an exemplary illustration of the obtained registration results in \figref{f:para-cont-im}. We report results for the course of the parameter continuation in \figref{f:para-cont-plots}.

\begin{table}
\caption
{
Quantitative analysis of the parameter continuation in $\beta$. We report results for the \iquote{hand images} and the \iquote{brain images} for different regularization schemes (see \figref{f:regprob}). The spatial grid size for the images is $\vect{n}_x = (512,512)^\T$. The number of the time points $n_t$ is chosen adaptively (see \secref{s:parameters} for details). We use the full set of stopping conditions (see \secref{s:stopping-criteria}) with a tolerance of $\tau_{\F{J}} = \num{1E-3}$. We report results for plain $H^1$- and $H^2$-regularization ($\gamma=0$; top block) as well as for the Stokes regularization scheme ($\gamma=1$; $H^1$-regularization; bottom block). We invert for a stationary velocity field (i.e.\ $n_c=1$). We report ($i$) the considered lower bound on the deformation gradient ($\epsilon_F$) or the grid angle ($\epsilon_{\theta}$), ($ii$) the number of the required estimation steps, ($iii$) the minimal value for the deformation gradient for the optimal regularization parameter, ($iv$) the computed optimal value for $\beta$ ($\beta^\star$), ($v$) the minimal change in $\beta$ ($\delta\beta_{\min}$), as well as ($vi$) the relative change in the $L^2$-distance ($\|m_R^h - m_1^h\|_{2,\text{rel}}$).
}
\tabadjust
\begin{tabular}{|c|ccccccc|}
\hline
\multicolumn{8}{|c|}{\cellcolor[gray]{0.8}Smoothness regularization}
\\
\hline
  data
& $\F{S}$
& $\epsilon_F$
& steps
& $\min(\det(\mat{F}_1^h))$
& $\beta^\star$
& $\delta\beta_{\min}$
& $\|m_R^h - m_1^h\|_{2,\text{rel}}$
\\
\hline
\multirow{2}{*}{\iquote{brain images}}
& $H^1$
& $\num{5E-2}$
& 9
& $\num{5.090636e-02}$
& $\num{3.109375e-01}$
& $\num{5.000000e-03}$
& $\num{7.014228e-01}$
\\
\cline{2-8}
& $H^2$
& $\num{5E-2}$
& 10
& $\num{5.113449e-02}$
& $\num{2.125000e-02}$
& $\num{5.000000e-04}$
& $\num{5.521585e-01}$
\\
\hline
\multirow{2}{*}{\iquote{hand images}}
& $H^1$
& $\num{1E-1}$
& 10
& $\num{1.125474e-01}$
& $\num{3.390625e-02}$
& $\num{5.000000e-04}$
& $\num{8.737472e-02}$
\\
\cline{2-8}
& $H^2$
& $\num{1E-1}$
& 12
& $\num{1.051038e-01}$
& $\num{2.687500e-04}$
& $\num{5.000000e-06}$
& $\num{6.833687e-02}$
\\
\hline
\hline
\multicolumn{8}{|c|}{\cellcolor[gray]{0.8}Stokes regularization}
\\
\hline
  data
& $\F{S}$
& $\epsilon_{\theta}$
& steps
& $\min(\det(\mat{F}_1^h))$
& $\beta^\star$
& $\delta\beta_{\min}$
& $\|m_R^h - m_1^h\|_{2,\text{rel}}$
\\
\hline
\iquote{hand images}
& $H^1$
& $\pi/16$
& 10
& $\num{9.999987e-01}$
& $\num{2.125000e-02}$
& $\num{5.000000e-04}$
& $\num{1.171091e-01}$
\\
\hline
\end{tabular}
\label{t:para-cont}
\end{table}

\ipoint{Observations} The most important observation is that we can precisely control the properties of our mapping without having to manually tune any parameters. We only have to decide on geometric bounds (the smallest tolerable deformation of a grid element or a bound on the shear angle of the grid cell) the decision on which is intuitive for practitioners.

The accuracy of our method (in space) is only limited by the grid resolution (i.e.\ how much frequencies we can resolve; this statement is confirmed by the experiments conducted in \secref{s:res-convergence}) as well as the defined bounds on the binary search used to estimate $\beta$ (see~\secref{s:parameters}). Clearly, the desired level of accuracy competes with the computational work load we are willing to invest.

For plain $H^1$- and $H^2$-regularization, we achieve an excellent agreement between $m_R$ and $m_1$ (see~\figref{f:para-cont-im}) with a reduction of the $L^2$-distance by approximately half an order of magnitude for the \iquote{brain images} and 1.5 orders of magnitude for the \iquote{hand images} (see \tabref{t:para-cont}). The discrepancy between the lower bound $\epsilon_F$ and $\min(\det(\mat{F}_h^1))$ for the obtained optimal value of $\beta$ is small. In particular, we are e.g.\ bounded from above by an absolute difference of $\num{1.1345E-3}$ for the \iquote{brain images} ($H^2$-regularization) and $\num{5.104E-3}$ for the \iquote{hand images} ($H^2$-regularization). These values are well above the attainable accuracy reported in \secref{s:res-convergence}.

For the results reported for the Stokes regularization scheme we can qualitatively (see~\figref{f:para-cont-im}) and quantitatively (see~\tabref{t:para-cont}) observe that enforcing incompressibility up to numerical accuracy is a too strong prior for the considered problem. However, the key observation and intention of this experiment is to demonstrate that we attain a deformation that is very well behaved (with $\det(\mat{F}_1^h)=1$). A direct comparison to the result obtained for the incompressible case reveals that the mapping is diffeomorphic but displays a large variation in the magnitude of the determinant of the deformation gradient (see the leftmost image in the middle row and bottom row as well as the corresponding maps for $\det(\mat{F}_1^h)$ in \figref{f:para-cont-im}). If we further decrease the bound on $\det(\mat{F}_1^h)$ we will loose control and generate a mapping that locally is close to being nondiffeomorphic. We again emphasize that the intention of this work is the study of algorithmic properties. We will address the practical benefit of exploiting a model of (near-)incompressible flow in a follow-up paper and refer to~\cite{Borzi:2002a, Chen:2011a, Mansi:2011a, Ruhnau:2007a} for potential applications. This exemplary result on real-world data demonstrates that it might be beneficial to consider a relaxation of the incompressibility constraint in order to improve the mismatch between the considered images while maintaining as much control on the deformation regularity as possible.

In future work, we will focus on improvements of the computational efficiency for estimating $\beta$. We have tested combining it with a grid continuation but could not observe strong improvements. We will also investigate the idea of providing a coarse estimation of $\beta$ via the spectral properties of the Hessian and from there on do a parameter continuation (see \secref{s:parameters} for additional comments).

\ipoint{Conclusion} We conclude that the designed framework is highly accurate, is stable, guarantees deformation regularity (assuming that the user-defined tolerance is sufficiently bounded away from irregularities) and does not require any additional tuning of parameters. The user merely has to provide a lower bound on an acceptable volume change or a bound on the distortion of a volume element (shear angle), the decision on which is intuitive for practitioners.

\section{Conclusions}
\label{s:conclusions}

We have presented numerical methods for large deformation diffeomorphic nonrigid image registration that\bipa\item operate in a black-box fashion, \item introduce novel algorithmic features (including a second-order Newton-Krylov method (a full Newton method and a Gauss-Newton approximation), spectral preconditioning, an efficient solver for Stokes problems, and a spectral Galerkin method in time), \item is stable and efficient and \item guarantees deformation regularity with an explicit control on the quality of the deformation.\eipa

In addition, we have conducted a detailed numerical study to demonstrate computational performance and numerical behavior on synthetic and real-world problems. The most important observations of our study are the following:
\begin{itemize}
\item The Newton-Krylov methods outperform the globalized Picard method (see \secref{s:res-convergence}) as we increase the image size and the registration fidelity.
\item We can enforce incompressibility with high accuracy. The numerical accuracy (in space) is only limited by the resolution of the data (see~\secref{s:res-convergence}).
\item We can compute deformations that are guaranteed to be regular (i.e.\ a local diffeomorphism) up to user specifications. Controlling the magnitude of $\det(\mat{F}_1)$ is not sufficient, as volume elements might still collapse (deformation field with strong shear). Therefore, we introduced a parameter continuation in $\beta$ that can be interfaced not only with lower bounds on the determinant of the deformation gradient but also with bounds on the geometric properties of the grid cells (in particular, the shear angle of a grid cell; see~\secref{s:res-parameter-continuation}).
\item The experiments reported in this study demonstrate that it is adequate to limit the inversion to stationary velocity fields when considering two-image registration problems. This observation is in accordance with results reported for other classes of large deformation image registration algorithms~\cite{Arsingny:2006a, Ashburner:2011a, Hernandez:2009a, Mansi:2011a, Vercauteren:2008a, Vercauteren:2009a}. We have additionally provided advice on how to decide on the number of unknowns in time (see~\secref{s:res-expansion}).
\end{itemize}

The control equation for the velocity is a space-time nonlinear elliptic system.  But the main cost in our formulation is the solution of transport problems to compute the image transformation and the adjoint variables. That is, we require two hyperbolic PDE solves for computing the gradient (which essentially corresponds to one Picard iteration or one outer iteration of a Newton-Krylov method; see~\algref{a:outer-iteration}) and an additional two hyperbolic PDE solves for evaluating the incremental control equation (Hessian matrix-vector product in Newton-Krylov methods; see~\algref{a:hessian-matvec}) in each inner iteration of the Krylov subspace method. Because we use a pseudospectral discretization in space, the elliptic solve for the Picard iteration is (for quadratic regularization models) only at the cost of a spectral diagonal scaling. For the Newton-Krylov methods, we have to solve a linear system using an iterative solver. The Picard scheme has a lower cost per iteration but requires more iterations than the Newton-Krylov scheme.

Our results demonstrate that there is a significant difference in stability, computational work and accuracy between the Picard and Newton-Krylov methods, especially when we require a high accuracy of the inversion. If we require an inaccurate solution or use a strong regularization, the differences between Picard and Newton-Krylov methods are less pronounced. Better preconditioning of the Hessian would make the Newton-Krylov approach preferable across the spectrum of accuracy requirements.  The Newton-Krylov approach is not significantly more complex, since we essentially use the same numerical tools that have been used for the solution of the first-order optimality conditions. Also, the individual building blocks ((incremental) forcing term, regularization operator $\D{A}$ and the projection operator $\D{K}$) that appear in the first- and second-order optimality system are very similar. Therefore, the difference of solving the first- or the second-order optimality conditions essentially amounts to interfacing a Krylov-subspace method to solve the saddle point problem.

By formulating the nonrigid image registration as a problem of optimal control, we target the design of a generic, biophysically constrained framework for large deformation diffeomorphic image registration. Further, there are many applications that do require incompressible or near-incompressible deformations, for example, in medical image analysis. Our framework provides such a technology.

We report results only in two dimensions. Nothing in our formulation and numerical approximation is specific to the 2D case. The next steps will be its extension to 3D, and to problems that have time sequences of images. For such cases, we expect to have to invert for a nonstationary velocity field. In addition, we aim at designing a framework that allows for a relaxation of the incompressibility constraint, as we observed in this study that incompressibility might be a too strong prior. We also observed (results not included in this study) that the use of an incompressibility constraint can promote shear. In a follow-up paper, we will target this problem by introducing a novel continuum mechanical model that allows us to control the shear inside the deformation field.

\section*{Acknowledgments}

We would like to thank Amir Gholaminejad, Georg Stadler and Bryan Quaife for helpful discussions and comments. We would like to thank Florian Tramnitzke for his initiative work on this project. This material is based upon work supported by AFOSR grants FA9550-12-10484 and FA9550-11-10339; by NSF grants CCF-1337393, OCI-1029022; by the U.S. Department of Energy, Office of Science, Office of Advanced Scientific Computing Research, Applied Mathematics program under Award Numbers DE-SC0010518, DE-SC0009286; and DE-FG02-08ER2585 and by NIH grant 10042242. Any opinions, findings, and conclusions or recommendations expressed herein are those of the authors and do not necessarily reflect the views of the AFOSR or the NSF.

\begin{appendix}

\section{Expansion in Time: Derivation}
\label{s:appendix-expansion}

This section summarizes modifications of the regularization operator as well as the (incremental) control equation on account of the expansion in time (see \secref{s:expansion}). Inserting~\eqref{e:expansion} into~\eqref{e:quadratic-regularization} yields
\begin{equation}\label{e:regularization-expansion}
\iut{\F{S}[\vect{v}]} = \half{\beta}\sum_{l=1}^{n_c}\sum_{l'=1}^{n_c} c_{ll'}
\iom{ \ipl{\D{B}[\vect{v}_l]}{\D{B}[\vect{v}_{l'}]} },
\quad \text{where}\,\,
c_{ll'}\defeq\displaystyle\iut{ b_l\,b_{l'} }\;\;\forall l,l'\in\mathbf{N}.
\end{equation}

Taking first and second variations with respect to the $l$-th expansion coefficient $\vect{v}_l$ yields the control equation
\begin{equation}\label{e:control-equation-expansion}
\beta\sum_{l'=1}^{n_c} c_{ll'} \D{A}[\vect{v}_{l'}]
+\iut{b_l\,\D{K}[\vect{f}]} = 0, \quad l = 1,\ldots,n_c,
\end{equation}

\noindent and the incremental control equation
\begin{equation}\label{e:incremental-control-equation-expansion}
\D{H}_l[\vect{\tilde{v}}_l] \defeq
\beta\sum_{l'=1}^{n_c} c_{ll'} \D{A}[\vect{\tilde{v}}_{l'}]
+\iut {b_l\,\D{K}[\vect{\tilde{f}}]}, \quad l = 1,\ldots,n_c,
\end{equation}

\noindent respectively. Accordingly, the operators $\D{B}$ and $\D{A}$ simply act on $\vect{v}_l$ instead of $\vect{v}$. The definition for these operators can be found in~\secref{s:optimality-conditions}.

We use a global basis on the unit time horizon for the expansion (see~\secref{s:expansion}). We use \emph{Chebyshev polynomials} as basis functions in~\eqref{e:expansion} on account of their excellent approximation properties as well as their orthogonality. The latter property considerably reduces the computational complexity, since $c_{ll'} = 0$ for all $l,l'$, $l'\not= l$, and $c_{ll} = 1$ (see \eqref{e:regularization-expansion}, \eqref{e:control-equation-expansion} and~\eqref{e:incremental-control-equation-expansion}). To avoid Runge's phenomenon (see e.\,g.\ \cite[p.~82ff.]{Boyd:2000a}) \iname{Chebyshev-Gauss-Labatto nodes} are used.

\section{Incompressibility Constraint: Elimination}
\label{s:elimination-details}

Here, we derive the elimination of $p$ and $\tilde{p}$ from the optimality systems for the Stokes regularization scheme. We only consider a quadratic $H^1$-regularization for the velocity $\vect{v}$. However, the same line of arguments applies to the $H^2$-regularization model.

Applying the divergence to~\eqref{e:first-order-control} results in\footnote{Since we discuss the implementation of the incompressibility constraint we set $\gamma=1$.}
\[
-\idiv \beta \ilap \vect{v} + \ilap p + \idiv \vect{f} = 0
\quad\text{in}\quad\Omega \times [0,1].
\]

\noindent Under the optimality assumption $\idiv \vect{v} = 0$ it follows from the definition of the vectorial Laplacian that $p = -\ilap^{-1} \left(\idiv \vect{f}\right)$. Inserting this expression into~\eqref{e:first-order-control} projects $\vect{v}$ onto the manifold of divergence-free velocity fields and as such eliminates~\eqref{e:first-order-incomp} (assuming that the initial $\vect{v}$ is divergence-free). Accordingly, we obtain the control equation (reduced gradient)
\begin{equation}\label{e:first-order-control-elim-kf}
-\beta\ilap\vect{v} + \D{K}[\vect{f}] = 0,
\quad\text{where}\;\; \D{K}[\vect{f}] \defeq
-\igrad(\ilap^{-1}(\idiv \vect{f})) + \vect{f},
\end{equation}

\noindent to replace~\eqref{e:first-order-control}. This expression is equivalent to~\eqref{e:eliminated-control} in the case when an $H^1$-regularization model is used (i.e.\ $\D{A} = -\ilap$). As stated above, the derivation also holds for the $H^2$-regularization operator. We only have to replace $-\beta\ilap\vect{v}$ by $\beta\ilap^2\vect{v}$. Computing the second variations of the weak form of the eliminated system yields the incremental control equation
\[
-\beta\D{A}[\vect{\tilde{v}}] + \D{K}[\vect{\tilde{f}}] = -\vect{\hat{g}},
\]

\noindent where $\hat{\vect{g}}$ is the reduced gradient in \eqref{e:first-order-control-elim-kf}.

\section{Relation to LDDMM}
\label{s:relation-to-lddmm}

In this section we relate our work to~\cite{Hart:2009a} and by that to approaches based on LDDMM~\cite{Ashburner:2011a, Beg:2005a, Dupuis:1998a, Miller:2004a, Trouve:1998a}. Since the work in~\cite{Beg:2005a, Hart:2009a} is based on first-order information, we only consider the reduced gradient in~\eqref{e:first-order-control} (setting $\gamma=0$). In weak form we have
\begin{align*}
\iut{\langle \vect{g}, \tilde{v}\rangle_{L^2(\Omega)}}
& =
\iut
{
\langle\beta(\D{B}\D{B}^\adj)[\vect{v}] + \vect{f},
\vect{\tilde{v}}\rangle_{L^2(\Omega)}
}
=
\iut{
\langle\beta\vect{v}+(\D{B}\D{B}^\adj)^{-1}[\vect{f}],
\vect{\tilde{v}}\rangle_{\fs{W}}.
}
\end{align*}

\noindent The expression
$
\beta\vect{v}+(\D{B}\D{B}^\adj)^{-1}[\vect{f}]
= \vect{v} + (\beta\D{B}\D{B}^\adj)^{-1}[\vect{f}]
$
is exactly the gradient in the function space $\fs{W}$ that has been used in~\cite{Hart:2009a}. This expression yields the \emph{preconditioned gradient descent scheme}
\[
\vect{v}_{k+1}^h = \vect{v}_k^h
-
\alpha_k((\beta\D{B}^h\D{B}^{h,\adj})^{-1}
[
\vect{f}^h_k
]
- \vect{v}_k^h),
\]

\noindent where $(\beta\D{B}^h\D{B}^{h,\adj})^{-1}[\vect{f}^h]$ is nothing but a Picard iterate (see~\eqref{e:picard}). Subtracting $\vect{v}_k^h$ translates this iterate into an update. This is exactly the formulation we have used in this work (see \secref{s:picard}) so that the considered first-order method is equivalent to the solver used in~\cite{Hart:2009a} (under the assumption that $\alpha_k=1$, i.e.\ if we neglect the line search). Accordingly, the same line of arguments used in~\cite{Hart:2009a} to relate their work to LDDMM apply to our numerical framework.

\section{Measures of Deformation Regularity}
\label{s:map-and-def-grad}

\subsection{Deformation Map}

To visualize the deformation pattern, $\vect{y}$ has to be inferred from $\vect{v}$. This can be done by solving
\begin{equation}\label{e:displacement}
\p_t \vect{u} + (\igrad \vect{u}) \vect{v} = \vect{v}
\;\;{\rm in}\;\; \Omega \times (0,1],
\qquad \vect{u} = 0 \;\;{\rm in}\;\;\Omega \times \{0\},
\end{equation}

\noindent with periodic boundary conditions on $\p\Omega$. Here, $\vect{u} : \bar{\Omega} \times [0,1] \rightarrow \ns{R}^d$, $(\vect{x},t)\mapsto \vect{u}(\vect{x},t)$, is a displacement field and $\vect{y} \defeq \vect{x} - \vect{u}_1$, $\vect{y} : \bar{\Omega} \rightarrow \ns{R}^d$, where $\vect{u}_1 \defeq \vect{u}(\cdot, t=1)$, $\vect{u}_1: \bar{\Omega} \rightarrow \ns{R}^{d}$, $\vect{x}\mapsto \vect{u}_1(\vect{x})$.

\ipoint{Visualization} As can be seen in the visualization of the deformed grids, the mapping $\vect{y}$ actually corresponds the inverse of the deformation map applied to an image. This reflects the fact that our model is formulated in an Eulerian frame of reference. Note that all images reported are high-resolution vector graphics. Zooming in in the digital version of the paper will reveal local properties of the deformation map.

\subsection{Deformation Gradient}

It is well known from calculus that the determinant of the Jacobian matrix $\det(\igrad \vect{y})$ can be used to assess invertibility of $\vect{y}$ as well as local volume change, provided that $\vect{y}\in C^2(\Omega)^d$. In the framework of continuum mechanics, we can obtain this information from the deformation tensor field $\mat{F} : \bar{\Omega} \times [0,1] \rightarrow \ns{R}^{d\times d}$, where $\mat{F}$ is related to $\vect{v}$ by
\begin{equation}\label{e:def-grad}
\partial_t \mat{F} + (\vect{v}\cdot \igrad) \mat{F}
= (\igrad \vect{v}) \mat{F} \;\; {\rm in} \;\; \Omega \times (0,1],
\qquad \mat{F} = \mat{I} \;\; {\rm in} \quad \Omega \times \{0\},
\end{equation}

\noindent with periodic boundary conditions on $\p\Omega$. Here, $\mat{I} = \operatorname{diag}(1,\ldots,1)\in\ns{R}^{d\times d}$; $\det(\mat{F}_1)$ is equivalent to $\det(\igrad \vect{y})$, where $\mat{F}_1 \defeq \mat{F}(\cdot, t=1)$, $\mat{F}_1 : \bar{\Omega} \rightarrow \ns{R}^{d \times d}$, $\vect{x} \mapsto \mat{F}_1(\vect{x})$.

\ipoint{Visualization} We limit the color map for the display of $\det(\mat{F}_1)$ to $[0,2]$. In particular, the color map ranges from black (compression: $\det(\mat{F}_1) \in (0,1)$; black corresponds to values of 0 or below (due to clipping), which represents a singularity or the loss of mass, respectively) to orange (mass conservation: $\det(\mat{F}_1) = 1$) to white (expansion: $\det(\mat{F}_1) > 1$; white represents values of 2 or greater (due to clipping)).

\end{appendix}



\begin{thebibliography}{10}

\bibitem{Arsingny:2006a}
{\sc V.~Arsigny, O.~Commowick, X.~Pennec, and N.~Ayache}, {\em A
  {L}og-{E}uclidean framework for statistics on diffeomorphisms}, in Lect Notes
  Comput Sc, vol.~4190, 2006, pp.~924--931.

\bibitem{Ashburner:2007a}
{\sc J.~Ashburner}, {\em A fast diffeomorphic image registration algorithm},
  NeuroImage, 38 (2007), pp.~95--113.

\bibitem{Ashburner:2011a}
{\sc J.~Ashburner and K.~J. Friston}, {\em Diffeomorphic registration using
  geodesic shooting and {G}auss-{N}ewton optimisation}, NeuroImage, 55 (2011),
  pp.~954--967.

\bibitem{Beg:2005a}
{\sc M.~F. Beg, M.~I. Miller, A.~Trouv\'e, and L.~Younes}, {\em Computing large
  deformation metric mappings via geodesic flows of diffeomorphisms}, Int J
  Comput Vis, 61 (2005), pp.~139--157.

\bibitem{Benzi:2005a}
{\sc M.~Benzi, G.~H. Golub, and J.~Liesen}, {\em Numerical solution of saddle
  point problems}, Acta Numerica, 14 (2005), pp.~1--137.

\bibitem{Benzi:2011a}
{\sc M.~Benzi, E.~Haber, and L.~Taralli}, {\em A preconditioning technique for
  a class of {PDE}-constrained optimization problems}, Adv Comput Math, 35
  (2011), pp.~149--173.

\bibitem{Biros:2008a}
{\sc G.~Biros and G.~Do\v{g}an}, {\em A multilevel algorithm for inverse
  problems with {PDE} constraints}, Inverse Probl, 24 (2008).

\bibitem{Biros:2005a}
{\sc G.~Biros and O.~Ghattas}, {\em Parallel
  {L}agrange-{N}ewton-{K}rylov-{S}chur methods for {PDE}-constrained
  optimization---{P}art {I}: {T}he {K}rylov-{S}chur solver}, SIAM J Sci Comput,
  27 (2005), pp.~687--713.

\bibitem{Biros:2005b}
\leavevmode\vrule height 2pt depth -1.6pt width 23pt, {\em Parallel
  {L}agrange-{N}ewton-{K}rylov-{S}chur methods for {PDE}-constrained
  optimization---{P}art {II}: {T}he {L}agrange-{N}ewton solver and its
  application to optimal control of steady viscous flows}, SIAM J Sci Comput,
  27 (2005), pp.~714--739.

\bibitem{Borzi:2002a}
{\sc A.~Borz\`i, K.~Ito, and K.~Kunisch}, {\em Optimal control formulation for
  determining optical flow}, SIAM J Sci Comput, 24 (2002), pp.~818--847.

\bibitem{Boyd:2000a}
{\sc J.~P. Boyd}, {\em Chebyshev and {F}ourier spectral methods}, Dover,
  Mineola, New York, US, 2000.

\bibitem{Brezzi:1991a}
{\sc F.~Brezzi and M.~Fortin}, eds., {\em Mixed and hybrid finite element
  methods}, Springer, 1991.

\bibitem{Broit:1981a}
{\sc C.~Broit}, {\em Optimal registration of deformed images}, PhD thesis,
  Computer and Information Science, University of Pennsylvania, Philadelphia,
  Pennsylvania, US, 1981.

\bibitem{Burger:2013a}
{\sc M.~Burger, J.~Modersitzki, and L.~Ruthotto}, {\em A hyperelastic
  regularization energy for image registration}, SIAM J Sci Comput, 35 (2013),
  pp.~B132--B148.

\bibitem{Byrd:2008a}
{\sc R.~H. Byrd, F.~E. Curtis, and J.~Nocedal}, {\em An inexact {SQP} method
  for equality constrained optimization}, SIAM J Optim, 19 (2008),
  pp.~351--369.

\bibitem{Chen:2011a}
{\sc K.~Chen and D.~A. Lorenz}, {\em Image sequence interpolation using optimal
  control}, J Math Imag Vis, 41 (2011), pp.~222--238.

\bibitem{Chen:2012a}
\leavevmode\vrule height 2pt depth -1.6pt width 23pt, {\em Image sequence
  interpolation based on optical flow, segmentation and optimal control}, IEEE
  T Image Process, 21 (2012), pp.~1020--1030.

\bibitem{Christensen:1994a}
{\sc G.~E. Christensen, R.~D. Rabbitt, and M.~I. Miller}, {\em {3D} brain
  mapping using a deformable neuroanatomy}, Phys Med Biol, 39 (1994),
  pp.~609--618.

\bibitem{Christensen:1996a}
\leavevmode\vrule height 2pt depth -1.6pt width 23pt, {\em Deformable templates
  using large deformation kinematics}, IEEE T Image Process, 5 (1996),
  pp.~1435--1447.

\bibitem{Dembo:1982a}
{\sc R.~S. Dembo, S.~C. Eisenstat, and T.~Steihaug}, {\em Inexact {N}ewton
  methods}, SIAM J Numer Anal, 19 (1982), pp.~400--408.

\bibitem{Dembo:1983a}
{\sc R.~S. Dembo and T.~Steihaug}, {\em Truncated-{N}ewton algorithms for
  large-scale unconstrained optimization}, Math Program, 26 (1983),
  pp.~190--212.

\bibitem{Droske:2003a}
{\sc M.~Droske and M.~Rumpf}, {\em A variational approach to non-rigid
  morphological registration}, SIAM Appl Math, 64 (2003), pp.~668--687.

\bibitem{Dupuis:1998a}
{\sc P.~Dupuis, U.~Gernander, and M.~I. Miller}, {\em Variational problems on
  flows of diffeomorphisms for image matching}, Qart Appl Math, 56 (1998),
  pp.~587--600.

\bibitem{Eisenstat:1996a}
{\sc S.~C. Eisentat and H.~F. Walker}, {\em Choosing the forcing terms in an
  inexact {N}ewton method}, SIAM J Sci Comput, 17 (1996), pp.~16--32.

\bibitem{Fischer:2002a}
{\sc B.~Fischer and J.~Modersitzki}, {\em Fast diffusion registration}, Contemp
  Math, 313 (2002), pp.~117--129.

\bibitem{Fischer:2003a}
\leavevmode\vrule height 2pt depth -1.6pt width 23pt, {\em Curvature based
  image registration}, J Math Imag Vis, 18 (2003), pp.~81--85.

\bibitem{Fischer:2008a}
\leavevmode\vrule height 2pt depth -1.6pt width 23pt, {\em Ill-posed medicine
  -- an introduction to image registration}, Inverse Probl, 24 (2008),
  pp.~1--16.

\bibitem{FrohnSchauf:2008a}
{\sc C.~Frohn-Schauf, S.~Henn, and K.~Witsch}, {\em Multigrid based total
  variation image registration}, Comut Visual Sci, 11 (2008), pp.~101--113.

\bibitem{Gill:1981a}
{\sc P.~E. Gill, W.~Murray, and M.~H. Wright}, {\em Practical optimization},
  Academic Press, Waltham, Massachusetts, US, 1981.

\bibitem{Golub:1996a}
{\sc G.~H. Golub and C.~F.~V. Loan}, {\em Matrix Computations}, The Johns
  Hopkins University Press, 3~ed., 1996.

\bibitem{Gooya:2012a}
{\sc A.~Gooya, K.~M. Pohl, M.~Bilello, L.~Cirillo, G.~Biros, E.~R. Melhem, and
  C.~Davatzikos}, {\em {GLISTR}: {G}lioma image segmentation and registration},
  IEEE T Med Imaging, 31 (2013), pp.~1941--1954.

\bibitem{Gunzburger:2003a}
{\sc M.~D. Gunzburger}, {\em Perspectives in flow control and optimization},
  SIAM, Philadelphia, Pennsylvania, US, 2003.

\bibitem{Gurtin:1981a}
{\sc M.~E. Gurtin}, {\em An introduction to continuum mechanics}, vol.~158 of
  Mathematics in Science and Engineering, Academic Press, 1981.

\bibitem{Haber:2001a}
{\sc E.~Haber and U.~M. Ascher}, {\em Preconditioned all-at-once methods for
  large, sparse parameter estimation problems}, Inverse Probl, 17 (2001),
  pp.~1847--1864.

\bibitem{Haber:2000a}
{\sc E.~Haber, U.~M. Ascher, and D.~Oldenburg}, {\em On optimization techniques
  for solving nonlinear inverse problems}, Inverse Probl, 16 (2000),
  pp.~1263--1280.

\bibitem{Haber:2010a}
{\sc E.~Haber, R.~Horesh, and J.~Modersitzki}, {\em Numerical optimization for
  constrained image registration}, Numer Linear Algebr, 17 (2010),
  pp.~343--359.

\bibitem{Haber:2004a}
{\sc E.~Haber and J.~Modersitzki}, {\em Numerical methods for volume preserving
  image registration}, Inverse Probl, 20 (2004), pp.~1621--1638.

\bibitem{Haber:2006a}
\leavevmode\vrule height 2pt depth -1.6pt width 23pt, {\em A multilevel method
  for image registration}, SIAM J Sci Comput, 27 (2006), pp.~1594--1607.

\bibitem{Haber:2007a}
\leavevmode\vrule height 2pt depth -1.6pt width 23pt, {\em Image registration
  with guaranteed displacement regularity}, Int J Comput Vis, 71 (2007),
  pp.~361--372.

\bibitem{Han:2014a}
{\sc L.~Hand, J.~H. Hipwell, B.~Eiben, D.~Barratt, M.~Modat, S.~Ourselin, and
  D.~J. Hawkes}, {\em A nonlinear biomechanical model based registration method
  for aligning prone and supine {MR} breast images}, IEEE T Med Imaging, 33
  (2014), pp.~682--694.

\bibitem{Hart:2009a}
{\sc G.~L. Hart, C.~Zach, and M.~Niethammer}, {\em An optimal control approach
  for deformable registration}, in Proc CVPR IEEE, 2009, pp.~9--16.

\bibitem{Henn:2005a}
{\sc S.~Henn}, {\em A multigrid method for a fourth-order diffusion equation
  with application to image processing}, SIAM J Sci Comput, 27 (2005),
  pp.~831--849.

\bibitem{Hernandez:2014a}
{\sc H.~Hernandez}, {\em Gauss-{N}ewton inspired preconditioned optimization in
  large deformation diffeomorphic metric mapping}, Phys Med Biol, 59 (2014),
  pp.~6085--6115.

\bibitem{Hernandez:2009a}
{\sc M.~Hernandez, M.~N. Bossa, and S.~Olmos}, {\em Registration of anatomical
  images using paths of diffeomorphisms parameterized with stationary vector
  field flows}, Int J Comput Vis, 85 (2009), pp.~291--306.

\bibitem{Hogea:2008a}
{\sc C.~Hogea, C.~Davatzikos, and G.~Biros}, {\em Brain-tumor interaction
  biophysical models for medical image registration}, SIAM J Sci Comput, 30
  (2008), pp.~3050--3072.

\bibitem{Horn:1981a}
{\sc B.~K.~P. Horn and B.~G. Shunck}, {\em Determining optical flow}, Artif
  Intell, 17 (1981), pp.~185--203.

\bibitem{Kalmoun:2011a}
{\sc E.~M. Kalmoun, L.~Garrido, and V.~Caselles}, {\em Line search multilevel
  optimization as computational methods for dense optical flow}, SIAM J Imag
  Sci, 4 (2011), pp.~695--722.

\bibitem{Lee:2010a}
{\sc E.~Lee and M.~Gunzburger}, {\em An optimal control formulation of an image
  registration problem}, J Math Imag Vis, 36 (2010), pp.~69--80.

\bibitem{Lee:2011a}
\leavevmode\vrule height 2pt depth -1.6pt width 23pt, {\em Analysis of finite
  element discretization of an optimal control fomualtion of the image
  registration problem}, SIAM J Numer Anal, 49 (2011), pp.~1321--1349.

\bibitem{Loeckx:2004a}
{\sc D.~Loeckx, F.~Maes, D.~Vandermeulen, and P.~Suetens}, {\em Nonrigid image
  registration using free-form deformations with a local rigidity constraint},
  in Lect Notes Comput Sc, vol.~3216, 2004, pp.~639--646.

\bibitem{Mang:2012b}
{\sc A.~Mang, A.~Toma, T.~A. Schuetz, S.~Becker, T.~Eckey, C.~Mohr,
  D.~Petersen, and T.~M. Buzug}, {\em Biophysical modeling of brain tumor
  progression: from unconditionally stable explicit time integration to an
  inverse problem with parabolic {PDE} constraints for model calibration}, Med
  Phys, 39 (2012), pp.~4444--4459.

\bibitem{Mansi:2011a}
{\sc T.~Mansi, X.~Pennec, M.~Sermesant, H.~Delingette, and N.~Ayache}, {\em
  {iLogDemons}: {A} demons-based registration algorithm for tracking
  incompressible elastic biological tissues}, Int J Comput Vis, 92 (2011),
  pp.~92--111.

\bibitem{Miller:2004a}
{\sc M.~I. Miller}, {\em Computational anatomy: {S}hape, growth and atrophy
  comparison via diffeomorphisms}, NeuroImage, 23 (2004), pp.~S19--S33.

\bibitem{Modersitzki:2004a}
{\sc J.~Modersitzki}, {\em Numerical methods for image registration}, Oxford
  University Press, New York, 2004.

\bibitem{Modersitzki:2008a}
\leavevmode\vrule height 2pt depth -1.6pt width 23pt, {\em {FLIRT} with
  rigidity---image registration with a local non-rigidity penalty}, Int J
  Comput Vis, 76 (2008), pp.~153--163.

\bibitem{Modersitzki:2009a}
\leavevmode\vrule height 2pt depth -1.6pt width 23pt, {\em {FAIR}: Flexible
  algorithms for image registration}, SIAM, Philadelphia, Pennsylvania, US,
  2009.

\bibitem{Museyko:2009a}
{\sc O.~Museyko, M.~Stiglmayr, K.~Klamroth, and G.~Leugering}, {\em On the
  application of the {M}onge-{K}antorovich problem to image registration}, SIAM
  J Imaging Sci, 2 (2009), pp.~1068--1097.

\bibitem{Nielsen:2002a}
{\sc M.~Nielsen, P.~Johansen, A.~D. Jackson, and B.~Lautrup}, {\em Brownian
  warps: {A} least committed prior for non-rigid registration}, in Lect Notes
  Comput Sc, vol.~2489, 2002, pp.~557--564.

\bibitem{Nocedal:2006a}
{\sc J.~Nocedal and S.~J. Wright}, {\em Numerical Optimization}, Springer, New
  York, New York, US, 2006.

\bibitem{Pennec:2005a}
{\sc X.~Pennec, R.~Stefanescu, V.~Arsigny, P.~Fillard, and N.~Ayache}, {\em
  Riemannian elasticity: {A} statistical regularization framework for
  non-linear registration}, in Lect Notes Comput Sc, vol.~3750, 2005,
  pp.~943--950.

\bibitem{Rohlfing:2003a}
{\sc T.~Rohlfing, C.~R. Maurer, D.~A. Bluemke, and M.~A. Jacobs}, {\em
  Volume-preserving nonrigid registration of {MR} breast images using free-form
  deformation with an incompressibility constraint}, IEEE T Med Imaging, 22
  (2003), pp.~730--741.

\bibitem{Ruhnau:2007a}
{\sc P.~Ruhnau and C.~Schn\"orr}, {\em Optical {S}tokes flow estimation: {A}n
  imaging-based control approach}, Exp Fluids, 42 (2007), pp.~61--78.

\bibitem{Rumpf:2009a}
{\sc M.~Rumpf and B.~Wirth}, {\em A nonlinear elastic shape averaging
  approach}, SIAM J Imaging Sci, 2 (2009), pp.~800--833.

\bibitem{Sdika:2008a}
{\sc M.~Sdika}, {\em A fast nonrigid image registration with constraints on the
  {J}acobian using large scale constrained optimization}, IEEE T Med Imaging,
  27 (2008), pp.~271--281.

\bibitem{Sotiras:2013a}
{\sc A.~Sotiras, C.~Davatzikos, and N.~Paragios}, {\em Deformable medical image
  registration: {A} survey}, IEEE T Med Imaging, 32 (2013), pp.~1153--1190.

\bibitem{Sundar:2009a}
{\sc H.~Sundar, C.~Davatzikos, and G.~Biros}, {\em Biomechanically constrained
  {4D} estimation of mycardial motion}, in Lect Notes Comput Sc, vol.~5762,
  2009, pp.~257--265.

\bibitem{Thirion:1998a}
{\sc J.~P. Thirion}, {\em Image matching as a diffusion process: {A}n analogy
  with {M}axwell's demons}, Med Image Anal, 2 (1998), pp.~243--260.

\bibitem{Trouve:1998a}
{\sc A.~Trouv\'e}, {\em Diffeomorphism groups and pattern matching in image
  analysis}, Int J Comput Vis, 28 (1998), pp.~213--221.

\bibitem{Vercauteren:2008a}
{\sc T.~Vercauteren, X.~Pennec, A.~Perchant, and N.~Ayache}, {\em Symmetric
  log-domain diffeomorphic registration: {A} demons-based approach}, in Lect
  Notes Comput Sc, vol.~5241, 2008, pp.~754--761.

\bibitem{Vercauteren:2009a}
\leavevmode\vrule height 2pt depth -1.6pt width 23pt, {\em Diffeomorphic
  demons: {E}fficient non-parametric image registration}, NeuroImage, 45
  (2009), pp.~S61--S72.

\bibitem{Vialard:2012a}
{\sc F.-X. Vialard, L.~Risser, D.~Rueckert, and C.~J. Cotter}, {\em
  Diffeomorphic {3D} image registration via geodesic shooting using an
  efficient adjoint calculation}, Int J Comput Vis, 97 (2012), pp.~229--241.

\bibitem{Vogel:2002a}
{\sc C.~R. Vogel}, {\em Computational methods for inverse problems}, SIAM,
  Philadelphia, Pennsylvania, US, 2002.

\bibitem{Yanovsky:2007a}
{\sc I.~Yanovsky, P.~M. Thompson, S.~Osher, and A.~D. Loew}, {\em Topology
  preserving log-unbiased nonlinear image registration: {T}heory and
  implementation}, in Proc CVPR IEEE, 2007, pp.~1--8.

\bibitem{Younes:2007a}
{\sc L.~Younes}, {\em Jacobi fields in groups of diffeomorphisms and
  applications}, Qart Appl Math, 650 (2007), pp.~113--134.

\end{thebibliography}
\end{document}